\documentclass[a4paper,11pt]{article}

\usepackage[a4paper, margin=2.97cm]{geometry}

\usepackage[T1]{fontenc}
\usepackage[utf8]{inputenc}
\usepackage{lmodern}
\usepackage[english]{babel}
\usepackage{microtype}

\usepackage{standalone}
\usepackage[dvipsnames]{xcolor}
\usepackage{tikz}
\usetikzlibrary{calc,shapes.arrows,decorations.pathreplacing,decorations.markings,arrows,positioning,patterns}
\usepackage{multicol}

\usepackage{graphicx}
\usepackage{caption}
\usepackage{subcaption}
\usepackage{hyperref}

\definecolor{rouge}{RGB}{255,77,77}
\definecolor{vert}{RGB}{0,178,102}
\definecolor{jaune}{RGB}{255,255,0}
\definecolor{violet}{RGB}{208,32,144}
\definecolor{orange}{RGB}{255,140,0}
\definecolor{bleu}{RGB}{0,0,205}

\usepackage{amsmath,amsfonts,amssymb}
\usepackage{amsthm}
\usepackage{dsfont}

\usepackage{cite}
\usepackage[capitalise]{cleveref}

\usepackage{enumitem}
\setenumerate{itemsep=1pt,topsep=3pt}
\setitemize{itemsep=1pt,topsep=3pt}

\graphicspath{img/}

\newcommand{\N}{\ensuremath{\mathbb{N}}}
\newcommand{\Z}{\ensuremath{\mathbb{Z}}}

\newcommand{\R}{\ensuremath{\mathbb{R}}}

\newcommand{\A}{\ensuremath{\mathcal{A}}}
\newcommand{\B}{\ensuremath{\mathcal{B}}}
\newcommand{\F}{\ensuremath{\mathcal{F}}}
\newcommand{\mC}{\ensuremath{\mathcal{C}}}

\newcommand{\orangesquare}{\ensuremath{ {\color{orange} \blacksquare} }}

\newcommand{\leftgen}{\ensuremath{ \leftarrow\!\! }}
\newcommand{\rightgen}{\ensuremath{ \rightarrow\!\! }}
\newcommand{\upgen}{\ensuremath{ \uparrow\!\! }}
\newcommand{\downgen}{\ensuremath{ \downarrow\!\! }}

\newtheorem{theorem}{Theorem}

\newtheorem{definition}[theorem]{Definition}
\newtheorem{example}[theorem]{Example}
\newtheorem{conjecture}[theorem]{Conjecture}
\newtheorem{lemma}[theorem]{Lemma}
\newtheorem{coro}[theorem]{Corollary}

\newtheorem{proposition}[theorem]{Proposition}

\newtheorem*{proposition*}{Proposition}
\newtheorem*{remark*}{Remark}
\newtheorem*{example*}{Example}

\title{The domino problem is undecidable on surface groups}

\author{Nathalie Aubrun, Sebasti\'an Barbieri and Etienne Moutot}

\begin{document}

\maketitle

\begin{abstract}
We show that the domino problem is undecidable on orbit graphs of non-deterministic substitutions which satisfy a technical property. As an application, we prove that the domino problem is undecidable for the fundamental group of any closed orientable surface of genus at least 2.
\end{abstract}

\begin{center}
  \begin{figure}[hb]
  	\centering
    ~\\~\\
  	\includegraphics[scale=1.3]{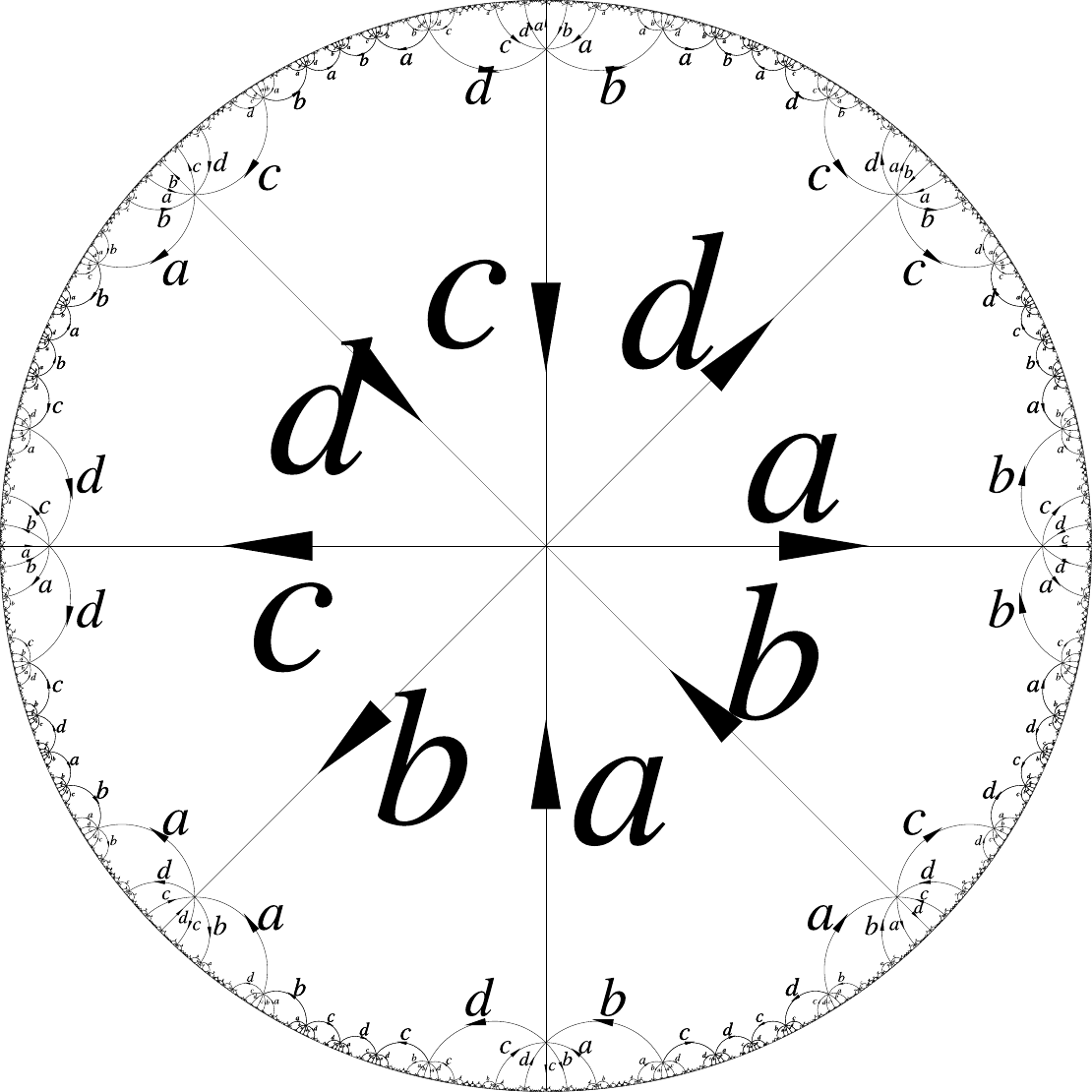}
  \end{figure}
\end{center}

\newpage

 \section{Introduction}

 Initially studied by Wang~\cite{Wang1961}, the domino problem asks whether there exists an algorithm which determines if a finite set of unit square tiles with colored edges --called Wang tiles-- can be used to tile the plane in such a way that edges of adjacent Wang tiles have the same color. It was originally conjectured by Wang that every set of Wang tiles which admits such a tiling of the plane, also admits a periodic tiling. Furthermore, he showed that this property would imply the decidability of the domino problem. However, Berger~\cite{BergerPhD} and later Robinson~\cite{Robinson1971} both combined their constructions of an aperiodic set of Wang tiles and a reduction from the halting problem of Turing machines to show that the domino problem was in fact undecidable.

The domino problem can be naturally extended to a much broader context. Let $\Gamma$ be a labeled directed infinite graph, $\A$ a finite set, and $F = \{p_1,\dots,p_n\}$ a finite list of colorings $p_i$ of vertices of finite connected subgraphs of $\Gamma$ by $\A$. The domino problem of $\Gamma$ asks whether there exists an algorithm which receives on input a finite list $F$ as above, and decides if there exists a coloring of the vertices of $\Gamma$ for which none of the $p_i$ embed as a colored labeled subgraph. Naturally, Wang's domino problem can be reinterpreted in this setting by letting $\Gamma$ be the bi-infinite square grid, $\A$ the set of Wang tiles, and $F$ the list of all horizontal or vertical pairs of tiles whose colors do not match.

A particularly interesting case is when $\Gamma$ is a labeled directed Cayley graph of a finitely generated group $G$ given by a set of its generators. In this case, there is a direct correspondence between colorings of the vertices of $\Gamma$ by $\A$ which avoid a list of forbidden colored subgraphs as described above, and subshifts of finite type (SFT), that is, closed and translation invariant subsets of $\A^G$ which are determined by a finite list of forbidden patterns. What is more, it can be shown (see~\cite{BertheRigo2018}) that the domino problems of all such Cayley graphs of $G$ are computationally many-one equivalent, and thus one may speak about the domino problem of $G$. In the particular case when $G = \Z$, the domino problem is decidable: every $\Z$-SFT can be represented by a labeled finite graph~\cite{LindMarcus1995}, and the existence of a configuration in the SFT (i.e. a bi-infinite word) is equivalent to the existence of a cycle in the graph. The case of $\Z^2$ coincides with the
formalism of Wang tiles, and is thus undecidable.

The domino problem on graphs other than $\Z$ and $\Z^2$ has been largely investigated. The undecidability for a graph which models the hyperbolic plane was settled by Kari~\cite{Kari}, and can also be obtained from the construction of a hierarchical aperiodic tiling on the hyperbolic plane by Goodman-Strauss~\cite{GoodmanStrauss2010}. There has also been research in the case of graphs which can be obtained by self-similar substitutions~\cite{Barbieri2016}. For finitely generated groups, the only groups where the domino problem is known to be decidable are virtually free groups. It is even conjectured that a group has decidable domino problem if and only if it is virtually free. Recent results support the conjecture: decidability of the domino problem is a quasi-isometry invariant for finitely presented groups~\cite{Cohen2014} -- i.e. a geometric property of the group -- and that the conjecture holds true for Baumslag-Solitar groups~\cite{AubrunKari2013}, polycyclic groups~\cite{Jeandel2015poly} and groups of
the form $G_1\times G_2$~\cite{Jeandel2015translation}.

The results of Aubrun and Kari~\cite{Kari, AubrunKari2013} share a common factor: the domino problem is shown to be undecidable on two specific structures that grow exponentially in a regular way --with integer or rational base. But what if the structure grows non-regularly? The reduction from the immortality problem of rational piecewise affine maps used in their work seems difficult to adapt in this case. A class of structures which can grow non-regularly is given by orbit graphs of non-deterministic substitutions (\cref{section.substitutions}). This class of structures includes the hyperbolic plane model of~\cite{Kari}, which can be though of as an orbit graph of the one-letter substitution $0\mapsto00$. Using an idea involving superposition of orbit graphs of substitutions, presented in~\cite{GS} and whose idea they attribute to Lorenzo Sadun, we show that the domino problem is undecidable on all orbit graphs of non-deterministic substitutions that satisfy a technical property (\cref{section.undecidability_orbit_graphs}).

As an application of the previous result we show that the domino problem for the fundamental group of the closed orientable surface of genus 2 is undecidable (\cref{section.surface}). An immediate corollary of that is that the fundamental group of any closed orientable surface of genus at least $2$ has undecidable domino problem.

Finally, we discuss the case of word-hyperbolic groups (\cref{section.remarks}) and show that if a famous conjecture of Gromov  --or a weaker version-- holds, then the only word-hyperbolic groups with decidable domino problem are the virtually free groups.


 \section{Definitions and properties}

 \subsection{Subshifts on graphs}\label{subsection.subshifts_graphs}

In this article we consider bounded degree infinite countable graphs with labels on the edges.
We define a \emph{graph $\Gamma$} to be a triple $(V_\Gamma,E_\Gamma,L_\Gamma)$ where $V_\Gamma$ is an infinite countable set of vertices, $E_\Gamma\subset V_\Gamma^2$ is the set of edges, such that for every vertex $v\in V_\Gamma$ $|\left\{u\in V_\Gamma\mid (u,v)\in E_\Gamma\text{ or }(v,u)\in E_\Gamma\right\}|<M$ where $M$ is some constant, and $L_\Gamma\colon E_\Gamma\to L$ is a labeling function which assigns to every edge a label in a finite set $L$.
Important examples of such graphs are Cayley graphs of finitely generated groups. More precisely, given a finitely generated group $G$ and a set of generators $\mathcal{S}$, its Cayley graph is given by $V_\Gamma = G$, $E_\Gamma = \{ (g,gs) \mid g \in G, s \in \mathcal{S}\}$, $L_\Gamma((g,gs)) = s$. Let $\Gamma=(V_\Gamma,E_\Gamma,L_\Gamma)$ be a graph as defined above. Let $S,T$ be two finite subsets of $V_\Gamma$. A mapping $\phi\colon S\to T$ is a \emph{label preserving graph isomorphism} if $\phi$ is a bijection and
 \begin{itemize}
  \item for all $u,v\in S$, $(u,v)\in E_\Gamma$ if and only if $(\phi(u),\phi(v))\in E_\Gamma$;
  \item for all $u,v\in S$, $L_\Gamma\left((u,v)\right)=L_\Gamma\left((\phi(u),\phi(v))\right)$.
 \end{itemize}

 Let $\A$ be a finite alphabet and $\Gamma$ a graph. The set of mappings from $V_\Gamma$ to $\A$, denoted $\A^\Gamma$, is the set of \emph{configurations} over $\Gamma$. Endowed with the prodiscrete topology, the set $\A^\Gamma$ is compact and metrizable.
 If $S\subset V_\Gamma$ is a finite and connected set of vertices, a \emph{pattern with support $S$} is a mapping $p\colon S\to\A$.
 A pattern $p\colon S\to \A$ \emph{appears} in a configuration $x\in\A^G$ (resp. in a pattern $p'\colon S'\to \A$) if there exists a finite set of vertices $T\subset V_\Gamma$ (resp. $T\subset S'$) and a label preserving graph isomorphism $\phi\colon S\to T$ such that $p_u=x_{\phi(u)}$ (resp. $p_u=p'_{\phi(u)}$) for every $u\in S$.
 In this case, we denote $p\sqsubset x$ (resp. $p\sqsubset p'$).

A \emph{subshift} $X_F\subset \A^\Gamma$ is a set of configurations that avoid some set of forbidden patterns $F$, i.e. $X_F:=\left\{ x\in A^\Gamma\mid \text{ no pattern of }F\text{ appears in }x\right\}$. This notion extends the classical definition of subshift for group actions to arbitrary graphs.
A \emph{subshift of finite type} (SFT) is a subshift for which $F$ can be chosen finite -- equivalently, an SFT may also be defined by a finite set of allowed patterns. In the case where the support of all the forbidden patterns in $F$  consist of two vertices connected by an edge, we say $X_F$ is a nearest neighbor subshift.

Given a graph $\Gamma$ and a finite alphabet $A$, a pattern as defined in the previous section can be encoded by a finite graph, which is an induced finite subgraph of $\Gamma$, with labels on edges and letters from $A$ on vertices. This is what is meant in the sequel by coding of a pattern.

Let $\Gamma$ be a graph in the previous sense. The \emph{domino problem for $\Gamma$} is defined as the set $\texttt{DP}(\Gamma)$ of codings of finite sets of forbidden patterns $F$ such that $X_F\neq\emptyset$. If the set $\texttt{DP}(\Gamma)$ is recursive, we say that \emph{$\Gamma$ has decidable domino problem}, and undecidable domino problem otherwise. In this formalism, Kari's result for the hyperbolic plane~\cite{Kari} is equivalent to the statement that all orbit graphs (see \cref{def:orbit_graph}) of the one-letter substitution $(\{0\},0\mapsto00)$ have undecidable domino problem.

\begin{theorem}[Kari~\cite{Kari}]\label{teorema_jarkko}
	For all orbit graphs of the substitution $(\{0\},0\mapsto00)$ the domino problem is undecidable.
\end{theorem}

Cayley graphs of finitely generated groups form an highly interesting class of graphs. Denote by $\Gamma(G,S)$ the Cayley graph of a group with generating set $S$. If $S$ and $S'$ are two finite generating sets of a group $G$, then $\texttt{DP}(\Gamma(G,S))$ and $\texttt{DP}(\Gamma(G,S'))$ are many-one equivalent~\cite{BertheRigo2018}. In particular, this means that the decidability of the domino problem on a Cayley graph does not depend on the choice of the generating set, thus we can legitimately speak of the decidability of the domino problem of a finitely generated group. A survey on the domino problem for finitely generated groups can be found in~\cite[Chapter 9]{BertheRigo2018}.

\section{Substitutions, orbits and tilings}\label{section.substitutions}

Inspired by~\cite{GS}, we associate a tiling of $\R^2$ with the orbit of an infinite word $w\in \A^{\Z}$ under the action of a substitution, in which every tile codes a production rule of the substitution.

\subsection{Parent functions}

In this section we define parent functions, which will be used to give precise descriptions of orbits of infinite words under the action of a substitution.

A \emph{parent function} $P\colon\Z\to\Z$ is an onto and non-decreasing function. In particular, such a function $P$ satisfies that for every $i\in\Z$, $P(i+1)-P(i)\in\{0,1\}$. Let $u=(u_i)_{i\in\Z}\in\left(\N \setminus \{0\}\right)^\Z$ be a bi-infinite sequence of positive integers. The \emph{accumulation function of $u$} is the function $\Delta\colon\Z\rightarrow\Z$ given by
\[ \Delta(i) = \begin{cases}
		\sum_{k=0}^{i-1} u_k & \text{if } i\geq 1\\
		0 & \text{if } i=0\\
		-\sum_{k=i}^{-1} u_k & \text{if } i\leq -1\\
	     \end{cases}. \]

 Note that the family of discrete intervals $\left(I_k\right)_{k\in\Z}$ where $I_k=[\Delta(k);\Delta(k+1)-1]$ forms a partition of $\Z$.
 If $P$ is a parent function, and if we define the sequence $u$ by $u_i=|P^{-1}(i)|$ for every $i\in\Z$, then we get that $P(j)=i$ for every $j\in[\Delta(i);\Delta(i+1)-1]$, where $\Delta$ is the accumulation function of $u$.

\subsection{Substitution systems}

A \emph{non-deterministic substitution} is a couple $(\A,R)$ where $\A$ is a finite alphabet and $R\subset \A\times \A^*$ is a finite set called the \emph{relation}, and whose elements are called \emph{production rules}. We say that an infinite word $\omega \in \A^{\Z}$ \emph{produces} the word $\omega'\in \A^{\Z}$ with respect to the parent function $P$ if for every $i\in\Z$, one has $(\omega_i,\omega'|_{[\Delta(i);\Delta(i+1)-1]})\in R$, where $\omega'|_{[\Delta(i);\Delta(i+1)-1]}$ is the finite subword of $\omega'$ that appears on indices $\{j\in\Z\mid P(j)=i\}$. In this case we shall extend the above notation and write $(\omega,\omega')\in R$. An \emph{orbit} of a non-deterministic substitution $(\A,R)$ is a set $\left\{ (\omega^i,P_i)\right\}_{i\in\Z}\in \left(\A^{\Z}\times\Z^{\Z}\right)^{\Z}$ such that for every $i\in\Z$, $P_i$ is a parent function, and the word $\omega^i$ produces the word $\omega^{i+1}$ with respect to $P_i$. A \emph{deterministic substitution} (or substitution for short) is
a non-deterministic substitution where the relation is a function. A non-deterministic substitution $(\A,R)$ \emph{has an expanding eigenvalue} if there exist $\lambda>1$ and $v\colon\A\to\R^{+}\setminus \{0\}$ such that for every $(a,w)\in R$,
$$
\lambda\cdot v(a)=\sum_{i=1}^{|w|} v(w_i).
$$

\begin{example} The substitution given by $\A = \{0\}$ and the production rule $0 \mapsto 00$ has an expanding eigenvalue. This can be verified by choosing $\lambda = 2$ and $v(0)=1$.\end{example}

\subsection{Orbits as tilings of \texorpdfstring{$\R^2$}{R2}}\label{subsection:orbits_tilings}

\begin{figure}[!ht]
\begin{center}
\begin{tikzpicture}[scale=1]
\draw[thick] (0,0) -- (0,3);
\draw[thick] (0,3) to [controls=+(45:1.7) and +(135:1.7)] (6,3);
\draw[thick] (6,3) -- (6,0);
\draw[thick] (0,0) to [controls=+(30:0.75) and +(150:0.75)] (2,0);
\draw[thick] (2,0) to [controls=+(30:0.25) and +(150:0.25)] (3,0);
\draw[thick] (5,0) to [controls=+(30:0.25) and +(150:0.25)] (6,0);
\foreach \Point in {(0,0), (0,3), (6,3), (6,0), (2,0), (3,0), (5,0)}{
    \node[scale=1] at \Point {\textbullet};
}
\draw node at (3,2) {$(a,w)$-tile};
\draw node at (-0.5,3.25) {\scalebox{0.8}{$(x,y)$}};
\draw node at (-0.5,-0.25) {\scalebox{0.8}{$(x,y-\log(\lambda))$}};

\draw[<->] (0.1,4.25) -- (5.9,4.25);
\draw node at (3,4.5) {$v(a)\cdot e^y$};

\draw[<->] (-0.5,0) -- (-0.5,3);
\draw node at (-1.25,1.5) {\scalebox{0.8}{$\log(\lambda)$}};

\draw[<->] (0.1,-0.5) -- (1.9,-0.5);
\draw node at (1,-1) {\scalebox{0.8}{$\frac{1}{\lambda}v(w_1)\cdot e^y$}};

\draw[<->] (2.1,-0.5) -- (2.9,-0.5);
\draw node at (2.5,-1) {\scalebox{0.8}{$\frac{1}{\lambda}v(w_2)\cdot e^y$}};

\draw node at (4,-0.5) {$\dots$};
\draw node at (4,0) {$\dots$};
\draw node at (4,-1) {$\dots$};

\draw[<->] (5.1,-0.5) -- (5.9,-0.5);
\draw node at (5.5,-1) {\scalebox{0.8}{$\frac{1}{\lambda}v(w_k)\cdot e^y$}};

\draw node at (8,2) {$(a,w_1\dots w_k)\in R$};
\end{tikzpicture}
\end{center}
\caption{An $(a,w)$-tile for some production rule $(a,w)\in R$ with $w=w_1\dots w_k$.}
\label{figure.ab_tile}
\end{figure}
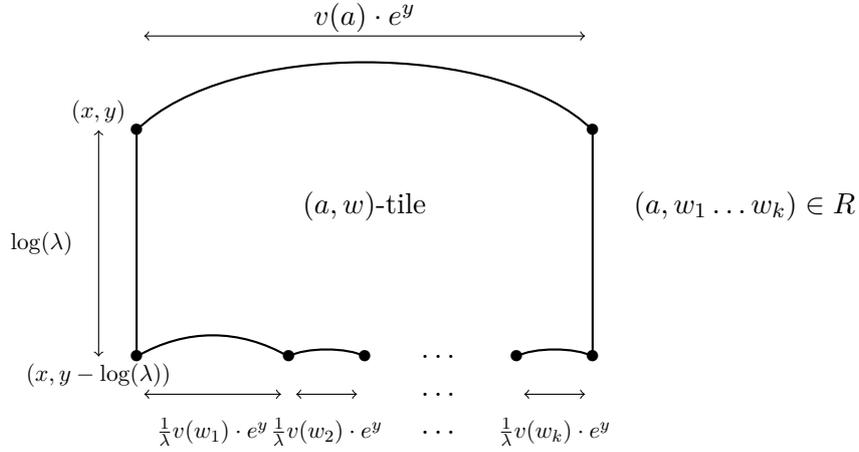

Let $(\A,R)$ be a non-deterministic substitution with an expanding eigenvalue $\lambda>1$ and $v\colon\A\to\R^{+}\setminus \{0\}$. For every production rule $(a,w)\in R$, define the \emph{$(a,w)$-tile in position $(x,y)\in\R^2$} as the square polygon with $|w|+3$ edges pictured in \cref{figure.ab_tile}, where $w=w_{1}\dots w_{k}$ (horizontal edges are curved to be more visible, but are in fact just straight lines).

\begin{remark*}\label{remark.tile_rectangular}
The length of the top edge and the sum of lengths of bottom edges of this tile are the same. Since $(\A,R)$ has an expanding eigenvalue $\lambda>1$ with $v$, one has
$$\sum_{j=1}^{k}\frac{1}{\lambda} v(w_i)\cdot e^y =\frac{e^y}{\lambda}\cdot\lambda\cdot v(a)=v(a)\cdot e^y,$$
so that the bottom right vertex $(x+\frac{1}{\lambda}\left(v(w_1)+\dots+v(w_k)\right) e^y,y-\log(\lambda))$ is indeed $(x+v(a)\cdot e^y,y-\log(\lambda))$.
\end{remark*}

The \emph{$(\A,R)$-tiles} is the set of all $(a,w)$-tiles in position $(x,y)$ for all possible $(a,w)\in R$ and $x,y\in\R$. Given an orbit $\Omega = \left\{ (\omega^i,P_i)\right\}_{i\in\Z}$ for $(\A,R)$ a \emph{tiling of $\R^2$ for $\Omega$} is a function $\Psi_{\Omega}\colon \Z^2\to\R^2$ such that for every $(i,j)\in \Z^2$ we have:

\begin{itemize}
	\item $\Psi_\Omega(i,j)=(x,y)$ if and only if $\Psi_\Omega(i,j+1)=(x+v((\omega^i)_j)\cdot e^y,y)$;
	\item $\Psi_\Omega(i,j)=(x,y)$ if and only if $\Psi_\Omega(i+1,\min{P^{-1}_{i+1}{(j)}})=(x,y-\log(\lambda))$.
\end{itemize}

\begin{figure}[!ht]
 	\begin{center}
 			\begin{tikzpicture}[scale=0.85]
	\draw[thick] (-8,3) to [controls=+(45:1.2) and +(135:1.2)] (-2,3);
	\draw[thick] (-5,0) to [controls=+(30:1) and +(150:1)] (-2,0);
	\draw[thick] (-4,-3) to [controls=+(30:0.75) and +(150:0.75)] (-2,-3);

	\draw[thick] (6,3) to [controls=+(45:1) and +(135:1)] (9,3);
	\draw[thick] (6,0) to [controls=+(30:1) and +(150:1)] (8,0);
	\draw[thick] (6,-3) to [controls=+(30:0.5) and +(150:0.5)] (7,-3);
	\draw node at (2,2) {$\omega^{-1}_0$-tile};
	\draw[thick] (-2,0) -- (-2,3);
	\draw[thick] (-2,3) to [controls=+(45:1.7) and +(135:1.7)] (6,3);
	\draw[thick] (6,3) -- (6,0);
	\draw[thick] (-2,0) to [controls=+(30:0.75) and +(150:0.75)] (0,0);
	\draw[thick] (0,0) to [controls=+(30:0.75) and +(150:0.75)] (2,0);
	\draw[thick] (2,0) to [controls=+(30:1) and +(150:1)] (5,0);
	\draw[thick] (5,0) to [controls=+(30:0.25) and +(150:0.25)] (6,0);
	\draw node at (-1,-1) {$\omega^{0}_{-2}$-tile};
	\draw[thick] (-2,0) -- (-2,-3);
	\draw[thick] (-2,-3) to [controls=+(45:0.5) and +(135:0.5)] (-1,-3);
	\draw[thick] (-1,-3) to [controls=+(45:0.5) and +(135:0.5)] (0,-3);
	\draw[thick] (0,-3) -- (0,0);
	\draw node at (1,-1) {$\omega^{0}_{-1}$-tile};
	\draw[thick] (0,0) -- (0,-3);
	\draw[thick] (0,-3) to [controls=+(45:0.5) and +(135:0.5)] (1,-3);
	\draw[thick] (1,-3) to [controls=+(45:0.5) and +(135:0.5)] (2,-3);
	\draw[thick] (2,-3) -- (2,0);
	\draw node at (3.5,-1) {$\omega^{0}_{0}$-tile};
	\draw[thick] (2,-3) to [controls=+(45:0.5) and +(135:0.5)] (3,-3);
	\draw[thick] (3,-3) to [controls=+(45:0.5) and +(135:0.5)] (5,-3);
	\draw[thick] (5,-3) -- (5,0);
	\draw node at (5.5,-1) {$\omega^{0}_{1}$-tile};
	\draw[thick] (5,-3) to [controls=+(45:0.25) and +(135:0.25)] (5.75,-3);
	\draw[thick] (5.75,-3) to [controls=+(45:0.15) and +(135:0.15)] (6,-3);
	\draw[thick] (6,-3) -- (6,0);

	\foreach \Point in {(2,0), (-2,3), (6,3)}{
			\node[scale=1] at \Point {\textbullet};
	}
	\draw node at (2,0.5) {$(0,0)$};
	\end{tikzpicture}
 	\end{center}
 	\caption{A tiling $\Psi_\Omega$ of an orbit into $\R^2$.}
 	\label{figure.embedding_orbit}
 \end{figure}
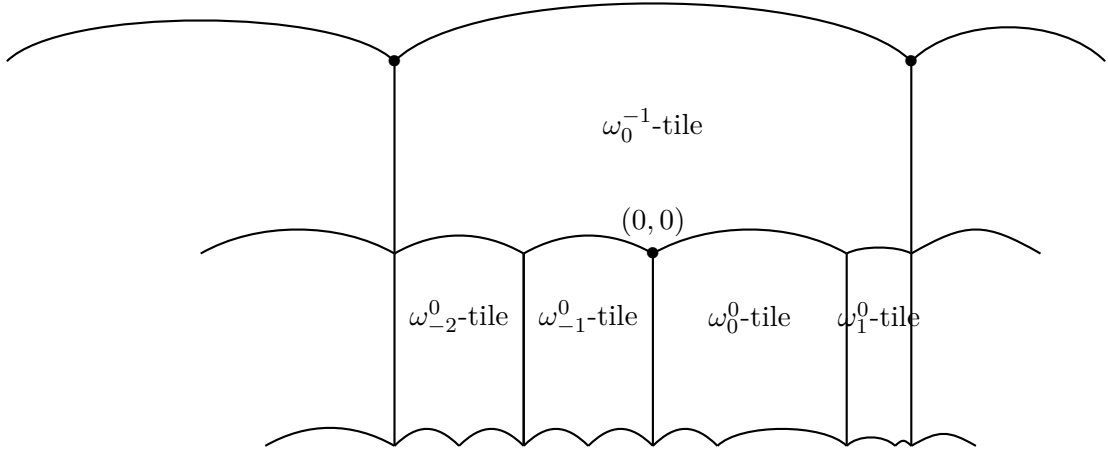

Note that by the previous remark, the collection of $(\A,R)$ obtained by putting an $((a^i)_j,a^{i+1}|_{[\Delta_{i+1}(j);\Delta_{i+1}(j+1)-1]})$-tile at position $\Psi_{\Omega}(i,j)$ defines a tiling of $\R^2$, that is, the collection of square polygons covers $\R^2$ and has pairwise disjoint interiors. See~\cref{figure.embedding_orbit}. Furthermore, fixing one position, say $\Psi_{\Omega}(0,0)=(0,0)$ defines the function $\Psi_{\Omega}$ completely. It follows that for a substitution with an expanding eigenvalue, there is always one tiling for it.

\begin{proposition}\label{proposition.orbits_implies_tiling}
	If a substitution $(\A,R)$ has an expanding eigenvalue, then for every orbit $\Omega$ of $(\A,R)$ there exists a tiling $\Psi_{\Omega}$ for $\Omega$.
\end{proposition}

%
%
%

We shall use that fact that tilings exist in~\cref{subsection.superposition} to define the superposition of two orbits from two different substitutions and prove the non-emptiness of a subshift of finite type.

\section{Undecidability of the domino problem on orbit graphs}\label{section.undecidability_orbit_graphs}

Let $(\A,R)$ be a non-deterministic substitution. Denote $M=\max_{(a,w)\in R} |w|$.

\begin{definition}
	\label{def:orbit_graph}
	The \emph{orbit graph} associated with the orbit $\Omega=\left\{ (\omega^i,P_i)\right\}_{i\in\Z}$ of $(\A,R)$ is the graph $\Gamma_\Omega$ with set of vertices $\Z^2$, edges $E_\Omega$ and labeling function $L_\Omega\colon E_\Omega\to\{\texttt{next}\}\cup [0;M-1]$ given by
	\begin{itemize}
		\item for every $i,j\in\Z$, $((i,j),(i,j+1))\in E_\Omega$ and $L_\Omega\left(((i,j),(i,j+1))\right)=\texttt{next}$;
		\item for every $i\in\Z$ and every $k\in[\Delta_{i+1}(j);\Delta_{i+1}(j+1)-1]$, $((i,j),(i+1,k))\in E_\Omega$ and $L_\Omega\left(((i,j),(i+1,k))\right)=k-\Delta_{i+1}(j)$,
	\end{itemize}
	where $\Delta_i$ is the accumulation function associated with $\left(|P_i^{-1}(j)|\right)_{j\in\Z}$ for every $i\in\Z$.
\end{definition}

Note that $\Gamma_{\Omega}$ depends uniquely upon the parent functions $\{P_i\}_{i\in \Z}$ and not on $\{\omega^i\}_{i\in \Z}$. However, we implicitly require that the sequence of parent functions defines an orbit $\Omega = \left\{ (\omega^i,P_i)\right\}_{i\in\Z}$ of $(\A,R)$.

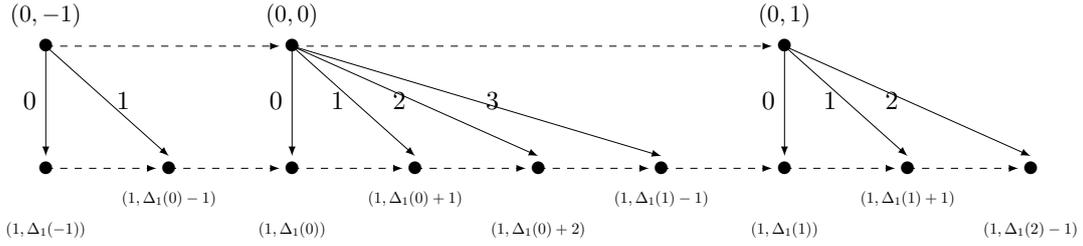
\begin{figure}[!ht]
	\begin{center}
		\scalebox{0.9}{\begin{tikzpicture}[scale=0.9]

\foreach \Point in {(-6,0), (-6,-2), (-4,-2), (-2,0), (-2,-2), (0,-2), (2,-2), (4,-2), (6,0), (6,-2), (8,-2), (10,-2)}{
		\node[scale=1.25] at \Point {\textbullet};
}

\draw[thin,->,>=latex,dashed] (-6,0) to (-2-0.2,0);
\draw[thin,->,>=latex,dashed] (-2,0) to (6-0.2,0);
\draw[thin,->,>=latex,dashed] (-6,-2) to (-4-0.2,-2);
\draw[thin,->,>=latex,dashed] (-4,-2) to (-2-0.2,-2);
\draw[thin,->,>=latex,dashed] (-2,-2) to (0-0.2,-2);
\draw[thin,->,>=latex,dashed] (0,-2) to (2-0.2,-2);
\draw[thin,->,>=latex,dashed] (2,-2) to (4-0.2,-2);
\draw[thin,->,>=latex,dashed] (4,-2) to (6-0.2,-2);
\draw[thin,->,>=latex,dashed] (6,-2) to (8-0.2,-2);
\draw[thin,->,>=latex,dashed] (8,-2) to (10-0.2,-2);

\node[scale=0.9]  at (-6,0.5) {$(0,-1)$};
\node[scale=0.6]  at (-6,-3) {$(1,\Delta_{1}(-1))$};
\node[scale=0.6]  at (-4,-2.5) {$(1,\Delta_{1}(0)-1)$};
\draw[thin,->,>=latex] (-6,0) -- node[left] {$0$} (-6,-1.8);
\draw[thin,->,>=latex] (-6,0) -- node[right] {$1$} (-4,-1.8);

\node[scale=0.9]  at (-2,0.5) {$(0,0)$};
\node[scale=0.6]  at (-2,-3) {$(1,\Delta_{1}(0))$};
\node[scale=0.6]  at (0,-2.5) {$(1,\Delta_{1}(0)+1)$};
\node[scale=0.6]  at (2,-3) {$(1,\Delta_{1}(0)+2)$};
\node[scale=0.6]  at (4,-2.5) {$(1,\Delta_{1}(1)-1)$};
\draw[thin,->,>=latex] (-2,0) -- node[left] {$0$} (-2,-1.8);
\draw[thin,->,>=latex] (-2,0) -- node[left] {$1$} (0,-2+0.2);
\draw[thin,->,>=latex] (-2,0) -- node[left] {$2$} (2,-2+0.2);
\draw[thin,->,>=latex] (-2,0) -- node[right] {$3$} (4,-2+0.2);

\node[scale=0.9]  at (6,0.5) {$(0,1)$};
\node[scale=0.6]  at (6,-3) {$(1,\Delta_{1}(1))$};
\node[scale=0.6]  at (8,-2.5) {$(1,\Delta_{1}(1)+1)$};
\node[scale=0.6]  at (10,-3) {$(1,\Delta_{1}(2)-1)$};
\draw[thin,->,>=latex] (6,0) -- node[left] {$0$} (6,-2+0.2);
\draw[thin,->,>=latex] (6,0) -- node[left] {$1$} (8,-2+0.2);
\draw[thin,->,>=latex] (6,0) -- node[left] {$2$} (10,-2+0.2);

\end{tikzpicture}}
	\end{center}
	\caption{Part of an orbit graph. Dashed arrow are edges of the graph labeled with $\texttt{next}$.}
	\label{figure.orbit_graph}
\end{figure}

The goal of this section is to show that the domino problem of any orbit graph associated to an orbit of a non-deterministic substitution with an expanding eigenvalue is undecidable.
In order to prove this we show a variation of the "Technical Lemma" of Cohen and Goodman-Strauss~\cite{GS}.
Their lemma takes two primitive expansive deterministic substitutions $(\A,\sigma)$ and $(\A',\tau)$ and produces a non-deterministic one that simulates the orbits of $(\A,\sigma)$ and $(\A',\tau)$ in its orbits.
Their proof uses the idea of superposing two tilings associated to the substitutions and coding their intersections. For our purposes, we shall consider any orbit $\Omega$ of a non-deterministic substitution $(\A,R)$ with an expanding eigenvalue $\lambda$ and construct a subshift of finite type $Y$ in $\Gamma_\Omega$ which encodes an orbit graph of the specific substitution $(\{0\}, 0\mapsto 00)$. For technical reasons that will become clear during the proof, we shall first consider the case where $\lambda > 2$ and then deduce the general case from this case.

\subsection{Superposition of orbits}\label{subsection.superposition}

Let us fix a non-deterministic substitution $(\A,R)$ with an expanding eigenvalue $\lambda>2$. Without loss of generality, we may choose the function $v\colon \A \to \R^{+}\setminus \{0\}$ associated to $\lambda$ such that $v(a)>4$ for each $a \in \A$.

Let $\Omega=\left\{ (\omega^i,P_i)\right\}_{i\in\Z}$ be an orbit of $(\A,R)$. We shall construct a finite alphabet $\B$ and a finite set of forbidden patterns $F$ such that the subshift $Y \subset \B^{\Gamma_{\Omega}}$ defined by the set of forbidden patterns $F$ has the following properties:
\begin{enumerate}
  \item $Y$ is non-empty,
  \item every configuration $y \in Y$ encodes an orbit graph of the substitution $(\{0\}, 0\mapsto 00)$.
\end{enumerate}

We first give an informal description of the alphabet $\B$. Consider an orbit $\Omega=\left\{ (\omega^i,P_i)\right\}_{i\in\Z}$ of $(\A,R)$ and $\Xi=\left\{ ((0^{\infty})^i,Q_i)\right\}_{i\in\Z}$ an orbit of $(\{0\}, 0\mapsto 00)$. By~\cref{proposition.orbits_implies_tiling} both of these orbits can be realized as tilings of $\R^2$. Symbols from $\B$ will encode non-empty finite regions of the tiling with $(\{0\}, 0\mapsto 00)$ that are witnessed by $(\A,R)$-tiles. These regions will be chosen in such a way that their union recovers the whole tiling and they are pairwise disjoint. More precisely, the alphabet $\B$ will consist of

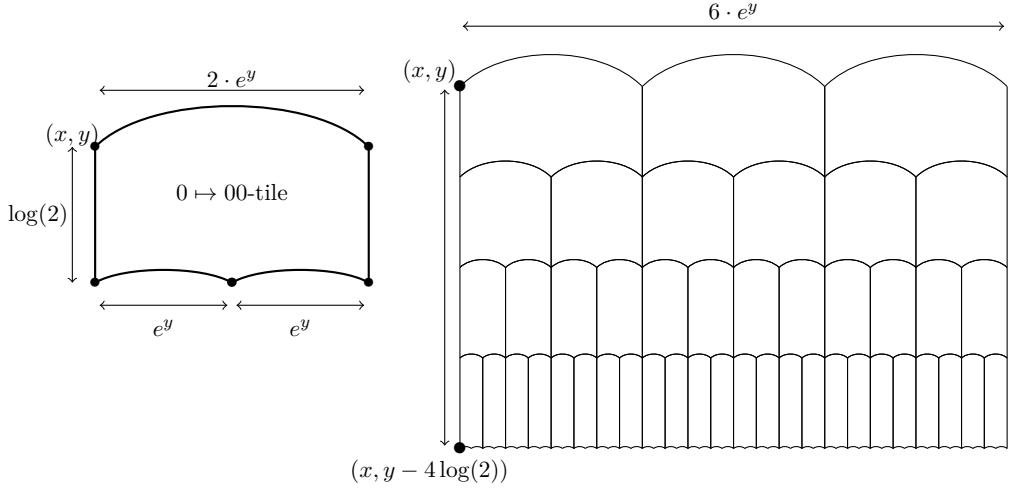
\begin{figure}[!ht]
  	\begin{center}
  			\begin{tikzpicture}

	\begin{scope}[scale=0.6]
	\draw[thick] (0,0) -- (0,3);
	\draw[thick] (0,3) to [controls=+(45:1.7) and +(135:1.7)] (6,3);
	\draw[thick] (6,3) -- (6,0);
	\draw[thick] (0,0) to [controls=+(30:0.75) and +(150:0.75)] (3,0);
	\draw[thick] (3,0) to [controls=+(30:0.75) and +(150:0.75)] (6,0);
	\foreach \Point in {(0,0), (0,3), (6,3), (6,0), (3,0), (3,0)}{
	    \node[scale=0.8] at \Point {\textbullet};
	}
	\draw node at (3,2) {\scalebox{0.8}{$0\mapsto00$-tile}};
	\draw node at (-0.5,3.25) {\scalebox{0.8}{$(x,y)$}};

	\draw[<->] (0.1,4.25) -- (5.9,4.25);
	\draw node at (3,4.5) {\scalebox{0.8}{$2\cdot e^y$}};

	\draw[<->] (-0.5,0) -- (-0.5,3);
	\draw node at (-1.25,1.5) {\scalebox{0.8}{$\log(2)$}};

	\draw[<->] (0.1,-0.5) -- (2.9,-0.5);
	\draw node at (1.5,-1) {\scalebox{0.8}{$e^y$}};

	\draw[<->] (3.1,-0.5) -- (5.9,-0.5);
	\draw node at (4.5,-1) {\scalebox{0.8}{$e^y$}};

	\end{scope}

	\begin{scope}[scale=0.4, shift={(12,3.5)}]
	\foreach \Point in {(0,0), (6,0), (12,0)}{
	\begin{scope}[shift={\Point}]
	\draw (0,0) -- (0,3);
	\draw (0,3) to [controls=+(45:2) and +(135:2)] (6,3);
	\draw (6,3) -- (6,0);
	\draw (0,0) to [controls=+(45:1) and +(135:1)] (3,0);
	\draw (3,0) to [controls=+(45:1) and +(135:1)] (6,0);
	\end{scope}
	}
	\foreach \Point in {(0,-3), (3,-3), (6,-3), (9,-3), (12,-3), (15,-3)}{
	\begin{scope}[shift={\Point},scale=0.5]
	\draw (0,0) -- (0,6);
	\draw (0,6) to [controls=+(45:2) and +(135:2)] (6,6);
	\draw (6,6) -- (6,0);
	\draw (0,0) to [controls=+(45:1) and +(135:1)] (3,0);
	\draw (3,0) to [controls=+(45:1) and +(135:1)] (6,0);
	\end{scope}
	}
	\foreach \Point in {(0,-6), (1.5,-6), (3,-6), (4.5,-6), (6,-6), (7.5,-6), (9,-6), (10.5,-6), (12,-6), (13.5,-6), (15,-6), (16.5,-6)}{
	\begin{scope}[shift={\Point},scale=0.25]
	\draw (0,0) -- (0,12);
	\draw (0,12) to [controls=+(45:2) and +(135:2)] (6,12);
	\draw (6,12) -- (6,0);
	\draw (0,0) to [controls=+(45:1) and +(135:1)] (3,0);
	\draw (3,0) to [controls=+(45:1) and +(135:1)] (6,0);
	\end{scope}
	}
	\foreach \Point in {(0,-9), (0.75,-9), (1.5,-9), (2.25,-9), (3,-9), (3.75,-9), (4.5,-9), (5.25,-9), (6,-9), (6.75,-9), (7.5,-9), (8.25,-9), (9,-9), (9.75,-9), (10.5,-9), (11.25,-9), (12,-9), (12.75,-9), (13.5,-9), (14.25,-9), (15,-9), (15.75,-9), (16.5,-9), (17.25,-9)}{
	\begin{scope}[shift={\Point},scale=0.125]
	\draw (0,0) -- (0,24);
	\draw (0,24) to [controls=+(45:2) and +(135:2)] (6,24);
	\draw (6,24) -- (6,0);
	\draw (0,0) to [controls=+(45:1) and +(135:1)] (3,0);
	\draw (3,0) to [controls=+(45:1) and +(135:1)] (6,0);
	\end{scope}
	}
	\foreach \Point in {(0,3),(0,-9)}{
	\node[scale=1] at \Point {\textbullet};
	}
	\draw node at (-1,3.5) {\scalebox{0.8}{$(x,y)$}};
	\draw node at (-1,-9.75) {\scalebox{0.8}{$(x,y-4\log(2))$}};
	\draw[<->] (0.1,5) -- (17.9,5);
	\draw node at (9,5.5) {\scalebox{0.8}{$6\cdot e^y$}};
	\draw[<->] (-0.5,2.9) -- (-0.5,-8.9);
	\end{scope}

	\end{tikzpicture}
  	\end{center}
  	\caption{A $0\mapsto00$-tile, and a $(3,4)$-block in position $(x,y)\in\R^2$.}
  	\label{figure.simplified_0mapsto00_tile}
  \end{figure}

\begin{itemize}
	\item A production rule $(a,w)\in R$ describing the type of $(\A,R)$-tile.
	\item Two integers $(h,t)$, describing a finite region of the tiling associated to an orbit of $(\{0\}, 0\mapsto 00)$. The integer $h$ represents the number of $(\{0\}, 0\mapsto 00)$ tiles than can fit vertically in the current type of $(\A,R)$-tile and $t$ is the number that fits horizontally on the top edge.
	\item A tuple of $|w|$ pairs of integers $[(b_0,s_0),(b_1,s_1)\dots,(b_{|w|-1},s_{|w|-1})]$ which describes how to locally paste the region with its neighboring regions. More precisely, it contains all information needed to recover the function $Q_i$ from the finite coded regions.
  Each $b_i$ represents the index of the $(\{0\}, 0\mapsto 00)$-tile that intersects the left corner of the $i$-th bottom edge of the $(\A,R)$-tile (starting from 0), and $s_i$ its binary label, depending if the vertex intersects the left or right child of $b_i$ (see \cref{figure.simplified_0mapsto00_tile}).
\end{itemize}

The \emph{$0\mapsto00$-tile in position $(x,y)\in\R^2$} is the square polygon whose five vertices have coordinates $(x,y)$, $(x,y-\log(2))$, $(x+e^y,y-\log(2))$, $(x+2\cdot e^y,y-\log(2))$ and $(x+2\cdot e^y,y)$ as pictured on the left of \cref{figure.simplified_0mapsto00_tile}.
The width of these tiles depends on their position --more precisely only on their second coordinate-- but their height does not and is always $\log(2)$.

By~\cref{proposition.orbits_implies_tiling} we can tile the plane with this family of tiles by putting tiles vertex to vertex, each tile having a left and a right neighbor, two children and one parent. In the sequel we will be interested in blocks of such tiles.
The \emph{$(h,t)$-block in position $(x,y)\in\R^2$} is a pattern of width $2te^y$ and height $h\log(2)$, filled in with tiles as pictured on \cref{figure.simplified_0mapsto00_tile}, and whose top left vertex has coordinates $(x,y)$.
Similarly, by~\cref{proposition.orbits_implies_tiling} we can also tile $\R^2$ with $(\A,R)$-tiles and speak of the $(a,w)$-tile at position $(x,y)$ as in~\cref{figure.ab_tile}.

\begin{figure}[!ht]
  	\begin{center}
  			\begin{tikzpicture}[scale=0.4]
	\foreach \Point in {(0,0),(6,0)}{
	\begin{scope}[color=black!50,shift={\Point}]
	\draw[fill=blue!30] (0,0) -- (0,3) to [controls=+(45:2) and +(135:2)] (6,3) -- (6,0) to [controls=+(135:1) and +(45:1)] (3,0) to [controls=+(135:1) and +(45:1)] (0,0) -- cycle;
	\end{scope}
	}

	\foreach \Point in {(12,0)}{
		\begin{scope}[color=black!50,shift={\Point}]
		\draw (0,0) -- (0,3) to [controls=+(45:2) and +(135:2)] (6,3) -- (6,0) to [controls=+(135:1) and +(45:1)] (3,0) to [controls=+(135:1) and +(45:1)] (0,0) -- cycle;
		\end{scope}
	}

	\foreach \Point in {(0,-3),(3,-3),(6,-3),(9,-3)}{
	\begin{scope}[color=black!50,shift={\Point},scale=0.5]
	\draw[fill=blue!30] (0,0) -- (0,6) to [controls=+(45:2) and +(135:2)] (6,6) -- (6,0) to [controls=+(135:1) and +(45:1)] (3,0) to [controls=+(135:1) and +(45:1)] (0,0) -- cycle;
	\end{scope}
	}

	\foreach \Point in {(12,-3),(15,-3)}{
		\begin{scope}[color=black!50,shift={\Point},scale=0.5]
		\draw (0,0) -- (0,6) to [controls=+(45:2) and +(135:2)] (6,6) -- (6,0) to [controls=+(135:1) and +(45:1)] (3,0) to [controls=+(135:1) and +(45:1)] (0,0) -- cycle;
		\end{scope}
	}

	\foreach \Point in {(0,-6), (1.5,-6), (3,-6), (4.5,-6), (6,-6), (7.5,-6), (9,-6), (10.5,-6)}{
	\begin{scope}[color=black!50,shift={\Point},scale=0.25]
	\draw[fill=blue!30] (0,0) -- (0,12) to [controls=+(45:2) and +(135:2)] (6,12) -- (6,0) to [controls=+(135:1) and +(45:1)] (3,0) to [controls=+(135:1) and +(45:1)] (0,0) -- cycle;
	\end{scope}
	}

	\foreach \Point in {(12,-6), (13.5,-6), (15,-6), (16.5,-6)}{
		\begin{scope}[color=black!50,shift={\Point},scale=0.25]
		\draw (0,0) -- (0,12) to [controls=+(45:2) and +(135:2)] (6,12) -- (6,0) to [controls=+(135:1) and +(45:1)] (3,0) to [controls=+(135:1) and +(45:1)] (0,0) -- cycle;
		\end{scope}
	}

	\foreach \Point in {(0,-9), (0.75,-9), (1.5,-9),(2.25,-9), (3,-9), (3.75,-9), (4.5,-9), (5.25,-9), (6,-9), (6.75,-9), (7.5,-9), (8.25,-9), (9,-9), (9.75,-9), (10.5,-9), (11.25,-9), (12,-9), (12.75,-9), (13.5,-9), (14.25,-9), (15,-9), (15.75,-9), (16.5,-9), (17.25,-9)}{
	\begin{scope}[color=black!50,shift={\Point},scale=0.125]
	\draw (0,0) -- (0,24) to [controls=+(45:2) and +(135:2)] (6,24) -- (6,0) to [controls=+(135:1) and +(45:1)] (3,0) to [controls=+(135:1) and +(45:1)] (0,0) -- cycle;
	\end{scope}
	}

  \begin{scope}[shift={(0.3,2.7)}]
    \draw node at (0,-6)    [opacity=0.6] {{\tiny \color{Blue} $0$}};
    \draw node at (1.5,-6)  [opacity=0.6] {{\tiny \color{Blue} $1$}};
    \draw node at (3,-6)    [opacity=0.6] {{\tiny \color{Blue} $2$}};
    \draw node at (4.5,-6)  [opacity=0.6] {{\tiny \color{Blue} $3$}};
    \draw node at (6,-6)    [opacity=0.6] {{\tiny \color{Blue} $4$}};
    \draw node at (7.5,-6)  [opacity=0.6] {{\tiny \color{Blue} $5$}};
    \draw node at (9,-6)    [opacity=0.6] {{\tiny \color{Blue} $6$}};
    \draw node at (10.5,-6) [opacity=0.6] {{\tiny \color{Blue} $7$}};
  \end{scope}

	\draw[thick,opacity=0.5,fill=black!25] (1.8,-8.5) -- (1.8,2.8) to [controls=+(45:2) and +(135:2)] (13.1,2.8) -- (13.1,-8.5) to [controls=+(135:1) and +(45:1)] (10,-8.5) to [controls=+(135:1.3) and +(45:1.3)] (4,-8.5) to [controls=+(135:0.7) and +(45:0.7)] (1.8,-8.5);


	\draw node at (8,-1.5) {{$(\mathcal{A},R)$-tile}};
	\draw node at (8,-2.5) {{$(a,w_1w_2w_3)$}};
	\draw[fill=white] (2.0,-4.5-0.5) rectangle + (0.5,0.5);
	\draw[fill=white] (3.5,-4.5-0.5) rectangle + (0.5,0.5);
	\draw[fill=white] (9.5,-4.5-0.5) rectangle + (0.5,0.5);
	\draw node at (2.25,-4.25-0.5) {{{\tiny $1$}}};
	\draw node at (3.75,-4.25-0.5) {{{\tiny $2$}}};
	\draw node at (9.75,-4.25-0.5) {{{\tiny $6$}}};
	\draw[->, thick] (2.25, -4.5-0.5) to (1.8,-8.5);
	\draw[->, thick] (3.75, -4.5-0.5) to (4,-8.5);
	\draw[->, thick] (9.75, -4.5-0.5) to (10,-8.5);

  \node[draw,circle,inner sep=1.5pt,fill] at (0,3) {};
  \draw node at (0,3+0.5) [left] {$(x,y)$};

  \node[draw,circle,inner sep=1.5pt,fill] at (1.8,2.8) {};
  \draw node at (1.8,2.8-0.5) [right] {$(x+\widetilde{x},y-\widetilde{y})$};

	\end{tikzpicture}
  	\end{center}
  	\caption{The blue $(3,2)$-block intersects the $(\A,R)$-tile in the manner described above.
    The bottom vertices of the $(\A,R)$-tile have horizontal coordinates corresponding to tiles on the last line of the $0\mapsto00$-block. Namely the 2nd (index 1), the 3rd (index 2) and the 7th (index 6).
    These vertices are respectively on the left, right and right child of these $0\mapsto00$-tiles.
    Therefore, the associated symbol of $\B$ is given by: $\left((a,w_1w_2w_3),(3,2),[(1,0),(2,1),(6,1)]\right)$.
    }
  	\label{figure.example_intersection_block}
  \end{figure}
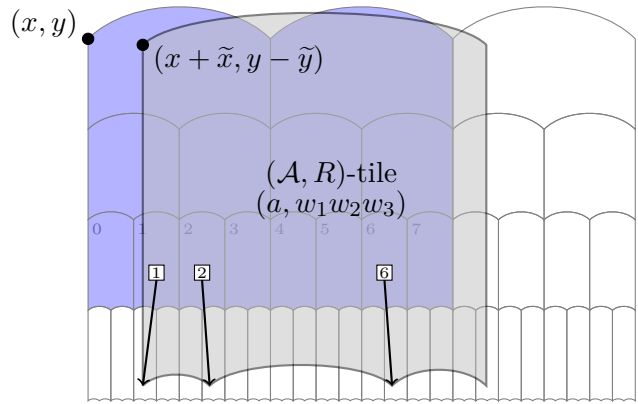

	Let $(x,y)\in\R^2$, $\widetilde{x}\in[0;2\cdot e^y[$ and $\widetilde{y}\in[0;\log(2)[$. We want to consider the largest values $(h,t)$ such that an $(\A,R)$-tile at position $(x+\widetilde{x},y-\widetilde{y})$ intersects the interior of the top-left tile of the $(h,t)$-block at $(x,y)$ and the bottom right corner $(x + 2te^y, y-h\log(2))$ of the $(h,t)$-block is contained in the $(\A,R)$-tile (see~\cref{figure.example_intersection_block}).

  We also need information of how to paste consecutive coded blocks. Each integer $b_i$ for $i \in \{0,\dots,|w|-1\}$ will code the number counted from left to right of the tile in the bottom row of the $(h,t)$-block which is the parent of the top-left tile of the block coded by the $i$-th son of $(a,w)$. The value $s_i \in \{0,1\}$ indicates whether the top-left tile of the block coded by the $i$-th son of $(a,w)$ is the left $(0)$ or right $(1)$ son (see the bottom of~\cref{figure.example_intersection_block}).

  \subsubsection*{Definition of the alphabet $\B$}

  We now define the alphabet $\B$ formally. A symbol \[b = \left((a,w), (h,t), [(b_0,s_0),\dots,(b_{|w|-1},s_{|w|-1})]\right) \] is in $\B$ if and only if $(a,w)\in R$ and there exists $(x,y)\in\mathbb{R}^2$, $\widetilde{x}\in\left[ 0,2\cdot e^y\right[$ and $\widetilde{y}\in\left[ 0,\log(2)\right[$ such that
  \begin{enumerate}
  	\item there is a $(a,w)$-tile in position $(x+\widetilde{x},y-\widetilde{y})$;
  	\item $h=\left\lfloor \frac{\log(\lambda)+\widetilde{y}}{\log(2)} \right\rfloor$;
  	\item $t=\left\lfloor \frac{\widetilde{x}+v(a)\cdot e^{y-\widetilde{y}}}{2\cdot e^{y}}\right\rfloor$;
  	\item For every $i \in \{0,\dots, |w|-1\}$,
  	\begin{itemize}
  		\item $b_i= \left\lfloor\frac{\widetilde{x}+e^{y-\widetilde{y}-\log(\lambda)}\sum_{k=1}^{i}v(w_k)}{2e^{y-(h-1)\log(2)}}\right\rfloor$;
  		\item $s_i = \left\lfloor\frac{\widetilde{x}+e^{y-\widetilde{y}-\log(\lambda)}\sum_{k=1}^{i}v(w_k)}{2e^{y-h\log(2)}}\right\rfloor \mod{2}.$
  	\end{itemize}
  \end{enumerate}
  The values $h$ and $t$ represent the height and width of the largest block of $0\mapsto00$-tiles that fit in the $(\A,R)$-tile as shown on~\cref{figure.example_intersection_block}. The numbers $b_i$ code the number of the $0\mapsto00$-tile on the bottom row of the $(h,t)$-block (from left to right starting at $0$) such that the horizontal coordinate of the $i$-th bottom vertex of the $(\A,R)$-tile is contained in it. The numbers $s_i$ satisfy that the tile indicated by $b_i$ is connected to the top-left tile coded by the $i$-th son of the $(\A,R)$-tile by the label $s_i$.

  Remark that as $\lambda >2$, we have $h\geq 1$. Furthermore, $h$ can take only two consecutive integer values. Similarly, for a given production rule $(a,w)\in R$, the bounds impose that $t$ is an integer satisfying $\left\lfloor\frac{v(a)}{4}\right\rfloor \leq t \leq \left\lfloor1+\frac{v(a)}{2}\right\rfloor$, as we chose the function $v\colon \A \to \R^{+}\setminus \{0\}$ such that for every $a \in \A$ $v(a)>4$, we get that $t \geq 1$. There are thus only finitely many possible pairs $(h,t)$. Finally, $b_i$ describes the index of the tile (starting from 0) on the last row of the $(h,t)$ block which contains the same vertical coordinate as the vertex corresponding to the $i$-th son of the $(a,w)$-tile and thus can take values in $[0;2^{h-1}(t+1)-1]$. As $s_i\in \{0,1\}$ we conclude that there are finitely many symbols in $\B$.

  \subsubsection*{Definition of the forbidden patterns $F$}

  All forbidden patterns in $F$ have supports which consist of three vertices $\{u,v,w\}$ such that $(u,v),(u,w)$ are edges, $L((u,v)) = \texttt{next}$ and $L((u,w)) = \ell$ for some $\ell$ appearing in the parent matching labels of the orbit graph. Denote by $$b^u = \left((a^u,z^u), (h^u,t^u), [(b_0^u,s_0^u),\dots,(b^u_{|w|-1},s^u_{|w|-1})]\right)$$ the symbol appearing in $u$ and do similarly for $v,w$. We say the pattern $p \colon \{u,v,w\} \to \B$ is in $F$ if and only if one of the following conditions does not hold:

	 \begin{enumerate}
	 	\item $a^w = (z^u)_{\ell+1}$;
	 	\item $h^u = h^v$;
	 	\item If $\ell < |z^u|-1$, then $2(b^u_{\ell+1}-b^u_{\ell})+{s^u_{\ell+1}-s^u_\ell} = t^w$.
	 	\item If $\ell = |z^u|-1$, then $2(2^{h^u-1}t^u+b^v_{0}-b^u_{|z^u|-1})+{s^v_{0}-s^u_{|z^u|-1}} = t^w$.
	 \end{enumerate}

\begin{figure}[!ht]
  	\begin{center}
  			\begin{tikzpicture}[scale=0.4]
	\foreach \Point in {(0,0), (6,0), (12,0), (18,0), (24,0)}{
	\begin{scope}[color=black!50,shift={\Point}]
	\draw (0,0) -- (0,3) to [controls=+(45:2) and +(135:2)] (6,3) -- (6,0) to [controls=+(135:1) and +(45:1)] (3,0) to [controls=+(135:1) and +(45:1)] (0,0);
	\end{scope}
	}
	\foreach \Point in {(0,0), (6,0), (12,0), (18,0)}{
	\begin{scope}[color=black!50,shift={\Point}]
	\draw[fill=red!15] (0,0) -- (0,3) to [controls=+(45:2) and +(135:2)] (6,3) -- (6,0) to [controls=+(135:1) and +(45:1)] (3,0) to [controls=+(135:1) and +(45:1)] (0,0);
	\end{scope}
	}
	\foreach \Point in {(0,-3), (3,-3), (6,-3), (9,-3), (12,-3), (15,-3), (18,-3), (21,-3), (24,-3), (27,-3)}{
	\begin{scope}[color=black!50,shift={\Point},scale=0.5]
	\draw (0,0) -- (0,6) to [controls=+(45:2) and +(135:2)] (6,6) -- (6,0) to [controls=+(135:1) and +(45:1)] (3,0) to [controls=+(135:1) and +(45:1)] (0,0) -- cycle;
	\end{scope}
	}
	\foreach \Point in {(0,-3), (21,-3),(3,-3), (6,-3), (9,-3), (12,-3), (15,-3), (18,-3)}{
	\begin{scope}[color=black!50,shift={\Point},scale=0.5]
	\draw[fill=red!15] (0,0) -- (0,6) to [controls=+(45:2) and +(135:2)] (6,6) -- (6,0) to [controls=+(135:1) and +(45:1)] (3,0) to [controls=+(135:1) and +(45:1)] (0,0) -- cycle;
	\end{scope}
	}
	\foreach \Point in {(0,-6), (1.5,-6), (3,-6), (4.5,-6), (6,-6), (7.5,-6), (9,-6), (12,-6), (13.5,-6), (15,-6), (16.5,-6), (18,-6), (19.5,-6), (21,-6), (22.5,-6), (24,-6), (25.5,-6), (27,-6), (28.5,-6)}{
	\begin{scope}[color=black!50,shift={\Point},scale=0.25]
	\draw (0,0) -- (0,12) to [controls=+(45:2) and +(135:2)] (6,12) -- (6,0) to [controls=+(135:1) and +(45:1)] (3,0) to [controls=+(135:1) and +(45:1)] (0,0) -- cycle;
	\end{scope}
	}
	\foreach \Point in {(0,-6), (1.5,-6), (22.5,-6), (3,-6), (4.5,-6), (6,-6), (7.5,-6), (9,-6), (10.5,-6), (12,-6), (13.5,-6), (15,-6), (16.5,-6), (18,-6), (19.5,-6), (21,-6)}{
	\begin{scope}[color=black!50,shift={\Point},scale=0.25]
	\draw[fill=red!15] (0,0) -- (0,12) to [controls=+(45:2) and +(135:2)] (6,12) -- (6,0) to [controls=+(135:1) and +(45:1)] (3,0) to [controls=+(135:1) and +(45:1)] (0,0) -- cycle;
	\end{scope}
	}
	\foreach \Point in {(0,-9), (0.75,-9), (1.5,-9), (2.25,-9), (3,-9), (3.75,-9), (4.5,-9), (5.25,-9), (6,-9), (6.75,-9), (7.5,-9), (8.25,-9), (9,-9), (9.75,-9), (10.5,-9), (11.25,-9), (12,-9), (12.75,-9), (13.5,-9), (14.25,-9), (15,-9), (15.75,-9), (16.5,-9), (17.25,-9), (18,-9), (18.75,-9), (19.5,-9), (20.25,-9), (21,-9), (21.75,-9), (22.5,-9), (23.25,-9), (24,-9), (24.75,-9), (25.5,-9), (26.25,-9), (27,-9), (27.75,-9), (28.5,-9), (29.25,-9)}{
	\begin{scope}[color=black!50,shift={\Point},scale=0.125]
	\draw (0,0) -- (0,24) to [controls=+(45:2) and +(135:2)] (6,24) -- (6,0) to [controls=+(135:1) and +(45:1)] (3,0) to [controls=+(135:1) and +(45:1)] (0,0) -- cycle;
	\end{scope}
	}
	\foreach \Point in {(0,-9), (0.75,-9), (1.5,-9), (2.25,-9), (3,-9), (3.75,-9), (4.5,-9), (5.25,-9), (6,-9), (6.75,-9), (7.5,-9), (8.25,-9), (9,-9), (9.75,-9), (10.5,-9), (11.25,-9), (12,-9), (12.75,-9), (13.5,-9), (14.25,-9), (15,-9), (15.75,-9), (16.5,-9), (17.25,-9), (18,-9), (18.75,-9), (19.5,-9), (20.25,-9), (21,-9), (21.75,-9), (22.5,-9), (23.25,-9)}{
	\begin{scope}[color=black!50,shift={\Point},scale=0.125]
	\draw[fill=red!15] (0,0) -- (0,24) to [controls=+(45:2) and +(135:2)] (6,24) -- (6,0) to [controls=+(135:1) and +(45:1)] (3,0) to [controls=+(135:1) and +(45:1)] (0,0) -- cycle;
	\end{scope}
	}
	\foreach \Point in {(0,-12), (0.375,-12), (0.75,-12), (1.125,-12), (1.5,-12), (1.875,-12), (2.25,-12), (2.625,-12), (3,-12), (3.375,-12), (3.75,-12), (4.125,-12), (4.5,-12), (4.875,-12), (5.25,-12), (5.625,-12), (6,-12), (6.375,-12), (6.75,-12), (7.125,-12), (7.5,-12), (7.875,-12), (8.25,-12), (8.625,-12), (9,-12), (9.375,-12), (9.75,-12), (10.125,-12), (10.5,-12), (10.875,-12), (11.25,-12), (11.625,-12), (12,-12), (12.375,-12), (12.75,-12), (13.125,-12), (13.5,-12), (13.875,-12), (14.25,-12), (14.625,-12), (15,-12), (15.375,-12), (15.75,-12), (16.125,-12), (16.5,-12), (16.875,-12), (17.25,-12), (17.625,-12), (18,-12), (18.375,-12), (18.75,-12), (19.125,-12), (19.5,-12), (19.875,-12), (20.25,-12), (20.625,-12), (21,-12), (21.375,-12), (21.75,-12), (22.125,-12), (22.5,-12), (22.875,-12), (23.25,-12), (23.625,-12), (24,-12), (24.375,-12), (24.75,-12), (25.125,-12), (25.5,-12), (25.875,-12), (26.25,-12), (26.625,-12), (27,-12), (27.375,-12), (27.75,-12), (28.125,-12), (28.5,-12), (28.875,-12), (29.25,-12), (29.
625,-12)}{
	\begin{scope}[color=black!30,shift={\Point},scale=0.0625]
	\draw (0,0) -- (0,48) to [controls=+(45:2) and +(135:2)] (6,48) -- (6,0) to [controls=+(135:1) and +(45:1)] (3,0) to [controls=+(135:1) and +(45:1)] (0,0) -- cycle;
	\end{scope}
	}

	\draw[thick,opacity=0.3,fill=black!25] (1.8,-11.5) -- (1.8,2.8) to [controls=+(45:2) and +(135:2)] (25.3,2.8) -- (25.3,-11.5) to [controls=+(135:1.5) and +(45:1.5)] (13.4,-11.5) to [controls=+(135:1.3) and +(45:1.3)] (7,-11.5) to [controls=+(135:1) and +(45:1)] (1.8,-11.5);

	\draw node at (13.55, -4.35) [opacity=0.6] {{\scalebox{1.6}{$u$}}};

	\draw[thick,opacity=0.3,fill=black!25] (7,-14.5) -- (7,-11.5) to [controls=+(45:1.3) and +(135:1.3)] (13.4,-11.5) -- (13.4,-14.5) to [controls=+(135:0.7) and +(45:0.7)] (11,-14.5) to [controls=+(135:0.9) and +(45:0.9)] (8.5,-14.5) to [controls=+(135:0.5) and +(45:0.5)] (7,-14.5);

	\draw node at (10.2, -13) [opacity=0.6] {{\scalebox{1}{$w$}}};

	\draw[->] (1.825, -7.5) to (1.75,-9);
	\draw[->] (7.125, -7.5) to (7.0,-9);
	\draw[->] (13.125, -7.5) to (13.25,-9);
	\draw[fill=white] (7.125-0.3,-7.5-0.3) rectangle + (0.6,0.6);
	\draw[fill=white] (1.825-0.3,-7.5-0.3) rectangle + (0.6,0.6);
	\draw[fill=white] (13.125-0.3,-7.5-0.3) rectangle + (0.6,0.6);
	\draw node at (7.125,-7.5) {{\scalebox{0.6}{$9$}}};
	\draw node at (1.825,-7.5) {{\scalebox{0.6}{$2$}}};
	\draw node at (13.125,-7.5) {{\scalebox{0.6}{$17$}}};

	\draw node at (12,5.5) {$\boldsymbol{t}$};
	\draw[decorate,decoration={brace,amplitude=6pt,raise=0.1cm}] (0.1,4.25) -- (23.9,4.25);
	\draw node at (-2,-3) {$\boldsymbol{h}$};
	\draw[decorate,decoration={brace,amplitude=6pt,raise=0.1cm}] (-0.5,-9) -- (-0.5,3);
	\draw node at (9.8,-9) {$\boldsymbol{17}$};
	\draw[decorate,decoration={brace,amplitude=6pt,raise=0.1cm}] (6.85,-10.5) -- (13,-10.5);

	\draw node at (9.8,-5.5) {$\boldsymbol{8}$};
	\draw[decorate,decoration={brace,amplitude=6pt,raise=0.1cm}] (6.75,-7) -- (12.75,-7);
	\end{tikzpicture}
  	\end{center}
  	\caption{Illustration of the third item in the definition of the set of forbidden patterns $F$: there are $8 = b_2 -b_1$ tiles in the bottom row of the top tile and $s_2 = 1, s_1 = 0$. Thus there must be $2(b_2-b_1)+(s_1-s_2)=17$ tiles on the top row of the pattern coded by the tile appearing below $u$ (which is called $w$). To make the picture smaller, the bottom tile is drawn shorter than it should be. If the rightmost tile is considered, we must add the number of tiles $2^{h-1}t$ to $b_0$ of the rightmost tile for the formula to add up.}
  	\label{figure.definition_rules}
  \end{figure}
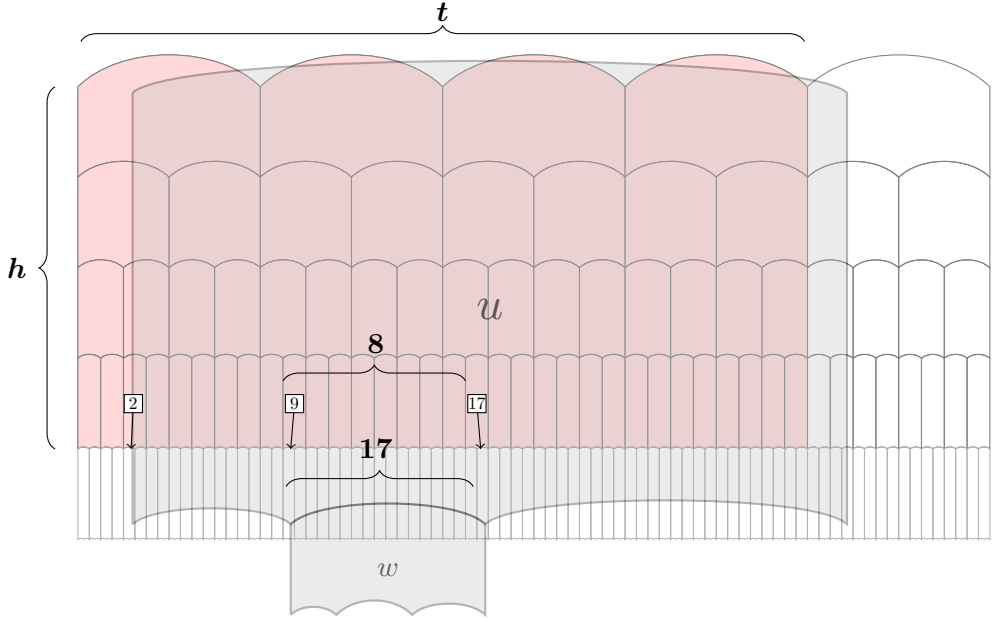

	 The first rule says that if the rule $(a,z_1z_2\dots z_k)$ appears in a vertex, then a rule starting with $z_{\ell+1}$ should appear in the son labeled with $\ell$. The second rule says any two symbols that lie in a row of the orbit graph have the same height $h$. The third and fourth rules say that if $w$ is the $\ell$-th son of $u$, then the length $t^w$ of the block appearing at $w$ must be consistent with the bottom row of the block appearing at $u$ (see~\cref{figure.definition_rules}).

	Consider an orbit $\Omega = \{(w^i,P_i)\}_{i \in \Z}$ of $(\A,R)$ and its associated orbit graph $\Gamma_{\Omega}$. We define $Y \subset \B^{\Gamma_{\Omega}}$ as the subshift consisting of all colorings of $\Gamma_{\Omega}$ by symbols of $\B$ where the patterns from $F$ do not appear.
  \begin{lemma}\label{prop_Ynonempty}
  	For every orbit $\Omega = \{(w^i,P_i)\}_{i \in \Z}$ of $(\A,R)$ the subshift $Y\subset \B^{\Gamma_{\Omega}}$ is non-empty.
  	\end{lemma}

  \begin{proof}
   By~\cref{proposition.orbits_implies_tiling} there exists a tiling $\Psi_{\Omega}\colon \Z^2 \to \R^2$ for $\Omega$. Similarly, fixing an orbit $\Xi$ of $(\{0\},0\mapsto 00)$ there is a tiling $\Psi_{\Xi}\colon \Z^2 \to \R^2$ for $\Xi$.

   We claim that for every $u =(i,j) \in \Z^2$, there is $u^* = (i^*,j^*) \in \Z^2$ such that if $\Psi_{\Xi}(u^*)= (x,y)$ then $\Psi_{\Omega}(u) = (x+\widetilde{x},y-\widetilde{y})$ for some $\widetilde{x} \in [0,2\cdot e^y[$ and $\widetilde{y} \in [0,\log(2)[$. Indeed, by definition of tiling, note that if $\Psi_{\Xi}(i_1,j_1)=(x_1,y_1)$ and $\Psi_{\Xi}(i_2,j_2) =(x_2,y_2)$ then $y_2= y_1-(i_2-i_1)\log(2)$. Therefore if we let $\Psi_{\Omega}(u) = (a,b)$ we can first find $i^*$ such that $\Psi_{\Xi}(i^*,k) = (\cdot,y) \in \R \times [b,b+\log(2)[$ for all $k \in \Z$. Furthermore, if $\Psi_{\Xi}(i^*,j_1) = (x_1,y_1)$ and $\Psi_{\Xi}(i^*,j_2) = (x_2,y_2)$ we have $y_1 = y_2 = y$ and $x_2-x_1 =2(j_2-j_1)e^{y}$. Therefore we can find $j^*$ such that $\Psi_{\Xi}(i^*,j^*) = (x,y)\in [a,a+2e^y[ \times [b,b+\log(2)[$. Hence, letting $u^* = (i^*,j^*)$ we have $\Psi_{\Omega}(u)-\Psi_{\Xi}(u*) =(\widetilde{x},-\widetilde{y})$ as required.

   Let us define a configuration $c\colon \Z^2 \to \B$. Let $c(i,j)$ be the symbol of $\B$ associated to the $((a^i)_j, a^{i+1}|_{[\Delta_{i+1}(j);\Delta_{i+1}(j+1)-1]})$-tile at position $\Psi_{\Xi}((i,j)^*)+(\widetilde{x},-\widetilde{y})$ as described in the definition of $\B$. We claim that $c \in Y$. We need to show that $c$ does not contain any forbidden pattern from $F$, i.e. that any pattern with one of the supports defining $F$ satisfies the four conditions described above.

   Let $u,v,w \in \Z^2$ such that $L((u,v))=\texttt{next}$ and $L((u,w))=\ell$ and consider the pattern $c|_{\{u,v,w\}}$.
   Denote $(\bar{x},\bar{y})=\Psi_{\Omega}(u)$, $(x,y)= \Psi_{\Omega}(u^*)$, $(\widetilde{x},-\widetilde{y})=(\bar{x}-x,\bar{y}-y)$ and the production rule appearing at $u$ be $(a,z_1\dots z_k)$ and thus $0 \leq \ell < k$.
   By definition of $c$ we have that $a^w = z_{\ell+1} = (z^u)_{\ell+1}$ and hence the first rule of $F$ holds. By definition of tiling we have that $\Psi_{\Omega}(v)=(\bar{x}+v(a)e^{\bar{y}},\bar{y})$ and hence if we have $u^* = (i_1^*,j_1^*)$ and $v^* = (i_2^*,j_2^*)$ then $i_2^* = i_1^*$. This implies that $h^u =  h^v$ and hence the second rule of $F$ holds.
   To simplify the notations for the remainder of the proof, we drop the superscripts for $u$, that is, we denote $h = h^u$, $b^u_i = b_i$ and $s^u_i = s_i$ and maintain the superscripts for the $v$ and $w$.

    On the one hand, by the Euclidean division algorithm, we have that for any $0 \leq \ell < k$:
	$$2b_{\ell} + s_{\ell} = \left\lfloor\frac{\widetilde{x}+e^{y-\widetilde{y}-\log(\lambda)}\sum_{r=1}^{\ell}v(z_r)}{2e^{y-h\log(2)}}\right\rfloor.$$
	Also, as $v(a)= e^{-\log(\lambda)\sum_{r=1}^{k}v(z_r)}$ we have that:
	\begin{align*}
	2(2^{h-1}t+b^v_{0})+s^v_{0}& = 2^ht + \left\lfloor\frac{\bar{x}+v(a)e^{y-\widetilde{y}}- x-2t^ue^{y}}{2e^{y-h\log(2)}}\right\rfloor\\
& = \left\lfloor\frac{\widetilde{x}+e^{y-\widetilde{y}}\sum_{r=1}^{k}v(z_r)e^{-\log(\lambda)}}{2e^{y-h\log(2)}}\right\rfloor\\
	\end{align*} and thus we shall denote $2(2^{h-1}t+b^v_{0})+s^v_{0}$ simply by $2b_{k}+s_{k}$ as it has the same expression as the numbers above.

On the other hand, we have $\Psi_{\Omega}(w) = (\bar{x}+e^{\bar{y}-\log(\lambda)}\sum_{k=1}^{\ell}v(z_k),\bar{y}-\log(\lambda))$. It is easy to verify that $\Psi_{\Xi}(w^*)=(x+2e^{y-h\log(2)}(2b_{\ell}+s_{\ell}),y-h\log(2))$. It follows that
\[\Psi_{\Omega}(w)-\Psi_{\Xi}(w^*) = \left(\widetilde{x}+e^{\bar{y}-\log(\lambda)}\sum_{k=1}^{\ell}v(z_k)-2e^{y-h\log(2)}(2b_{\ell}+s_{\ell}), ~ -\left(\widetilde{y}+\log(\lambda)-h\log(2)\right)\right) \]
   and thus

   \begin{align*}
   t^w & = \left\lfloor \frac{\widetilde{x}+e^{\bar{y}-\log(\lambda)}\sum_{k=1}^{\ell}v(z_k)-2e^{y-h\log(2)}(2b_{\ell}+s_{\ell}) + v(z_{\ell+1})e^{\bar{y}-\log(\lambda)} }{2e^{y-h\log(2)}}       \right\rfloor\\
   & =  \left\lfloor \frac{\widetilde{x}+e^{\bar{y}-\log(\lambda)}\sum_{k=1}^{\ell+1}v(z_k)}{2e^{y-h\log(2)}}       \right\rfloor - (2b_{\ell}+s_{\ell})  \\
   & =  (2b_{\ell+1}+s_{\ell+1}) - (2b_{\ell}+s_{\ell}) \\
   & = 2(b_{\ell+1} - b_{\ell})+s_{\ell+1}-s_{\ell}.   \\
 \end{align*} Therefore, conditions $3$ and $4$ are also satisfied, which means that $c|_{\{u,v,w\}} \notin F$. It follows that $c \in Y$ and hence $Y$ is non-empty. \end{proof}

  \subsection{Simulation of orbits of $(\{0\},0\mapsto 00)$ on $(\A,R)$.}\label{subsection.simulation}

\begin{figure}[!ht]
  	\begin{center}
  			\begin{tikzpicture}

	\begin{scope}[scale=0.4, shift={(0,0)}]
	\foreach \Point in {(0,0), (6,0), (12,0)}{
	\begin{scope}[shift={\Point}]
	\draw (0,0) -- (0,3);
	\draw (0,3) to [controls=+(45:2) and +(135:2)] (6,3);
	\draw (6,3) -- (6,0);
	\draw (0,0) to [controls=+(45:1) and +(135:1)] (3,0);
	\draw (3,0) to [controls=+(45:1) and +(135:1)] (6,0);
	\end{scope}
	}
	\foreach \Point in {(0,-3), (3,-3), (6,-3), (9,-3), (12,-3), (15,-3)}{
	\begin{scope}[shift={\Point},scale=0.5]
	\draw (0,0) -- (0,6);
	\draw (0,6) to [controls=+(45:2) and +(135:2)] (6,6);
	\draw (6,6) -- (6,0);
	\draw (0,0) to [controls=+(45:1) and +(135:1)] (3,0);
	\draw (3,0) to [controls=+(45:1) and +(135:1)] (6,0);
	\end{scope}
	}
	\foreach \Point in {(0,-6), (1.5,-6), (3,-6), (4.5,-6), (6,-6), (7.5,-6), (9,-6), (10.5,-6), (12,-6), (13.5,-6), (15,-6), (16.5,-6)}{
	\begin{scope}[shift={\Point},scale=0.25]
	\draw (0,0) -- (0,12);
	\draw (0,12) to [controls=+(45:2) and +(135:2)] (6,12);
	\draw (6,12) -- (6,0);
	\draw (0,0) to [controls=+(45:1) and +(135:1)] (3,0);
	\draw (3,0) to [controls=+(45:1) and +(135:1)] (6,0);
	\end{scope}
	}
	\end{scope}

	\begin{scope}[scale=0.4, shift={(19.5,-1.5)}]
	\foreach \i in {0,6,12}{
		\node[scale=1] at (\i+3,3) {\textbullet};
		\draw[->, thick] (\i+3,3) -- (\i+1.5,0.3);

		\draw[->, thick] (\i+3,3) -- (\i+4.5,0.3);
		\node at (\i+1.8,1.5) {$0$};
		\node at (\i+4.2,1.5) {$1$};
	}
	\foreach \i in {0,6}{
		\draw[thin,->,>=latex,dashed] (\i+3,3) -- (\i+8.7,3);
	}
	\foreach \i in {0,3,6,9,12,15}{
		\node[scale=1] at (\i+1.5,0) {\textbullet};
		\draw[->, thick] (\i+1.5,0) -- (\i+0.75,-2.7);
		\draw[->, thick] (\i+1.5,0) -- (\i+2.25,-2.7);
		\node at (\i+0.75,-1.5) {$0$};
		\node at (\i+2.25,-1.5) {$1$};
	}
	\foreach \i in {0,3,6,9,12}{
		\draw[thin,->,>=latex,dashed] (\i+1.5,0) -- (\i+4.2,0);
	}
	\foreach \i in {0,1.5,3,4.5,6,7.5,9,10.5,12,13.5,15,16.5}{
		\node[scale=1] at (\i+0.75,-3) {\textbullet};
	}
	\foreach \i in {0,1.5,3,4.5,6,7.5,9,10.5,12,13.5,15}{
		\draw[thin,->,>=latex,dashed] (\i+0.75,-3) -- (\i+2,-3);
	}
	\end{scope}

	\end{tikzpicture}
  	\end{center}
  	\caption{A $(3,3)$-block and its associated $\Gamma_{(3,3)}$ graph. The $\texttt{next}$ edges are shown as dashed lines.}
  	\label{figure.correspondence_graph_tile}
  \end{figure}
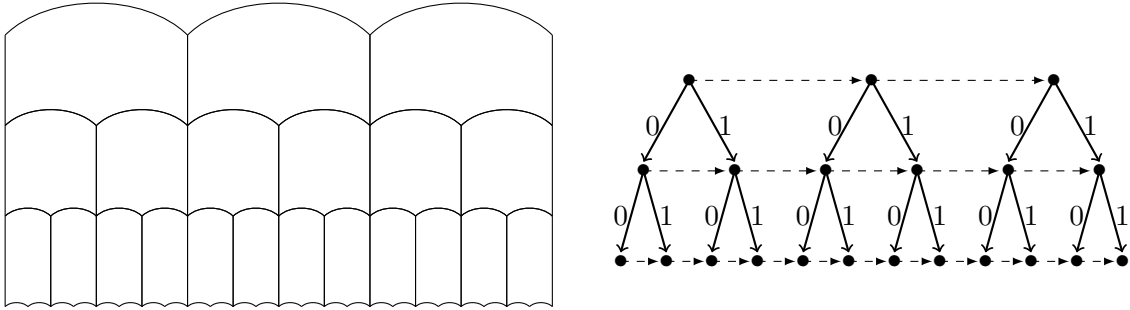

  For every $b \in \B$ we can associate a finite graph $\Gamma_{b} = (V_b,E_b,L_b)$ which appears as an induced subgraph on any orbit graph of $(\{0\},0\mapsto 00)$ as follows: Let $(h,t)$ be the second coordinate of $b$, the vertex set is $V_b = \{ (i,j) \mid i \in [0;h-1], j\in [0,t2^{i}-1] \}$ and the edges have labels given by $L_b(((i,j),(i,j+1))) = \texttt{next}$ for each $i$ and $j<t2^{i}-1$ and $L_b(((i-1,\lfloor \frac{j}{2} \rfloor ) ,(i,j) )) = j \mod{2}$ for every $i \geq 1$. See~\cref{figure.correspondence_graph_tile}.

  \begin{remark*}
  	For every $b \in \B$ the associated graph $\Gamma_b$ is non-empty. As $\lambda >2$ and $v(a)>4$ for every $a \in \A$ we have that the numbers $(h,t)$ associated to every $b \in \B$ are both larger than $1$.
  \end{remark*}

  More generally, given a finite connected subset $S \subset \Gamma$ and a pattern $q\colon S \to \B$ which appears in some configuration of $Y$ we can associate a finite subgraph $\Gamma_q$ by pasting together the graphs $(\Gamma_{q(s)})_{s \in S}$ in the following way:

  \begin{enumerate}
  	\item Whenever $u,v \in S$ are connected by a $\texttt{next}$ edge from $u$ to $v$, we connect $\Gamma_{q(u)}$ to $\Gamma_{q(v)}$ by joining the rightmost vertices of $\Gamma_{q(u)}$ to the leftmost vertices of $\Gamma_{q(v)}$ with $\texttt{next}$ edges. More precisely, if $q(u)$ codes an $(h,t)$-block, then for every $i \in [0;h-1]$ we connect the vertex $(i,t2^{i}-1)$ of $\Gamma_{q(u)}$ to $(i,0)$ of $\Gamma_{q(v)}$ by a $\texttt{next}$ edge.
  	\item Whenever $u,w \in S$ are connected by an edge with label $i$, we look at the coordinate $(b_i,s_i)$ of $q(u)$ and connect the left-top vertex of $\Gamma_{q(w)}$ to $b_i$-th vertex on from the left on the bottom row of $\Gamma_{q(u)}$ using an $s_i$-edge and then connect all vertices on the top row of $\Gamma_{q(w)}$ to the bottom row of $\Gamma_{q(u)}$ alternating $0-1$ edges. More precisely, if $q(u)$ codes an $(h,t)$-block then for each $j$ we connect vertex $(h-1,b_i+\lfloor\frac{s_i+j}{2}\rfloor)$ of the bottom row of $\Gamma_{q(u)}$ to vertex $(0,j)$ from the top row of $\Gamma_{q(w)}$ with a label $s_i+j \mod{2}$. If $(h-1,b_i+\lfloor\frac{s_i+j}{2}\rfloor)$ does not appear in the bottom row of $\Gamma_{q(u)}$ and $u$ is connected to some vertex $v$ by a $\texttt{next}$ label, then the vertex $(h-1,b_i+\lfloor\frac{s_i+j}{2}\rfloor)$ gets replaced by vertex $(h-1,b_i+\lfloor\frac{s_i+j}{2}\rfloor - 2^{h-1}t)$ of $\Gamma_{q(v)}$.
  \end{enumerate}

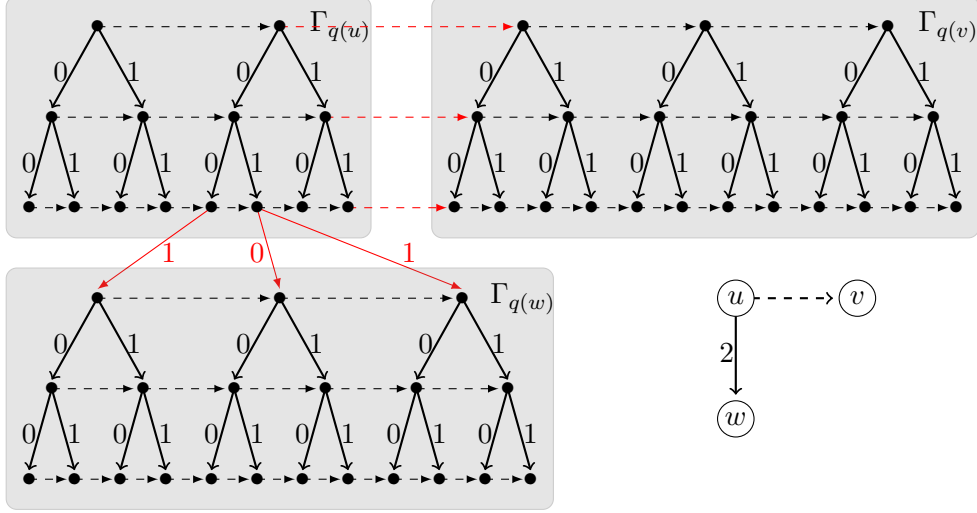
\begin{figure}[!ht]
	\begin{center}
			\begin{tikzpicture}

	\begin{scope}[scale=0.4, shift={(14,0)}]
	
	\draw[fill = gray, opacity = 0.2, rounded corners] (0,-4) rectangle (18,4);
	\node at (17,3) {$\Gamma_{q(v)}$};
	\foreach \i in {0,6,12}{
		\node[scale=1] at (\i+3,3) {\textbullet};
		\draw[->, thick] (\i+3,3) -- (\i+1.5,0.3);
		
		\draw[->, thick] (\i+3,3) -- (\i+4.5,0.3);
		\node at (\i+1.8,1.5) {$0$};
		\node at (\i+4.2,1.5) {$1$};
	}
	\foreach \i in {0,6}{
		\draw[thin,->,>=latex,dashed] (\i+3,3) -- (\i+8.7,3);
	}
	\foreach \i in {0,3,6,9,12,15}{
		\node[scale=1] at (\i+1.5,0) {\textbullet};
		\draw[->, thick] (\i+1.5,0) -- (\i+0.75,-2.7);
		\draw[->, thick] (\i+1.5,0) -- (\i+2.25,-2.7);
		\node at (\i+0.75,-1.5) {$0$};
		\node at (\i+2.25,-1.5) {$1$};
	}
	\foreach \i in {0,3,6,9,12}{
		\draw[thin,->,>=latex,dashed] (\i+1.5,0) -- (\i+4.2,0);
	}
	\foreach \i in {0,1.5,3,4.5,6,7.5,9,10.5,12,13.5,15,16.5}{
		\node[scale=1] at (\i+0.75,-3) {\textbullet};
	}
	\foreach \i in {0,1.5,3,4.5,6,7.5,9,10.5,12,13.5,15}{
		\draw[thin,->,>=latex,dashed] (\i+0.75,-3) -- (\i+2,-3);
	}
	\end{scope}
	
		\begin{scope}[scale=0.4, shift={(0,0)}]
		
		\draw[fill = gray, opacity = 0.2, rounded corners] (0,-4) rectangle (12,4);
		\node at (11,3) {$\Gamma_{q(u)}$};
		
		\draw[thin,->,>=latex,dashed, red] (6+3,3) -- (6+8.7+2,3);
		\draw[thin,->,>=latex,dashed, red] (9+1.5,0) -- (9+4.2+2,0);
		\draw[thin,->,>=latex,dashed, red] (10.5+0.75,-3) -- (10.5+2+2,-3);
		\draw[thin,->,>=latex, red] (6+0.75,-3) -- (3,-6+0.3);
		\draw[thin,->,>=latex, red] (7.5+0.75,-3) -- (9,-6+0.3);
		\draw[thin,->,>=latex, red] (7.5+0.75,-3) -- (15,-6+0.3);
		\node at (4.6+0.75,-4.5) {\color{red}$1$};
		\node at (7.5+0.75,-4.5) {\color{red}$0$};
		\node at (12.5+0.75,-4.5) {\color{red}$1$};
		\foreach \i in {0,6}{
			\node[scale=1] at (\i+3,3) {\textbullet};
			\draw[->, thick] (\i+3,3) -- (\i+1.5,0.3);
			
			\draw[->, thick] (\i+3,3) -- (\i+4.5,0.3);
			\node at (\i+1.8,1.5) {$0$};
			\node at (\i+4.2,1.5) {$1$};
		}
		\foreach \i in {0}{
			\draw[thin,->,>=latex,dashed] (\i+3,3) -- (\i+8.7,3);
		}
		\foreach \i in {0,3,6,9}{
			\node[scale=1] at (\i+1.5,0) {\textbullet};
			\draw[->, thick] (\i+1.5,0) -- (\i+0.75,-2.7);
			\draw[->, thick] (\i+1.5,0) -- (\i+2.25,-2.7);
			\node at (\i+0.75,-1.5) {$0$};
			\node at (\i+2.25,-1.5) {$1$};
		}
		\foreach \i in {0,3,6}{
			\draw[thin,->,>=latex,dashed] (\i+1.5,0) -- (\i+4.2,0);
		}
		\foreach \i in {0,1.5,3,4.5,6,7.5,9,10.5}{
			\node[scale=1] at (\i+0.75,-3) {\textbullet};
		}
		\foreach \i in {0,1.5,3,4.5,6,7.5,9}{
			\draw[thin,->,>=latex,dashed] (\i+0.75,-3) -- (\i+2,-3);
		}
		\end{scope}
		
		\begin{scope}[scale=0.4, shift={(0,-9)}]
		
		\draw[fill = gray, opacity = 0.2, rounded corners] (0,-4) rectangle (18,4);
		\node at (17,3) {$\Gamma_{q(w)}$};
		
		\foreach \i in {0,6,12}{
			\node[scale=1] at (\i+3,3) {\textbullet};
			\draw[->, thick] (\i+3,3) -- (\i+1.5,0.3);
			
			\draw[->, thick] (\i+3,3) -- (\i+4.5,0.3);
			\node at (\i+1.8,1.5) {$0$};
			\node at (\i+4.2,1.5) {$1$};
		}
		\foreach \i in {0,6}{
			\draw[thin,->,>=latex,dashed] (\i+3,3) -- (\i+8.7,3);
		}
		\foreach \i in {0,3,6,9,12,15}{
			\node[scale=1] at (\i+1.5,0) {\textbullet};
			\draw[->, thick] (\i+1.5,0) -- (\i+0.75,-2.7);
			\draw[->, thick] (\i+1.5,0) -- (\i+2.25,-2.7);
			\node at (\i+0.75,-1.5) {$0$};
			\node at (\i+2.25,-1.5) {$1$};
		}
		\foreach \i in {0,3,6,9,12}{
			\draw[thin,->,>=latex,dashed] (\i+1.5,0) -- (\i+4.2,0);
		}
		\foreach \i in {0,1.5,3,4.5,6,7.5,9,10.5,12,13.5,15,16.5}{
			\node[scale=1] at (\i+0.75,-3) {\textbullet};
		}
		\foreach \i in {0,1.5,3,4.5,6,7.5,9,10.5,12,13.5,15}{
			\draw[thin,->,>=latex,dashed] (\i+0.75,-3) -- (\i+2,-3);
		}
		\end{scope}
		
		\begin{scope}[scale=0.4, shift={(24,-6)}]
		
		\draw[->, thick, dashed] (0,0) -- (3.2,0);
		\draw[->, thick] (0,0) -- (0,-3.2);
		\node at (-0.3,-1.8) {$2$};
		\draw[fill = white] (0,0) circle (0.6);
		\draw[fill = white] (4,0) circle (0.6);
		\draw[fill = white] (0,-4) circle (0.6);
		\node[scale=1] at (0,0) {$u$};
		\node[scale=1] at (4,0) {{$v$}};
		\node[scale=1] at (0,-4) {{$w$}};
		\end{scope}

	\end{tikzpicture}
	\end{center}
	\caption{The rules for pasting graphs.}
	\label{figure.metagraph}
\end{figure}

  \begin{remark*}
  	The pasting rules above are consistent because no patterns from $F$ appear in $q$. More precisely, if two vertices are connected by a $\texttt{next}$ edge the blocks they code have the same height by rule $2$ of $F$ and thus the first rule is coherent. If two vertices are connected by an $i$-edge then the sites where the graphs are pasted do not overlap and cover everything by rules $3$ and $4$ of $F$. We illustrate the pasting rules in~\cref{figure.metagraph}.
  \end{remark*}

  Let $\Sigma$ be a finite alphabet and $F_{\Sigma}$ a set of nearest neighbor forbidden patterns on the orbit graph of $(\{0\},0\mapsto 00)$ over the alphabet $\Sigma$. We define ${\B}_{\Sigma}$ as the set of pairs $(b,p_b)$ such that $b \in \B$ and $p_b\colon \Gamma_{b}\to \Sigma$ is a pattern. Also, for a pattern ${p}$ on $\Gamma_{\Omega}$ with alphabet ${B}_{\Sigma}$ denote by $\pi_{\B}({p})$ the restriction to the first coordinate of ${\B}_{\Sigma}$. Also denote by $q({p}) \colon \Gamma_{\pi_{\B}({p})} \to \Sigma$ the pattern over $(\{0\},0\mapsto 00)$ whose support is the graph $\Gamma_{\pi_{\B}({p})}$ and is obtained by pasting together the corresponding patterns $p_b$ on the second coordinate of ${B}_{\Sigma}$.

  Define ${F}_{\B,{\Sigma}}$ as the set of all patterns $p$ over the alphabet ${\B}_{\Sigma}$ which have supports which consist in three vertices $\{u,v,w\}$ in $\Gamma_{\Omega}$ such that $(u,v),(u,w)$ are edges, $L((u,v)) = \texttt{next}$ and $L((u,w)) = \ell$ for some $\ell$ appearing in the parent matching labels of the orbit graph $\Gamma_{\Omega}$, and that satisfy one of the following two properties:

  \begin{enumerate}
  	\item The pattern $\pi_{\B}({p})$ obtained by restricting ${p}$ to the first coordinate of ${\B}_{\Sigma}$ is in $F$;
  	\item The pattern $q({p})$ obtained by pasting the patterns of ${p}$ described by the second coordinate of ${\B}_{\Sigma}$ contains a forbidden pattern from $F_{\Sigma}$.
  \end{enumerate}

   Clearly ${F}_{\B,{\Sigma}}$ has finitely many patterns (up to label preserving graph isomorphism). For any orbit $\Omega$ of $(\A,R)$ we define the subshift of finite type ${Y}_{\Sigma} \subset ({B}_{\Sigma})^{\Gamma_{\Omega}}$ as the set of all colorings of $\Gamma_{\Omega}$ by ${B}_{\Sigma}$ where no pattern from $F_{\B,{\Sigma}}$ appears.

  \begin{lemma}\label{proposition_reduction}
  	Let $\Omega$ and $\Xi$ be orbits of $(\A,R)$ and $(\{0\},0\mapsto 00)$ respectively. Let $\Gamma_{\Omega}$, $\Gamma_{\Xi}$ be orbit graphs of $\Omega$ and $\Xi$ respectively. Let $X_{\Sigma}$ be the subshift on $\Gamma_{\Xi}$ with alphabet $\Sigma$ defined by the nearest neighbor forbidden patterns $F_{\Sigma}$ and let ${Y}_{\Sigma} \subset ({B}_{\Sigma})^{\Gamma_{\Omega}}$ be defined as above. Then ${Y}_{\Sigma} = \emptyset$ if and only if $X_{\Sigma} = \emptyset$.
  \end{lemma}

  \begin{proof}
  	Assume there exists $\widetilde{y}\in {Y}_{\Sigma}$. Let $\widetilde{y}|_{n}$ be the restriction of $\widetilde{y}$ to the vertices $[-n,n]^2$ in $\Gamma_{\Omega}$. By definition of ${F}_{\B_{\Sigma}}$ the pattern $q(\widetilde{y}|_{n})$ does not contain any pattern from $F_{\Sigma}$. By a standard compactness argument, the sequence of patterns $(q(\widetilde{y}|_{n}))_{n \in \N}$ subconverges to a configuration $x\in \Sigma^{\Gamma_{\Xi}}$ which does not contain any pattern from $F_{\Sigma}$ and thus $x \in X_{\Sigma} \neq \emptyset$.

  	Conversely, let $x \in X_{\Sigma}$. By~\cref{prop_Ynonempty} there exists a configuration $y \in Y$. By identifying for each vertex $v \in \Gamma_{\Omega}$ the graphs $\Gamma_{y(v)}$ as a partition of the vertices of $\Gamma_{\Xi}$, we can construct a second coordinate $p_{x,y,v} = x|_{\Gamma_{y(v)}}$ which satisfies the second rule of ${F}_{\B,\Sigma}$. By definition $\widetilde{y} = (y(v),p_{x,y,v})$ is in ${Y}_{\Sigma}$ which is thus non-empty.
  \end{proof}

  \begin{remark*}
  	The alphabet ${B}_{\Sigma}$ and the set of forbidden patterns ${F}_{\B,\Sigma}$ which define ${Y}_{\Sigma}$ only depend upon $\Sigma$, $F_{\Sigma}$ and the substitution $(\A,R)$. It does not depend upon the choice of orbit $\Omega$ of $(\A,R)$.
  \end{remark*}

\begin{theorem}
	\label{th:DP_orbit}
  \begin{sloppypar}
	The domino problem is undecidable on any orbit graph of a non-deterministic substitution with an expanding eigenvalue.
  \end{sloppypar}
\end{theorem}

	\begin{proof}
	 For clarity, let us first assume that the expanding eigenvalue $\lambda$ associated to $(\A,R)$ satisfies $\lambda >2$. Let $\Sigma$ and $F_{\Sigma}$ be respectively an alphabet and a nearest neighbor set of forbidden patterns for an orbit graph $\Gamma_{\Xi}$ of an orbit $\Xi$ of $(\{0\},0\mapsto 00)$ which define a nearest neighbor SFT $X_{\Sigma}$. By~\cref{proposition_reduction} we know that $X_{\Sigma} = \emptyset$ if and only if ${Y}_{\Sigma} = \emptyset$. Furthermore, we claim that the alphabet and set of forbidden patterns which define ${Y}_{\Sigma}$ can be constructed effectively from $\Sigma$ and $F_{\Sigma}$. Indeed, the subshift $Y$ does not depend upon $\Sigma$ and thus its alphabet $\B$ and forbidden patterns $F$ can be hard-coded in the algorithm. It is easy to see that from $\B$ one can effectively construct the alphabet ${B}_{\Sigma}$ and the forbidden patterns ${F}_{\B,\Sigma}$ which define~${Y}_{\Sigma}$.

	 These two facts together show that if $\texttt{DP}(\Gamma_{\Omega})$ is decidable and $\lambda >2$, then so is $\texttt{DP}(\Gamma_{\Xi})$. Using the result of Kari (\cref{teorema_jarkko}) we have that $\texttt{DP}(\Gamma_{\Xi})$ is undecidable, hence $\texttt{DP}(\Gamma_{\Omega})$ is also undecidable.

	 We shall now deal with the general case where $1 <\lambda \leq 2$. For an integer $m \geq 1$ we define the relation $R^m$ recursively by:
	 \begin{itemize}
	 	\item $R^1 = R$.
	 	\item $R^{k+1}$ is the set of all pairs $(a,(c^1_{1}\dots c^{1}_{\ell_1})(c^2_{1}\dots c^{2}_{\ell_2}) \dots (c^1_{k}\dots c^{1}_{\ell_k}))$ in $\A\times\A^*$ for which there is a pair $(a,b_1\dots b_k) \in R^k$ such that $(b_i,c^i_{1}\dots c^{i}_{\ell_i}) \in R$ for each $i \in \{1,\hdots, k\}$.
	 \end{itemize}

	 In other words, $R^m$ is the set of all relations that can be obtained by starting with a symbol $a \in \A$ and replacing $m$ times each letter by the right hand side of a production rule of $R$. Let $n \in \N$ such that $\lambda^n > 2$ and note that the substitution $(\A,R^n)$ has the expanding eigenvalue $\lambda^n > 2$.

	Let $\Omega = \{(w^i,P_i)\}_{i \in \Z}$ be an orbit of $(\A,R)$. We have that for each $k \in \{0, \hdots, n-1\}$
  \[\Omega^{n,k} := \left\{\left(w^{in+k}, P_{in+k-(n-1)} \circ \dots \circ P_{in+k-1}\circ P_{in+k}\right)\right\}_{i \in \Z} \]
  is an orbit of $(\A,R^n)$. As before, let $\Sigma$ and $F_{\Sigma}$ be respectively an alphabet and a nearest neighbor set of forbidden patterns which define a nearest neighbor SFT $X_{\Sigma}$. Let ${Y}_{\Sigma}^{n,k}$ be the subshift $Y_{\Sigma}$ we constructed above, but now for the substitution $(\A,R^n)$ and orbit $\Omega^{n,k}$. Denote by ${\B}_{\Sigma}^n$ and ${F}_{\B,\Sigma}^n$ the alphabet and set of forbidden patterns of ${Y}_{\Sigma}^{n,k}$ respectively. By~\cref{proposition_reduction} we have that ${Y}_{\Sigma}^{n,k} = \emptyset$ if and only if $X_{\Sigma} = \emptyset$.

	We shall construct a subshift $Z$ on $\Gamma_{\Omega}$ which encodes a copy of ${{Y}_{\Sigma}^{k,n}}$ for each $k \in \{0, \hdots, n-1\}$. Consider again the alphabet ${\B}_{\Sigma}^n$. For every pattern $p \in {F}_{\B,\Sigma}^n$ with support $\{u,v,w\}$ such that $L((u,v))=\texttt{next}$ and $L((u,v)) = \ell$
  we define the set of patterns $F_{p}$ such that every $q \in F_{p}$ has support $\{v,u_1,u_2,\dots,u_{n}=w_0,w_1,\dots w_\ell\}$ such that $L((u_1,v))=\texttt{next}$, for every $i \in \{1, \hdots, n\}$, $L((u_i,u_{i+1}))=0$ and for every $j \in \{0,\hdots,\ell-1\}$, $L(w_i,w_{i+1})=\texttt{next}$ and every pattern $q$ in $F_{p}$ has the property that $q(u_1) = p(u)$, $q(v)=p(v)$ and $q(w_{\ell})=p_{\ell}$ (See~\cref{figure.reduction}).

	\begin{figure}[!ht]
		\begin{center}
				\begin{tikzpicture}

		\begin{scope}[scale=0.4, shift={(16,0)}]

		\draw[->, thick, dashed] (0,0) -- (3.2,0);
		\draw[->, thick] (0,0) -- (-1.6,-2.4);
		\draw[->, thick] (-2,-3) -- (-3.6,-5.4);
		\draw[->, thick] (-4.4,-6.6) -- (-5.6,-8.4);
		\node at (-0.3,-1.8) {$0$};
		\node at (-2.3,-4.8) {$0$};
		\node at (-4.3,-7.8) {$0$};
		\draw[fill = red!30] (0,0) circle (0.6);
		\draw[fill = blue!30] (4,0) circle (0.6);
		\draw[fill = white] (-2,-3) circle (0.6);
		\draw[fill = white] (-6,-9) circle (0.6);
		\node[scale=1] at (0,0) {$u_1$};
		\node[scale=1] at (4,0) {{$v$}};
		\node[scale=1] at (-2,-3) {{$u_2$}};
		\node[scale=1] at (-4,-6) {{$\dots$}};
		\node[scale=1] at (-6,-9) {{$u_n$}};
		\node[scale=1] at (-8,-9) {{$w_0=$}};
		\foreach \i in {-6,-2,2,6}{
		\draw[->, thick, dashed] (\i+0.8,-9) -- (\i+3.2,-9);}
		\foreach \i in {-2,2,10}{
		\draw[fill = white] (\i,-9) circle (0.6);}
		\draw[fill = green!30] (10,-9) circle (0.6);
		\node[scale=1] at (-2,-9) {{$w_1$}};
		\node[scale=1] at (2,-9) {{$w_2$}};
		\node[scale=1] at (6,-9) {{\tiny ${\dots}$}};
		\node[scale=1] at (10,-9) {{$w_{\ell}$}};

		\end{scope}

		\begin{scope}[scale=0.4, shift={(0,-2.5)}]

		\draw[->, thick, dashed] (0,0) -- (3.2,0);
		\draw[->, thick] (0,0) -- (0,-3.2);
		\node at (-0.3,-1.8) {$\ell$};
		\draw[fill = red!30] (0,0) circle (0.6);
		\draw[fill = blue!30] (4,0) circle (0.6);
		\draw[fill = green!30] (0,-4) circle (0.6);
		\node[scale=1] at (0,0) {$u$};
		\node[scale=1] at (4,0) {{$v$}};
		\node[scale=1] at (0,-4) {{$w$}};
		\end{scope}

	\end{tikzpicture}
		\end{center}
		\caption{On the left a pattern $p \in {F}_{\B,\Sigma}^n$. The corresponding patterns have the support shown on the right and coincide with $p$ in the three colored vertices.}
		\label{figure.reduction}
	\end{figure}
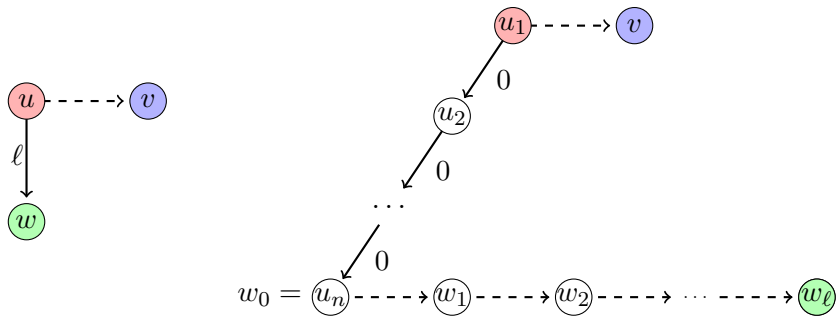

	Clearly each set $F_{p}$ is finite for each $p$. We define $F_{Z} := \bigcup_{p\in {F}_{\B,\Sigma}^n}F_{p}$. As ${F}_{\B,\Sigma}^n$ is finite, we conclude that $F_{Z}$ is finite. It is easy to see that it can be effectively constructed from $\widetilde{F}_{\B,\Sigma}^n$. We claim that $Z \subset ({B}_{\Sigma}^n)^{\Gamma_{\Omega}} = \emptyset$ if and only if $X_{\Sigma} = \emptyset$.

	Indeed, suppose $Z \neq \emptyset$ and let $z \in Z$. We can define a configuration $y \in ({B}_{\Sigma}^n)^{\Gamma_{\Omega^{n,0}}}$ by setting $y(i,j) = z(i\cdot n,j)$. It follows from the definition of $F_Z$ that no patterns from $\widetilde{F}_{\Sigma}^n$ appear in $y$ and hence $y \in {Y}_{\Sigma}^{0,n}$. In turn, this implies that $X_{\Sigma}\neq \emptyset$. Conversely, if $X_{\Sigma}\neq \emptyset$ we have that each ${Y}_{\Sigma}^{k,n}$ is non-empty. Let $y^{(k)} \in {Y}_{\Sigma}^{k,n}$ and define \[z(i,j) = y^{(i\mod{n})}\left( \left\lfloor \frac{i}{n}\right\rfloor ,j\right).\] From the definition of $F_Z$ it follows that no forbidden patterns appear in $z$ and hence $z \in Z$. It follows that if $\texttt{DP}(\Gamma_{\Omega})$ is decidable, then so is $\texttt{DP}(\Gamma_{\Xi})$. Using the result of Kari (\cref{teorema_jarkko}) we have that $\texttt{DP}(\Gamma_{\Xi})$ is undecidable, hence $\texttt{DP}(\Gamma_{\Omega})$ is also undecidable.\end{proof}

\section{The domino problem for surface groups}\label{section.surface}

A fundamental result of geometry is that up to homeomorphism, closed orientable surfaces are completely classified by their genus $g$: any such surface is either homeomorphic to a sphere or to a finite connected sum of tori. In this section we shall classify the domino problem of their fundamental groups.

\subsection{Surface groups}

We say that a group $G$ is a \emph{surface group} if it is isomorphic to the fundamental group of a closed orientable surface of genus $g \geq 1$. The case $g = 0$ is not very interesting because the sphere is simply connected and thus its fundamental group is trivial. The surface group of genus $g$ admits the following presentation

$$G_g \cong \langle a_1,b_1,\dots,a_g,b_g\mid [a_1,b_1]\dots[a_g,b_g]\rangle,$$

where $[a,b]=aba^{-1}b^{-1}$ is the commutator of $a$ and $b$.

\medskip

In the case $g = 1$ the corresponding surface group is the fundamental group of the torus, namely $\Z^2 \cong \langle a,b \mid aba^{-1}b^{-1} \rangle$, and hence by Berger's result~\cite{BergerPhD} its associated domino problem is undecidable.

The domino problem for a finitely generated group is known to be independent of the choice of generating set and a commensurability invariant~\cite[Corollary 9.53]{BertheRigo2018}. It turns out that all surface groups of genus $g\geq2$ are commensurable (see~\cite[Proposition 6.7]{ClimenhagaKatok} for a recent reference). By combining these two facts, it would be enough to prove the undecidability of the domino problem for just the surface group of genus $2$. In the sequel, we shall denote by $G$ the surface group of genus~$2$, i.e. the group with finite presentation
$$G \cong \langle a,b,c,d\mid [a,b][c,d]\rangle,$$
and denote by $S$ its generating set $\{a,b,c,d\}$.

\medskip

The Cayley graph of the surface group of genus~$2$ associated with the presentation above is not an orbit graph of some substitution with an expanding eigenvalue, but can be seen as such just by assigning labels to the edges. Moreover we shall see that these labels can be obtained \textit{locally}, which means that we can code it through an SFT.

\subsection{Finding a substitution in the surface group of genus $2$}\label{subsection.bijection}

The goal of this section is to make the parallel between the Cayley graph of the surface group $\mC_G:=\Gamma(G,S)$ and the orbit graph of a particular substitution.

\begin{figure}[htp]
	\centering
	\begin{tikzpicture}

\newcommand*{\octSize}{1.7}%
\newcommand*{\newOpacity}{0.3}%

\coordinate (n0) at (0,1*\octSize);
\coordinate (n1) at ({-1/sqrt(2)*\octSize}, {1/sqrt(2)*\octSize});
\coordinate (n2) at (-1*\octSize,0);
\coordinate (n3) at ({-1*\octSize/sqrt(2)}, {-1*\octSize/sqrt(2)});
\coordinate (n4) at (0,-1*\octSize);
\coordinate (n5) at ({1/sqrt(2)*\octSize}, {-1/sqrt(2)*\octSize});
\coordinate (n6) at (1*\octSize,0);
\coordinate (n7) at ({1/sqrt(2)*\octSize}, {1/sqrt(2)*\octSize});

\draw[->, >=stealth] (n0) -- (n1) node [pos=0.5, above] {$a$} ;
\draw[->, >=stealth] (n1) -- (n2) node [pos=0.5, left] {$b$} ;
\draw[->, >=stealth] (n2) -- (n3) node [pos=0.5, left] {$a^{-1}$} ;
\draw[->, >=stealth] (n3) -- (n4) node [pos=0.5, below left] {$b^{-1}$} ;
\draw[->, >=stealth] (n4) -- (n5) node [pos=0.5, below] {$c$} ;
\draw[->, >=stealth] (n5) -- (n6) node [pos=0.5, right] {$d$} ;
\draw[->, >=stealth] (n6) -- (n7) node [pos=0.5, right] {$c^{-1}$} ;
\draw[->, >=stealth] (n7) -- (n0) node [pos=0.5, above right] {$d^{-1}$} ;

\draw[->, >=stealth, line width=0.05mm, opacity=\newOpacity] (n0) -- (n2) ;
\draw[->, >=stealth, line width=0.5mm] (n0) -- (n3) node [pos=0.5, right] {$aba^{-1}$} ;
\draw[->, >=stealth, line width=0.05mm, opacity=\newOpacity] (n0) -- (n4) ;
\draw[->, >=stealth, line width=0.05mm, opacity=\newOpacity] (n0) -- (n5) ;
\draw[->, >=stealth, line width=0.05mm, opacity=\newOpacity] (n0) -- (n6) ;

\draw[->, >=stealth, line width=0.05mm, opacity=\newOpacity] (n1) -- (n3) ;
\draw[->, >=stealth, line width=0.05mm, opacity=\newOpacity] (n1) -- (n4) ;
\draw[->, >=stealth, line width=0.05mm, opacity=\newOpacity] (n1) -- (n5) ;
\draw[->, >=stealth, line width=0.05mm, opacity=\newOpacity] (n1) -- (n6) ;
\draw[->, >=stealth, line width=0.05mm, opacity=\newOpacity] (n1) -- (n7) ;

\draw[->, >=stealth, line width=0.05mm, opacity=\newOpacity] (n2) -- (n4) ;
\draw[->, >=stealth, line width=0.05mm, opacity=\newOpacity] (n2) -- (n5) ;
\draw[->, >=stealth, line width=0.05mm, opacity=\newOpacity] (n2) -- (n6) ;
\draw[->, >=stealth, line width=0.05mm, opacity=\newOpacity] (n2) -- (n7) ;

\end{tikzpicture}
	\caption{An elementary cycle.}
	\label{fig:diagonals}
\end{figure}
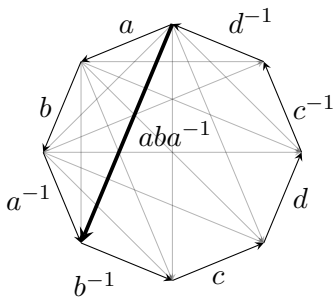

The group $G$ has only one relation $[a,b][c,d]=1_G$. Thus the only minimal cycles of the Cayley graph are cyclic permutations of $[a,b][c,d]$. We call them \emph{elementary cycles}. Moreover, any edge in the Cayley graph is part of at least one elementary cycle, since all generators and their inverses appear in the relation. Consider the Cayley graph $\Gamma(G, D)$ with
\[\mbox{$D = \left\{w ~|~ w \text{ subword of a cyclic permutation of } [a,b][c,d]\right\}$}.\]
This corresponds to adding all cords in every elementary cycle of the Cayley graph $\Gamma(G, S)$ (see \cref{fig:diagonals}). Let $d$ be the distance on $G$ given by:
\[ d(g, h) = \min \{|w| ~|~ w\in D^*, gw =_G h \} .\] In other words, $d(g, h)$ is the smallest number of elementary cycles that must be crossed to go from $g$ to $h$ in $\mC_G$. Let $B_i = \{g\in G ~|~ d(1_G, g) \leq i\}$ be the ball of radius $i$ and $C_i = \{g\in G ~|~ d(1_G, g) = i\}$ be the sphere of radius $i$, so that $B_{i+1} \setminus B_i = C_{i+1}$. Since $S \subset D$ and $\langle S \rangle = G$ we have that the $B_i$ partition $G$.

Consider $C_i$ for $i \geq 1$. In $\Gamma(G,S)$, every element of $C_i$ has exactly two neighbors in $C_i$ and either ($\texttt{a}$) one or ($\texttt{b}$) no neighbors in $C_{i-1}$. We must therefore have that there are $5$ and $6$ neighbors in $C_{i+1}$ for types ($\texttt{a}$) and ($\texttt{b}$) respectively. More specifically, if we now consider $\Gamma(G,D)$, it can be verified that the sequence of elements of $C_{i+1}$ that is obtained by adding elementary cycles to an element of type ($\texttt{a}$) has the type sequence $\texttt{a}\texttt{b}^5\texttt{a}\texttt{b}^5\texttt{a}\texttt{b}^5\texttt{a}\texttt{b}^5\texttt{a}\texttt{b}^4$ and the sequence of types for an element of type $(\texttt{b})$ is $\texttt{a}\texttt{b}^5\texttt{a}\texttt{b}^5\texttt{a}\texttt{b}^5\texttt{a}\texttt{b}^5\texttt{a}\texttt{b}^5\texttt{a}\texttt{b}^4$. This leads us to define the substitution $s\colon \{\texttt{a},\texttt{b}\} \to \{\texttt{a},\texttt{b}\}^*$ given by
\[
\left\{
\begin{aligned}
s(\texttt{a}) &= (\texttt{a}\texttt{b}^5)^4 \texttt{a}\texttt{b}^4\\
s(\texttt{b}) &= (\texttt{a}\texttt{b}^5)^5 \texttt{a}\texttt{b}^4.
\end{aligned}
\right.
\] \label{definition.substitution}

From now on, we fix $\Omega=\left(\omega^i,P_i\right)_{i\in\Z}$ an orbit of the substitution $s$ defined above, and denote by $\Gamma$ its associated orbit graph. Let us note that $s$ admits an expanding eigenvalue ($\lambda = 17+12\sqrt{2}$ and $v(\texttt{b})/v(\texttt{a}) = \frac{1+\sqrt{2}}{2}$).

The similarities between the two graphs will allow us to perform a reduction from the domino problem on $\Gamma$ (shown to be undecidable in \cref{section.undecidability_orbit_graphs}) to the domino problem on the surface group of genus 2. In order to do this reduction, all we need is a computable map which sends sets of pattern codings over $\Gamma$ into sets of pattern codings over $\mC_G$ such that the sets defining a non-empty subshift are mapped into sets defining a non-empty subshift and vice-versa. This is not trivial to do, because some of the edges are "lost" going from $\Gamma$ to $\mC_G$. In order to "recover" those edges, we shall construct an SFT $X$ over $G$ which recovers the lost information locally and use it to build the bijection needed for the reduction.
Note that technically we don't need an SFT to do so, a computable bijection would be enough. However doing it with an SFT provides a locally computable mapping, which is a nice bonus.

\subsubsection*{Definition of \texorpdfstring{$X$}{X}}

To define the SFT $X$, we introduce a notion of directions that will correspond to following edges of the orbit graph. More formally, let us first consider the general alphabet $\A_0$, consisting of the tuples
\[  \left( c, (h_1, d_1), (h_2, d_2), \hdots , (h_8, d_8) \right) \]
such that:
\begin{itemize}
	\item $c\in\{\blacksquare, \square\}$ is a color,
	\item $(h_1, \hdots, h_8)$ is a permutation of $S\cup S^{-1} = \{a,a^{-1},b,b^{-1},c,c^{-1},d,d^{-1}\}$,
	\item $d_1,\hdots,d_8\in\{\leftarrow,\rightarrow,\uparrow,\downarrow_1,\downarrow_2,\downarrow_3,\downarrow_4,\downarrow_5,\downarrow_6\}$ the directions associated to each generator.
\end{itemize}
Let $x\in \A_0^G$ be a configuration over $\A_0$. 
For every $g\in G$, if the first coordinate of $x_g$ is $c=\blacksquare$ (resp. $\square$), we call $x_g$ a black (resp. white) \emph{cell}.
We also denote by $(h,d)\sqsubset x_g$ the fact that the cell $x_g$ has direction $d$ associated to the generator $h$.

The alphabet $\A_1\subseteq \A_0$ is made of three types of elements with more precise directions imposed, depending on the color $c$:
\[  ( \blacksquare, (h_1, \leftarrow), (h_2, \rightarrow), (h_3, \uparrow), (h_4, \downarrow_1), (h_5, \downarrow_2), (h_6, \downarrow_3), (h_7, \downarrow_4), (h_8, \downarrow_5) ) \]
\[  ( \square, (h_1, \leftarrow), (h_2, \rightarrow), (h_3, \downarrow_1), (h_4, \downarrow_2), (h_5, \downarrow_3), (h_6, \downarrow_4), (h_7, \downarrow_5), (h_8, \downarrow_6) ) \]
Black cells have directions \emph{left}, \emph{right}, \emph{up} and \emph{down}, whereas whites ones have only \emph{left}, \emph{right} and \emph{down}.
Note that for both cells, \emph{up}, \emph{left} and \emph{right} are unique. We can then define their top, left and right neighbors.

\begin{definition}
  Let $x\in \A_1^G$ be a configuration over $\A_1$ and $g\in G$. We define:
  \begin{itemize}
  \item $g h_1$ the \emph{left neighbor of $g$} in $x$, denoted by $\leftgen_x(g)$,
  \item $g h_2$ is the \emph{right neighbor of $g$} in $x$, denoted by $\rightgen_x(g)$,
  \item If $x_g$ is a black cell, $g h_3$ is the \emph{top neighbor of $g$} in $x$, denoted by $\upgen_x(g)$,
  \item $g h_{3+i}$ for $i\in\{1,...,5\}$,  (resp. $g h_{2+i}$ for $i\in\{1,...,6\}$ for a white cell) is the \emph{$i$-th bottom neighbor of $g$} in $x$, denoted by $\downgen_{i,x}(g)$.
  \end{itemize}
\end{definition}

\bigskip

\begin{figure}[htp]
	\centering
	\begin{tikzpicture}

\draw (1,0) rectangle (1.5,0.5) node[pos=.5] {\large $*$};

\draw[fill=black] (-2,-1.5) rectangle (-1.5,-1) node[pos=.5, color=white] {\large $a$};
 \draw[->, >=stealth] (-1.5, -1.25) -- (-1, -1.25) node [pos=0.5, above] {$g_2$} ;
\draw[fill=white] (-1,-1.5) rectangle (-0.5,-1)  node[pos=.5] {\large $b$};
 \draw[->, >=stealth] (-0.5, -1.25) -- (0, -1.25) node [pos=0.5, above] {$g_3$} ;
\draw[fill=white] (0,-1.5) rectangle (0.5,-1);
 \draw[->, >=stealth] (0.5, -1.25) -- (1, -1.25) node [pos=0.5, above] {$g_4$} ;
 
\draw[fill=white] (1,-1.5) rectangle (1.5,-1);
 \draw[->, >=stealth] (1.5, -1.25) -- (2, -1.25) node [pos=0.5, above] {$g_5$} ;

\begin{scope}[shift={(4,0)}]
\draw[fill=white] (-2,-1.5) rectangle (-1.5,-1);
 \draw[->, >=stealth] (-1.5, -1.25) -- (-1, -1.25) node [pos=0.5, above] {$g_6$} ;
\draw[fill=white] (-1,-1.5) rectangle (-0.5,-1);
 \draw[->, >=stealth] (-0.5, -1.25) -- (0, -1.25) node [pos=0.5, above] {$g_7$} ;
\draw[fill=black] (0,-1.5) rectangle (0.5,-1);
\end{scope}

 \draw[->, >=stealth] (1.25, 0) -- (-1.75, -1) node [pos=0.5, above] {$g_1$} ;
 \draw[->, >=stealth] (4.25, -1) -- (1.25, 0) node [pos=0.5, above] {$g_8$} ;

\end{tikzpicture}
	\quad
	\begin{tikzpicture}

\draw (0,0) rectangle (0.5,0.5) node[pos=.5] {\large $*$};
\draw[<-, >=stealth] (0.5, 0.25) --  (1, 0.25) node [pos=0.5, above] {$g_8$} ;
\draw (1,0) rectangle (1.5,0.5) node[pos=.5] {\large $*$};

\draw[fill=black] (-2,-1.5) rectangle (-1.5,-1) node[pos=.5, color=white] {\large $a$};
 \draw[->, >=stealth] (-1.5, -1.25) -- (-1, -1.25) node [pos=0.5, above] {$g_2$} ;
\draw[fill=white] (-1,-1.5) rectangle (-0.5,-1)  node[pos=.5] {\large $b$};
 \draw[->, >=stealth] (-0.5, -1.25) -- (0, -1.25) node [pos=0.5, above] {$g_3$} ;
\draw[fill=white] (0,-1.5) rectangle (0.5,-1);
 \draw[->, >=stealth] (0.5, -1.25) -- (1, -1.25) node [pos=0.5, above] {$g_4$} ;

\begin{scope}[shift={(3,0)}]
\draw[fill=white] (-2,-1.5) rectangle (-1.5,-1);
 \draw[->, >=stealth] (-1.5, -1.25) -- (-1, -1.25) node [pos=0.5, above] {$g_5$} ;
\draw[fill=white] (-1,-1.5) rectangle (-0.5,-1);
 \draw[->, >=stealth] (-0.5, -1.25) -- (0, -1.25) node [pos=0.5, above] {$g_6$} ;
\draw[fill=black] (0,-1.5) rectangle (0.5,-1);
\end{scope}

 \draw[->, >=stealth] (0.25, 0) -- (-1.75, -1) node [pos=0.5, above] {$g_1$} ;
 \draw[->, >=stealth] (3.25, -1) -- (1.25, 0) node [pos=0.5, above] {$g_7$} ;


\end{tikzpicture}
	\caption{The two possible types of colorings of cycles. There are no color constraints \mbox{on $\boxed{*}$}, and the cycle $g_1 \hdots g_8$ is any cyclic permutation of $[a,b][c,d]$. }
	\label{fig:cycle_color}
\end{figure}
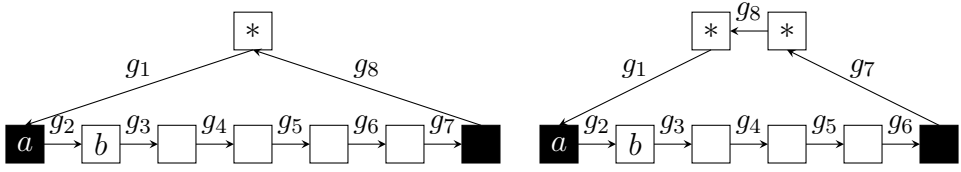

We forbid all patterns whose support is an elementary circle and that do not have the colors shown on \cref{fig:cycle_color}.
We also impose the orientations to be as drawn.
For example, the right of \scalebox{0.8}{
\begin{tikzpicture}

\draw[fill=black] (-2,-1.5) rectangle (-1.5,-1) node[pos=.5, color=white] {\large $a$};

\end{tikzpicture}
} is $g_2$, its top is $g_1^{-1}$, and the other directions of $a$ are not constrained by this cycle.
Similarly, the left of \scalebox{0.8}{
\begin{tikzpicture}

\draw[fill=white] (-1,-1.5) rectangle (-0.5,-1)  node[pos=.5] {\large $b$};

\end{tikzpicture}
} is $g_2^{-1}$, its right $g_3$ and other directions unconstrained.
To do so, we call $\F_1$ the set of all elementary cycles that are not of the form of \cref{fig:cycle_color}.

We add the constraint that directions must be consistent between adjacent cells, by forbidding the finite set $\F_2$ , which is the set of patterns on the support $\{ 1_G, h \}$ for $h\in S$, such that $x_{1_G}$ and $x_h$ are linked by mismatching directions. That is,
\[
  \F_2 = \left\{
  \text{pattern $p$ of support}\{ 1_G, h \} ~\middle|~
  \begin{split}
     \leftgen_p(1_G) = h &\text{ and } \rightgen_p(h) \neq 1_G  \text{ or}\\
  	 \rightgen_p(1_G) = h &\text{ and } \leftgen_p(h) \neq 1_G  \text{ or}\\
  	 \upgen_p(1_G) = h &\text{ and } \forall i, \downgen_{i,p}(h) \neq 1_G  \text{ or}\\
  	 \exists i, \downgen_{i,p}(1_G) = h &\text{ and } \upgen_p(h) \neq 1_G
  \end{split}
  \right\}.
\]
We define $X$ as the set of all configurations over $\A_1$ where no forbidden patterns from $\F_1\cup\F_2$ appear. As $\F_1\cup\F_2$ is finite, $X$ is an SFT.

\subsubsection*{Non-emptiness of \texorpdfstring{$X$}{X}}

We construct a configuration $x$ in the SFT $X\subset \A_1^G$ as the limit of a sequence of configurations $(y_n)_{n\in\N}$ of another SFT $X_2\subset \left(\A_1 \cup \{\texttt{orange}\}\right)^G$, where
\[ \texttt{orange}:=\left( \orangesquare, (a, \downarrow_1), (a^{-1}, \downarrow_2), (b, \downarrow_3), (b^{-1}, \downarrow_4), (c, \downarrow_5), (c^{-1}, \downarrow_6), (d, \downarrow_7), (d^{-1}, \downarrow_8) \right) ,\]
and we extend the definition of neighbors consistently.
$X_2$ is defined by $\F = \F_1\cup\F_2$ the same finite set of forbidden patterns as $X$. Intuitively, because the letter $\texttt{orange}$ has only bottom neighbors, the presence of an orange cell creates \textit{rings} (see \cref{fig:B1} and \cref{lemma:ring_patterns}).

\begin{definition} \label{def:ring-line}
  We say $L \subset G$ is a set of \emph{left-right neighbors} of $x$ if we can access all its elements by taking only their left and right neighbors, i.e. for every $g\in L$, we have \mbox{$L=\{\hdots, \leftgen_x^3(g),\leftgen_x^2(g), \leftgen_x(g), g, \rightgen_x(g), \rightgen_x^2(g), \rightgen_x^3(g), \hdots \}$}.

  If $L$ is finite, it is called a \emph{ring}, if it is infinite it is called a \emph{line}.
\end{definition}

\begin{lemma} \label{lemma:ring_patterns}
  For all $i$, there exists a pattern $p_i\in\left(\A_1 \cup \{\texttt{orange}\}\right)^{B_i}$ containing no forbidden patterns of $\F$ and such that $(p_i)_g$ is an orange cell if and only if $g=1_G$.
\end{lemma}
\begin{proof}
By induction on $i$, we prove a stronger statement:
\begin{center}
$\mathcal{H}_i:$ "There exists a coloring of $B_i$, in which the orange tile appears, but only at the origin. Moreover, in this coloring, $C_j$ is a ring for all $j\leq i$."
\end{center}

For $i=1$, apply the first cycle of \cref{fig:cycle_color} eight times, and from the orange origin, get the sphere of radius 1, which is a cycle as stated (see \cref{fig:B1}).

\begin{figure}[htp]
	\centering
	\begin{tikzpicture}

\newcommand*{\squareFactor}{1.9}%

\draw (0,0) circle (2.5);

\foreach \i in {0, 7.5, ..., 359} {
	\begin{scope} [rotate around={\i:(0,0)}]
		\draw [fill=white] (-2.5-0.25/\squareFactor,-0.25/\squareFactor) rectangle (-2.5+0.25/\squareFactor,0.25/\squareFactor);
	\end{scope}
}

\foreach \i in {0, 45, ..., 359} {
	\begin{scope} [rotate around={\i:(0,0)}]
	    \draw (0,0) -- (-2.5,0);
		\draw [fill=black] (-2.5-0.25/\squareFactor,-0.25/\squareFactor) rectangle (-2.5+0.25/\squareFactor,0.25/\squareFactor);
	\end{scope}
}

\draw [fill=orange] (-0.2,-0.2) rectangle (0.2,0.2);

\end{tikzpicture}
	\caption{Coloring of $B_1$.}
	\label{fig:B1}
\end{figure}
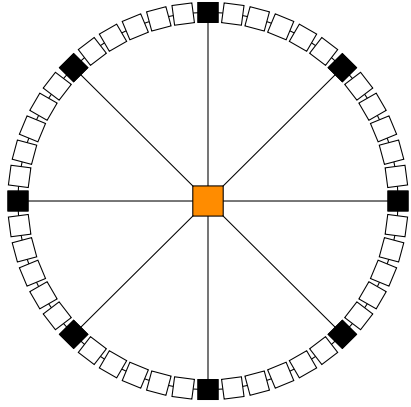

Now suppose we have a coloring of $B_i$ as in the statement.
We will use the cycles of \cref{fig:cycle_color} on the ring $C_i$ to build $C_{i+1}$.
We are sure that all the cells on $C_i$ are only black and whites due to the induction hypothesis on the orange cell.
Each of the black cells on $C_i$ must have 5 \emph{bottom} cells, and each white one needs 6.
We proceed iteratively, starting from any cell $c$ and any \emph{down} generator $g_1$ of this cell.
As the two possible cycles start the same, we put the colors $\blacksquare$, $\square$, $\square$, $\square$, $\square$ following the generators $g_1, g_2, g_3, g_4, g_5$, with the consistent orientation.
For the two next colors, it depends where $g_8$ leads.
If it lead to $C_{i+1}$, we are in the first case of \cref{fig:Bi+1}, and we use the colors of the first cycle. We start the process again but with cell $c$ and generator $g_8^{-1}$.
If $g_8$ leads to $C_i$, we are in the second case of \cref{fig:Bi+1}, and we use the corresponding colors. We then start again but with the cell $c'$ and generator $g_7^{-1}$.
We continue this process until all cells of $C_i$ have their bottom neighbors colored.

\begin{figure}[htp]
\centering
\begin{tikzpicture}

\foreach \i in {50, 63, ..., 130} {
	\begin{scope} [rotate around={4.5+\i:(0,0)}]
		\draw [fill=white] (-2.5,0) -- (-2,0);
	\end{scope}
}

\draw node at (2,-1) {$C_i$};

\draw node at (0,0) {$B_{i-1}$};

\draw [fill=black, opacity=0.3] (0,0) circle (2);
\draw (0,0) circle (2);
\draw [fill=black, opacity=0.15] (0,0) circle (2.5);
\draw (0,0) circle (2.5);

\begin{scope}[shift={(-1.4,-2.5)}]
	\draw[fill=black] (-2,-1.5) rectangle (-1.5,-1) node[pos=.5, color=white] {};
	 \draw[->, >=stealth] (-1.5, -1.25) -- (-1, -1.25) node [pos=0.5, above] {$g_2$} ;
	\draw[fill=white] (-1,-1.5) rectangle (-0.5,-1)  node[pos=.5] {};
	 \draw[->, >=stealth] (-0.5, -1.25) -- (0, -1.25) node [pos=0.5, above] {$g_3$} ;
	\draw[fill=white] (0,-1.5) rectangle (0.5,-1);
	 \draw[->, >=stealth] (0.5, -1.25) -- (1, -1.25) node [pos=0.5, above] {$g_4$} ;
	 
	\draw[fill=white] (1,-1.5) rectangle (1.5,-1);
	 \draw[->, >=stealth] (1.5, -1.25) -- (2, -1.25) node [pos=0.5, above] {$g_5$} ;

	\begin{scope}[shift={(4,0)}]
	\draw[fill=white] (-2,-1.5) rectangle (-1.5,-1);
	 \draw[->, >=stealth] (-1.5, -1.25) -- (-1, -1.25) node [pos=0.5, above] {$g_6$} ;
	\draw[fill=white] (-1,-1.5) rectangle (-0.5,-1);
	 \draw[->, >=stealth] (-0.5, -1.25) -- (0, -1.25) node [pos=0.5, above] {$g_7$} ;
	\draw[fill=black] (0,-1.5) rectangle (0.5,-1);
	\end{scope}

	 \draw[->, >=stealth] (1.25, 0) -- (-1.75, -1) node [pos=0.5, above] {$g_1$} ;
	 \draw[->, >=stealth] (4.25, -1) -- (1.25, 0) node [pos=0.5, above] {$g_8$} ;
\end{scope}

\draw node at (-0.1,-2.25) {c};

\end{tikzpicture}
~~~~~~~~~~
\begin{tikzpicture}

\foreach \i in {50, 63, ..., 130} {
	\begin{scope} [rotate around={4.5+\i:(0,0)}]
		\draw [fill=white] (-2.5,0) -- (-2,0);
	\end{scope}
}

\draw node at (2,-1) {$C_i$};

\draw node at (0,0) {$B_{i-1}$};

\draw [fill=black, opacity=0.3] (0,0) circle (2);
\draw (0,0) circle (2);
\draw [fill=black, opacity=0.15] (0,0) circle (2.5);
\draw (0,0) circle (2.5);

\begin{scope}[shift={(-0.6,-2.5)}]
	\draw[fill=black] (-2,-1.5) rectangle (-1.5,-1) node[pos=.5, color=white] {\large $a$};
	 \draw[->, >=stealth] (-1.5, -1.25) -- (-1, -1.25) node [pos=0.5, above] {$g_2$} ;
	\draw[fill=white] (-1,-1.5) rectangle (-0.5,-1)  node[pos=.5] {\large $b$};
	 \draw[->, >=stealth] (-0.5, -1.25) -- (0, -1.25) node [pos=0.5, above] {$g_3$} ;
	\draw[fill=white] (0,-1.5) rectangle (0.5,-1);
	 \draw[->, >=stealth] (0.5, -1.25) -- (1, -1.25) node [pos=0.5, above] {$g_4$} ;

	\begin{scope}[shift={(3,0)}]
	\draw[fill=white] (-2,-1.5) rectangle (-1.5,-1);
	 \draw[->, >=stealth] (-1.5, -1.25) -- (-1, -1.25) node [pos=0.5, above] {$g_5$} ;
	\draw[fill=white] (-1,-1.5) rectangle (-0.5,-1);
	 \draw[->, >=stealth] (-0.5, -1.25) -- (0, -1.25) node [pos=0.5, above] {$g_6$} ;
	\draw[fill=black] (0,-1.5) rectangle (0.5,-1);
	\end{scope}

	 \draw[->, >=stealth] (0.45, 0) -- (-1.75, -1) node [pos=0.5, above] {$g_1$} ;
	 \draw[->, >=stealth] (3.25, -1) -- (1.05, 0.038) node [pos=0.5, above] {$g_7$} ;
	 \draw[->, >=stealth] (1.05, 0.038) to [controls=+(55:-0.25) and +(140:-0.25)] (0.45, 0);
	 \draw node at (0.8, -0.3) {$g_8$};
\end{scope}

\draw node at (-0.1,-2.25) {c};
\draw node at (0.42,-2.18) {c'};

\end{tikzpicture}
\caption{From $B_i$ to $B_{i+1}$.}
\label{fig:Bi+1}
\end{figure}
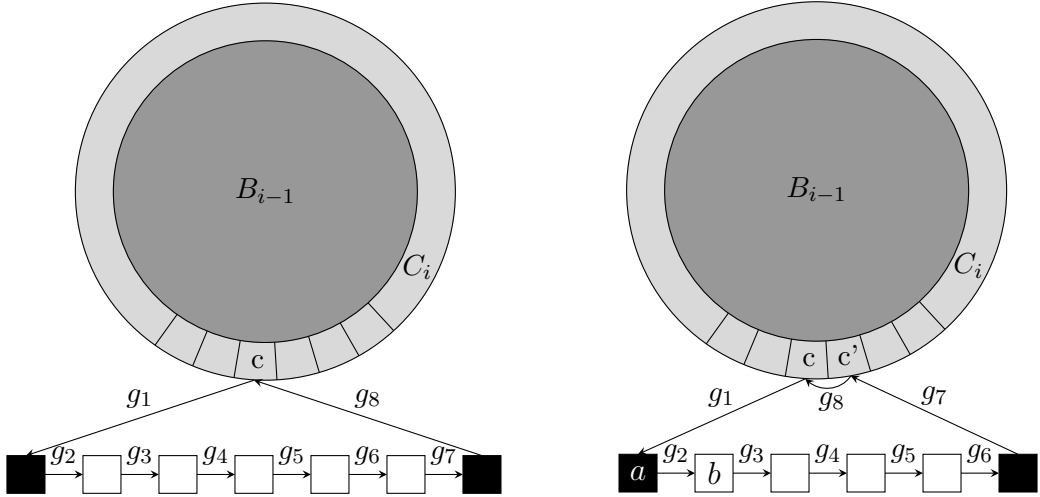

With this process, we colored a new ring, which is exactly $C_{i+1}$.
Indeed, the newly colored cells are in $C_{i+1}$, because one cycle separates them from $C_i$ and there are no other cells in $C_{i+1}$ because adding one cycle to these will increase the distance to $i+2$.

Because $B_{i+1} = B_i \cup C_{i+1}$, we now have colored $B_{i+1}$. We have not placed any new orange tile, so the only one is the one from $B_i$ i.e., by induction hypothesis, the origin. Therefore, the statement is proved for $i+1$.

\end{proof}

By compactness of $\left(\A_1 \cup \{\texttt{orange}\}\right)^G$ there exists a configuration $\widetilde{x} \in X_2$ which coincides with the pattern $p_i$ on $B_i$ for all $i \in \N$. In particular, \texttt{orange} appears only at the origin. Note that $X$ consists of all configurations of $X_2$ where \texttt{orange} does not appear. As \texttt{orange} only appears at the origin of $\widetilde{x}$, we can find arbitrarily large regions of $\widetilde{x}$ where \texttt{orange}~does not appear at all. By shift-invariance and compactness of $X_2$, there exists a configuration $x$ where it does not appear and therefore $x \in X\neq\emptyset$.

\subsubsection*{Configurations of \texorpdfstring{$X$}{X}}

Now that we know that $X$ is non-empty, we take a look at the properties of its configurations which are useful for our reduction.

We first show that without the orange tile, configurations cannot have rings, they can only have infinite lines. More precisely, we prove the contrapositive:
\begin{lemma} \label{lemma:ring}
  If there is a ring in a configuration of $X_2$, then \texttt{orange} must appear.
\end{lemma}

In particular, it means that configurations with a ring are not in $X$.

\begin{proof}
  Let $C \subset G$ be a ring of $x\in X_2$. As patterns from $\F_1$ do not appear in $X_2$, unless $C$ is a singleton and $x|_{C} = \texttt{orange}$, it must contain at least eight elements and at least two of them must be black cells and hence have top neighbors. The key point is that $C_1:=~\upgen_x(C)=\{\upgen_x(g)~|~g\in C\}$ is also a ring, but with strictly less elements.
  If there is an orange cell in $C_1$, it is done.
  If not, because all cycles are colored like \cref{fig:cycle_color}, we know that the top neighbors of $C$ are organized as a ring (we can "stick" cycles all around $C$). And this ring is strictly smaller than the previous one, because for each 7 our 6 cells of $C$ we have 1 or 2 cells in $C_1$.

\begin{figure}[htp]
	\centering
	\scalebox{0.8}{\begin{tikzpicture}


\begin{scope} [rotate around={50+-10/1.5:(0,0)}]
	\draw [->, >=stealth] (-4.75,0) -- (-2.65,-1.9);
\end{scope}

\begin{scope} [rotate around={50+50/1.5:(0,0)}]
	\draw [->, >=stealth] (-4.75,0) -- (-3.25,0);
\end{scope}

\begin{scope} [rotate around={50+110/1.5:(0,0)}]
	\draw [->, >=stealth] (-4.75,0) -- (-2.65,1.9);
\end{scope}

\begin{scope} [rotate around={50+160/1.5:(0,0)}]
	\draw [->, >=stealth] (-4.75,0) -- (-3.05,-1.1);
\end{scope}

\begin{scope} [rotate around={50+220/1.5:(0,0)}]
	\draw [->, >=stealth] (-4.75,0) -- (-3.05,+1.1);
\end{scope}

\draw (0,0) circle (5);

\foreach \i in {-20, -10, 0, 10, ..., 230} {
	\begin{scope} [rotate around={50+\i/1.5:(0,0)}]
		\draw [fill=white] (-4.75,-0.25) rectangle (-5.25,0.25);
	\end{scope}
}

\foreach \i in {-10, 50, 110, 160, 220} {
	\begin{scope} [rotate around={50+\i/1.5:(0,0)}]
		\draw [fill=black] (-4.75,-0.25) rectangle (-5.25,0.25);
	\end{scope}
}

\draw (0,0) circle (3);

\foreach \i in {50, 190} {
	\begin{scope} [rotate around={50+\i/1.5:(0,0)}]
		\draw [fill=white] (-2.75,-0.25) rectangle (-3.25,0.25);
	\end{scope}
}


\draw node at (4.4,3.3) {\Large$C$};
\draw node at (3,1.8) {\Large$C_1$};

\end{tikzpicture}}
	\caption{Top neighbors of rings are smaller rings.}
	\label{fig:rings_converge}
\end{figure}
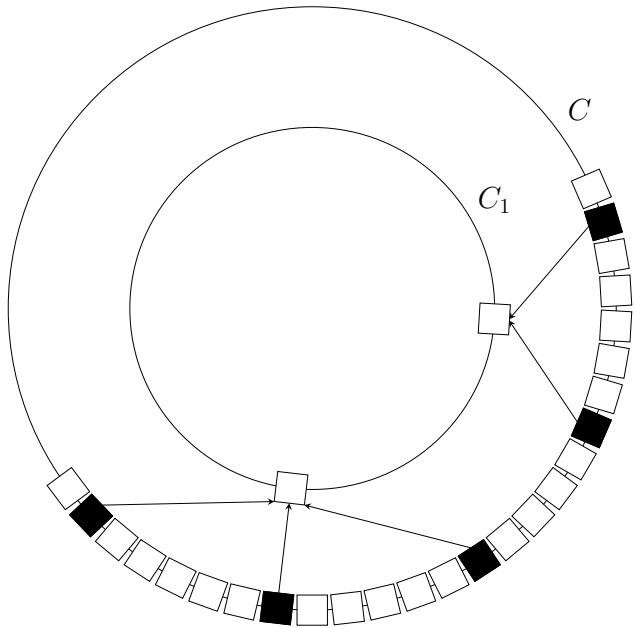

  Iterating the previous process, we reduce the size of $C$ which is finite. The process forcefully ends with the ring of size one. Note that $\{\texttt{orange}\}$ is a ring with one element. It is also the only one possible, since for any other cell $x_g$, $\rightarrow_x(g)\neq g$.
\end{proof}

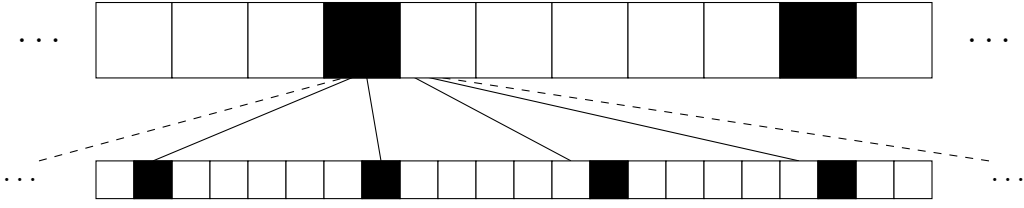
\begin{figure}[htp]
	\centering
	\begin{tikzpicture}

\newcommand*{\topdown}{0.4}%
\newcommand*{\bottomshift}{0}%

\foreach \i in {0.5, 3.5}{
	\draw[] (0.25+\i-\bottomshift, 8.5) -- (3.55, 0.08+10-\topdown) ;
}
\draw[dashed] (-0.75-\bottomshift, 8.5) -- (3.55, 0.08+10-\topdown) ;

\foreach \i in {6, 9}{
	\draw[] (0.25+\i-\bottomshift, 8.5) -- (4.05, 0.08+10-\topdown) ;
}
\draw[dashed] (11.75-\bottomshift, 8.5) -- (4.05, 0.08+10-\topdown) ;


\draw (-0.75,10.5-\topdown) node {\Large $\hdots$};
\draw (11.75,10.5-\topdown) node {\Large $\hdots$};

\foreach \i in {0,1,...,10}{
	\draw[fill=white] (0+\i,10-\topdown) rectangle (1+\i,11-\topdown) ;
}

\foreach \i in {3, 9}{
	\draw[fill=black] (0+\i,10-\topdown) rectangle (1+\i,11-\topdown) ;
}


\draw (-1,8.25) node {$\hdots$};
\draw (12,8.25) node {$\hdots$};

\foreach \i in {0,0.5,...,10.5}{
	\draw[fill=white] (0+\i-\bottomshift,8) rectangle (0.5+\i-\bottomshift,8.5) ;
}

\foreach \i in {0.5, 3.5, 6.5, 9.5}{
	\draw[fill=black] (0+\i-\bottomshift,8) rectangle (0.5+\i-\bottomshift,8.5) ;
}

\end{tikzpicture}
	\caption{Infinite lines of $X$.}
	\label{fig:lines}
\end{figure}

Because $X$ does not have the orange cell in its alphabet $\A_2$, there cannot be any rings in its configurations by \cref{lemma:ring}.
It means that starting from any element, one can take its right neighbor infinitely many times and never loop on the initial element.
This forms infinite lines (in the sense of \cref{def:ring-line}), which are all above and below the others (\mbox{\cref{fig:lines}}), thanks to the way cycles are colored.

These lines show the orbit graph structure: each black cell has 5 (black) bottom children with 24 whites on the line, and each white cell has 6 bottom children with 29 whites on the line. Exactly the same way as in the orbit graph $\Gamma$.

Moreover, we can show that these lines induce a height function on $G$: when going down, one never comes back to an upper line. This is a corollary of the following lemma.

\begin{lemma} \label{lemma:cycles}
  Let $x\in X, g\in G$ and $a_1, \hdots, a_k \in \{\leftgen_x, \rightgen_x, \upgen_x, \downgen_{1,x}, \downgen_{2,x}, \hdots \}, k>0$ such that $a_k \circ \hdots \circ a_1 (g) = g$. Then
  \[
  \left|\left\{i\in\{1,\dots,k\}\mid a_i=\upgen_x\right\}\right|
  =
  \left|\left\{j\in\{1,\dots,k\}\mid \exists k, a_j=\downgen_{k,x}\right\}\right|.
  \]
\end{lemma}

\begin{figure}[htp]
  	\centering
	\begin{subfigure}{0.33\textwidth}
	\centering
	\begin{tikzpicture}[scale=0.37]
 
\draw[thin,dashed] (0,4) -> (2,6);   
\draw[thin,postaction={decorate,decoration={markings,
    mark=at position .9 with {\arrow[scale=1.2]{latex}}}}] (0,4) -> (-1,2);  
\draw[thin,postaction={decorate,decoration={markings,
    mark=at position .9 with {\arrow[scale=1.2]{latex}}}}] (-1,2) -> (0,0);   
\draw[thin,postaction={decorate,decoration={markings,
    mark=at position .9 with {\arrow[scale=1.2]{latex}}}}] (0,0) -> (2,-2);        
\draw[thin,postaction={decorate,decoration={markings,
    mark=at position .9 with {\arrow[scale=1.2]{latex}}}}] (2,-2) -- (4,-3);
\draw[thin,postaction={decorate,decoration={markings,
    mark=at position .9 with {\arrow[scale=1.2]{latex}}}}] (4,-3) -- (6,-3);
\draw[thin,postaction={decorate,decoration={markings,
    mark=at position .9 with {\arrow[scale=1.2]{latex}}}}] (6,-3) -- (8,-2);
\draw[thin,postaction={decorate,decoration={markings,
    mark=at position .9 with {\arrow[scale=1.2]{latex}}}}] (8,-2) -- (10,0);
\draw[thin,postaction={decorate,decoration={markings,
    mark=at position .9 with {\arrow[scale=1.2]{latex}}}}] (10,0) -- (11,2);   
\draw[thin,postaction={decorate,decoration={markings,
    mark=at position .9 with {\arrow[scale=1.2]{latex}}}}] (11,2) -- (10,4);
\draw[thin,dashed] (10,4) -> (8,6);   

\begin{scope}[shift={(0.5,0.1)},scale=0.9,color=bleu]
\draw[thin,postaction={decorate,decoration={markings,
    mark=at position .9 with {\arrow[scale=1.2]{latex}}}}] (0,0) -> (2,-2);  
\draw[thin,postaction={decorate,decoration={markings,
    mark=at position .9 with {\arrow[scale=1.2]{latex}}}}] (2,-2) -- (4,-3); 
\draw[thin,postaction={decorate,decoration={markings,
    mark=at position .9 with {\arrow[scale=1.2]{latex}}}}] (4,-3) -- (6,-3); 
\draw[thin,postaction={decorate,decoration={markings,
    mark=at position .9 with {\arrow[scale=1.2]{latex}}}}] (6,-3) -- (8,-2);
\draw[thin,postaction={decorate,decoration={markings,
    mark=at position .9 with {\arrow[scale=1.2]{latex}}}}] (8,-2) -- (10,0);   
\draw[thin,postaction={decorate,decoration={markings,
    mark=at position .9 with {\arrow[scale=1.2]{latex}}}}] (10,0) -- (7,2);
\draw[thin,postaction={decorate,decoration={markings,
    mark=at position .9 with {\arrow[scale=1.2]{latex}}}}] (7,2) -- (3.5,2);     
\draw[thin,postaction={decorate,decoration={markings,
    mark=at position .9 with {\arrow[scale=1.2]{latex}}}}] (3.5,2) -- (0,0);  
    
\node[fill=white] at (2,-2) {};     
\node[fill=white] at (4,-3) {};      
\node[fill=white] at (6,-3) {};    
\node[fill=white] at (8,-2) {};  
\node[fill=white] at (10,0) {};
\node[fill=white] at (7,2) {}; 
\node[fill=white] at (3.5,2) {}; 
\node[fill=white] at (0,0) {};

\draw[red, very thick, ->] (5.25,1.5) -- (5.25,3); 
\draw[red, very thick, ->] (1.75,0.5) -- (1.75,2); 
\end{scope}

\begin{scope}[shift={(0.5,0.5)},scale=0.9,color=vert]
\draw[thin,dashed] (0,4) -> (2,6);  
\draw[thin,postaction={decorate,decoration={markings,
    mark=at position .9 with {\arrow[scale=1.2]{latex}}}}] (0,4) -> (-1,2);  
\draw[thin,postaction={decorate,decoration={markings,
    mark=at position .9 with {\arrow[scale=1.2]{latex}}}}] (-1,2) -- (0,0); 
\draw[thin,postaction={decorate,decoration={markings,
    mark=at position .9 with {\arrow[scale=1.2]{latex}}}}] (0,0) -- (3.5,2); 
\draw[thin,postaction={decorate,decoration={markings,
    mark=at position .9 with {\arrow[scale=1.2]{latex}}}}] (3.5,2) -- (7,2);    
\node[fill=white] at (0,4) {};
\node[fill=white] at (-1,2) {};
\node[fill=white] at (0,0) {};
\node[fill=white] at (3.5,2) {};
\node[fill=white] at (7,2) {};

\end{scope}

\begin{scope}[shift={(0.7,0.4)},scale=0.9,color=vert]  
\draw[thin,postaction={decorate,decoration={markings,
    mark=at position .9 with {\arrow[scale=1.2]{latex}}}}] (7,2) -- (10,0);     
\draw[thin,postaction={decorate,decoration={markings,
    mark=at position .9 with {\arrow[scale=1.2]{latex}}}}] (10,0) -- (11,2);  
\draw[thin,postaction={decorate,decoration={markings,
    mark=at position .9 with {\arrow[scale=1.2]{latex}}}}] (11,2) -- (10,4);  
\draw[thin,dashed] (10,4) -> (8,6);       
\node[fill=white] at (4.5,1) {};
\node[fill=white] at (7,2) {};
\node[fill=white] at (10,0) {};
\node[fill=white] at (11,2) {};   
\node[fill=white] at (10,4) {};    

\end{scope}

\node[fill=white] at (0,4) {};
\node[fill=white] at (-1,2) {};
\node[fill=white] at (0,0) {};
\node[fill=white] at (2,-2) {};
\node[fill=white] at (4,-3) {};
\node[fill=white] at (6,-3) {};
\node[fill=white] at (8,-2) {};
\node[fill=white] at (10,0) {};
\node[fill=white] at (11,2) {};
\node[fill=white] at (10,4) {};

\node[color=bleu] at (5,-2) {$\boldsymbol{\gamma_0}$};
\node[color=vert] at (2.5,4.5) {$\boldsymbol{\gamma'}$};
\node[] at (11,0.5) {$\boldsymbol{\gamma}$};

\node[color=rouge] at (2.75,2.5) {$\boldsymbol{a_i}$};

\node[color=rouge] at (5.25,0.5) {$\boldsymbol{a_j}$};
\end{tikzpicture}
	\caption{$a_i,a_j\notin\gamma_0\cap\gamma$}
	\label{fig:lemma5.4_case1}
	\end{subfigure}%
	\begin{subfigure}{0.33\textwidth}
	\centering
  	\begin{tikzpicture}[scale=0.37]
 
\draw[thin,dashed] (0,4) -> (2,6);   
\draw[thin,postaction={decorate,decoration={markings,
    mark=at position .9 with {\arrow[scale=1.2]{latex}}}}] (0,4) -> (-1,2);  
\draw[thin,postaction={decorate,decoration={markings,
    mark=at position .9 with {\arrow[scale=1.2]{latex}}}}] (-1,2) -> (0,0);   
\draw[thin,postaction={decorate,decoration={markings,
    mark=at position .9 with {\arrow[scale=1.2]{latex}}}}] (0,0) -> (2,-2);        
\draw[thin,postaction={decorate,decoration={markings,
    mark=at position .9 with {\arrow[scale=1.2]{latex}}}}] (2,-2) -- (4,-3);
\draw[thin,postaction={decorate,decoration={markings,
    mark=at position .9 with {\arrow[scale=1.2]{latex}}}}] (4,-3) -- (6,-3);
\draw[thin,postaction={decorate,decoration={markings,
    mark=at position .9 with {\arrow[scale=1.2]{latex}}}}] (6,-3) -- (8,-2);
\draw[thin,postaction={decorate,decoration={markings,
    mark=at position .9 with {\arrow[scale=1.2]{latex}}}}] (8,-2) -- (10,0);
\draw[thin,postaction={decorate,decoration={markings,
    mark=at position .9 with {\arrow[scale=1.2]{latex}}}}] (10,0) -- (11,2);   
\draw[thin,postaction={decorate,decoration={markings,
    mark=at position .9 with {\arrow[scale=1.2]{latex}}}}] (11,2) -- (10,4);
\draw[thin,dashed] (10,4) -> (8,6);   

\begin{scope}[shift={(0.5,0.1)},scale=0.9,color=bleu]
\draw[thin,postaction={decorate,decoration={markings,
    mark=at position .9 with {\arrow[scale=1.2]{latex}}}}] (0,0) -> (2,-2);  
\draw[thin,postaction={decorate,decoration={markings,
    mark=at position .9 with {\arrow[scale=1.2]{latex}}}}] (2,-2) -- (4,-3); 
\draw[thin,postaction={decorate,decoration={markings,
    mark=at position .9 with {\arrow[scale=1.2]{latex}}}}] (4,-3) -- (6,-3); 
\draw[thin,postaction={decorate,decoration={markings,
    mark=at position .9 with {\arrow[scale=1.2]{latex}}}}] (6,-3) -- (8,-2);
\draw[thin,postaction={decorate,decoration={markings,
    mark=at position .9 with {\arrow[scale=1.2]{latex}}}}] (8,-2) -- (10,0);   
\draw[thin,postaction={decorate,decoration={markings,
    mark=at position .9 with {\arrow[scale=1.2]{latex}}}}] (10,0) -- (7,2);
\draw[thin,postaction={decorate,decoration={markings,
    mark=at position .9 with {\arrow[scale=1.2]{latex}}}}] (7,2) -- (3.5,2);     
\draw[thin,postaction={decorate,decoration={markings,
    mark=at position .9 with {\arrow[scale=1.2]{latex}}}}] (3.5,2) -- (0,0);  
    
\node[fill=white] at (2,-2) {};     
\node[fill=white] at (4,-3) {};      
\node[fill=white] at (6,-3) {};    
\node[fill=white] at (8,-2) {};  
\node[fill=white] at (10,0) {};
\node[fill=white] at (7,2) {}; 
\node[fill=white] at (3.5,2) {}; 
\node[fill=white] at (0,0) {};

\draw[red, very thick, ->] (7,-2) -- (7,-3.5); 
\draw[red, very thick, ->] (9,-2) -- (9,-0.5); 
\end{scope}

\begin{scope}[shift={(0.5,0.5)},scale=0.9,color=vert]
\draw[thin,dashed] (0,4) -> (2,6);  
\draw[thin,postaction={decorate,decoration={markings,
    mark=at position .9 with {\arrow[scale=1.2]{latex}}}}] (0,4) -> (-1,2);  
\draw[thin,postaction={decorate,decoration={markings,
    mark=at position .9 with {\arrow[scale=1.2]{latex}}}}] (-1,2) -- (0,0); 
\draw[thin,postaction={decorate,decoration={markings,
    mark=at position .9 with {\arrow[scale=1.2]{latex}}}}] (0,0) -- (3.5,2); 
\draw[thin,postaction={decorate,decoration={markings,
    mark=at position .9 with {\arrow[scale=1.2]{latex}}}}] (3.5,2) -- (7,2);    
\node[fill=white] at (0,4) {};
\node[fill=white] at (-1,2) {};
\node[fill=white] at (0,0) {};
\node[fill=white] at (3.5,2) {};
\node[fill=white] at (7,2) {};
\end{scope}

\begin{scope}[shift={(0.7,0.4)},scale=0.9,color=vert]  
\draw[thin,postaction={decorate,decoration={markings,
    mark=at position .9 with {\arrow[scale=1.2]{latex}}}}] (7,2) -- (10,0);     
\draw[thin,postaction={decorate,decoration={markings,
    mark=at position .9 with {\arrow[scale=1.2]{latex}}}}] (10,0) -- (11,2);  
\draw[thin,postaction={decorate,decoration={markings,
    mark=at position .9 with {\arrow[scale=1.2]{latex}}}}] (11,2) -- (10,4);  
\draw[thin,dashed] (10,4) -> (8,6);       
\node[fill=white] at (4.5,1) {};
\node[fill=white] at (7,2) {};
\node[fill=white] at (10,0) {};
\node[fill=white] at (11,2) {};   
\node[fill=white] at (10,4) {};     

\end{scope}

\node[fill=white] at (0,4) {};
\node[fill=white] at (-1,2) {};
\node[fill=white] at (0,0) {};
\node[fill=white] at (2,-2) {};
\node[fill=white] at (4,-3) {};
\node[fill=white] at (6,-3) {};
\node[fill=white] at (8,-2) {};
\node[fill=white] at (10,0) {};
\node[fill=white] at (11,2) {};
\node[fill=white] at (10,4) {};

\node[color=bleu] at (5,-2) {$\boldsymbol{\gamma_0}$};
\node[color=vert] at (2.5,4.5) {$\boldsymbol{\gamma'}$};
\node[] at (11,0.5) {$\boldsymbol{\gamma}$};

\node[color=rouge] at (8.5,0) {$\boldsymbol{a_i}$};

\node[color=rouge] at (8,-3.5) {$\boldsymbol{a_j}$};
\end{tikzpicture}
	\caption{$a_i,a_j\in\gamma_0\cap\gamma$}
	\label{fig:lemma5.4_case2}
	\end{subfigure}%
	\begin{subfigure}{0.33\textwidth}
	\centering
  	\begin{tikzpicture}[scale=0.37]
 
\draw[thin,dashed] (0,4) -> (2,6);   
\draw[thin,postaction={decorate,decoration={markings,
    mark=at position .9 with {\arrow[scale=1.2]{latex}}}}] (0,4) -> (-1,2);  
\draw[thin,postaction={decorate,decoration={markings,
    mark=at position .9 with {\arrow[scale=1.2]{latex}}}}] (-1,2) -> (0,0);   
\draw[thin,postaction={decorate,decoration={markings,
    mark=at position .9 with {\arrow[scale=1.2]{latex}}}}] (0,0) -> (2,-2);        
\draw[thin,postaction={decorate,decoration={markings,
    mark=at position .9 with {\arrow[scale=1.2]{latex}}}}] (2,-2) -- (4,-3);
\draw[thin,postaction={decorate,decoration={markings,
    mark=at position .9 with {\arrow[scale=1.2]{latex}}}}] (4,-3) -- (6,-3);
\draw[thin,postaction={decorate,decoration={markings,
    mark=at position .9 with {\arrow[scale=1.2]{latex}}}}] (6,-3) -- (8,-2);
\draw[thin,postaction={decorate,decoration={markings,
    mark=at position .9 with {\arrow[scale=1.2]{latex}}}}] (8,-2) -- (10,0);
\draw[thin,postaction={decorate,decoration={markings,
    mark=at position .9 with {\arrow[scale=1.2]{latex}}}}] (10,0) -- (11,2);   
\draw[thin,postaction={decorate,decoration={markings,
    mark=at position .9 with {\arrow[scale=1.2]{latex}}}}] (11,2) -- (10,4);
\draw[thin,dashed] (10,4) -> (8,6);   

\begin{scope}[shift={(0.5,0.1)},scale=0.9,color=bleu]
\draw[thin,postaction={decorate,decoration={markings,
    mark=at position .9 with {\arrow[scale=1.2]{latex}}}}] (0,0) -> (2,-2);  
\draw[thin,postaction={decorate,decoration={markings,
    mark=at position .9 with {\arrow[scale=1.2]{latex}}}}] (2,-2) -- (4,-3); 
\draw[thin,postaction={decorate,decoration={markings,
    mark=at position .9 with {\arrow[scale=1.2]{latex}}}}] (4,-3) -- (6,-3); 
\draw[thin,postaction={decorate,decoration={markings,
    mark=at position .9 with {\arrow[scale=1.2]{latex}}}}] (6,-3) -- (8,-2);
\draw[thin,postaction={decorate,decoration={markings,
    mark=at position .9 with {\arrow[scale=1.2]{latex}}}}] (8,-2) -- (10,0);   
\draw[thin,postaction={decorate,decoration={markings,
    mark=at position .9 with {\arrow[scale=1.2]{latex}}}}] (10,0) -- (7,2);
\draw[thin,postaction={decorate,decoration={markings,
    mark=at position .9 with {\arrow[scale=1.2]{latex}}}}] (7,2) -- (3.5,2);     
\draw[thin,postaction={decorate,decoration={markings,
    mark=at position .9 with {\arrow[scale=1.2]{latex}}}}] (3.5,2) -- (0,0);  
    
\node[fill=white] at (2,-2) {};     
\node[fill=white] at (4,-3) {};      
\node[fill=white] at (6,-3) {};    
\node[fill=white] at (8,-2) {};  
\node[fill=white] at (10,0) {};
\node[fill=white] at (7,2) {}; 
\node[fill=white] at (3.5,2) {}; 
\node[fill=white] at (0,0) {};

\draw[red, very thick, ->] (8.75,1.75) -- (8.75,0.25); 

\draw[red, very thick, ->] (9,-2) -- (9,-0.5); 
\end{scope}

\begin{scope}[shift={(0.5,0.5)},scale=0.9,color=vert]
\draw[thin,dashed] (0,4) -> (2,6);  
\draw[thin,postaction={decorate,decoration={markings,
    mark=at position .9 with {\arrow[scale=1.2]{latex}}}}] (0,4) -> (-1,2);  
\draw[thin,postaction={decorate,decoration={markings,
    mark=at position .9 with {\arrow[scale=1.2]{latex}}}}] (-1,2) -- (0,0); 
\draw[thin,postaction={decorate,decoration={markings,
    mark=at position .9 with {\arrow[scale=1.2]{latex}}}}] (0,0) -- (3.5,2); 
\draw[thin,postaction={decorate,decoration={markings,
    mark=at position .9 with {\arrow[scale=1.2]{latex}}}}] (3.5,2) -- (7,2);    
\node[fill=white] at (0,4) {};
\node[fill=white] at (-1,2) {};
\node[fill=white] at (0,0) {};
\node[fill=white] at (3.5,2) {};
\node[fill=white] at (7,2) {};
\end{scope}

\begin{scope}[shift={(0.7,0.4)},scale=0.9,color=vert]  
\draw[thin,postaction={decorate,decoration={markings,
    mark=at position .9 with {\arrow[scale=1.2]{latex}}}}] (7,2) -- (10,0);     
\draw[thin,postaction={decorate,decoration={markings,
    mark=at position .9 with {\arrow[scale=1.2]{latex}}}}] (10,0) -- (11,2);  
\draw[thin,postaction={decorate,decoration={markings,
    mark=at position .9 with {\arrow[scale=1.2]{latex}}}}] (11,2) -- (10,4);  
\draw[thin,dashed] (10,4) -> (8,6);       
\node[fill=white] at (4.5,1) {};
\node[fill=white] at (7,2) {};
\node[fill=white] at (10,0) {};
\node[fill=white] at (11,2) {};   
\node[fill=white] at (10,4) {};  

\end{scope}

\node[fill=white] at (0,4) {};
\node[fill=white] at (-1,2) {};
\node[fill=white] at (0,0) {};
\node[fill=white] at (2,-2) {};
\node[fill=white] at (4,-3) {};
\node[fill=white] at (6,-3) {};
\node[fill=white] at (8,-2) {};
\node[fill=white] at (10,0) {};
\node[fill=white] at (11,2) {};
\node[fill=white] at (10,4) {};

\node[color=bleu] at (5,-2) {$\boldsymbol{\gamma_0}$};
\node[color=vert] at (2.5,4.5) {$\boldsymbol{\gamma'}$};
\node[] at (11,0.5) {$\boldsymbol{\gamma}$};

\node[color=rouge] at (7.5,-0.5) {$\boldsymbol{a_i}$};

\node[color=rouge] at (9,2) {$\boldsymbol{a_j}$};
\end{tikzpicture}
	\caption{$a_i\in\gamma_0\cap\gamma$ and $a_j~\notin\gamma_0\cap\gamma$}
	\label{fig:lemma5.4_case3}
	\end{subfigure}%
  	\caption{Possible cases for the induction. The cycle $\gamma'$ (in green) is obtained from the cycle $\gamma$ (in black) by deletion of one elementary cycle $\gamma_0$ (in blue).}
  	\label{fig:three_cases_lemma5.4}
  \end{figure}
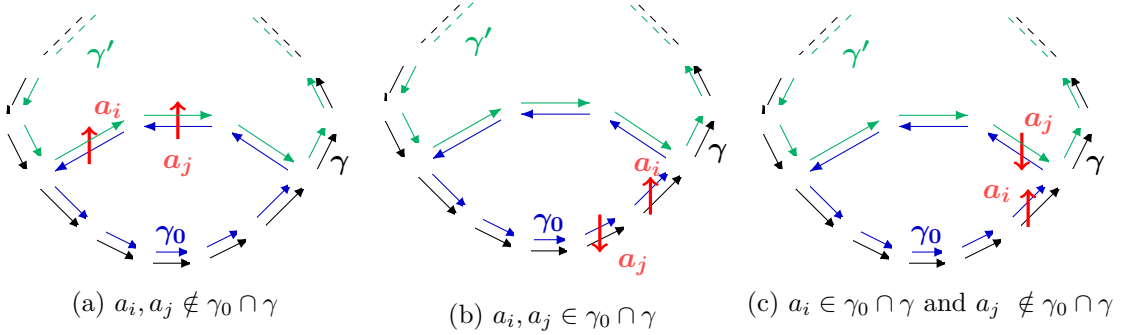

\begin{proof}
Since $a_k \circ \hdots \circ a_1 (g) = g$, the sequence of moves $a_1 \hdots a_k$ gives a cycle $\gamma$ starting from the vertex $g$ in the Cayley graph $\mC_G$ of $G$.
By abuse of notation, we will also call $a_i$ the labels of the edges in $\mC_G$ (thus from now we think of $a_i$ as an element of $S\cup S^{-1}$). So we are given a word $w=a_1\dots a_k\in (S\cup S^{-1})^k$ which represents the identity $1_G$. Since the word problem of the surface group of genus $2$ can be solved by Dehn's algorithm, this implies that we can obtain a finite sequence of words $w= w_0,w_1,w_2,\dots, w_N = 1_G$ such that $|w_i|>|w_{i+1}|$ and $w_{i+1}$ is obtained by $w_i$ by replacing the leftmost cyclical subword of $[a,b][c,d]$ of length at least $5$ by the inverse of its complement --for instance, the word $ba^{-1}b^{-1}cd$ can be rewritten as $ba^{-1}b^{-1}cd(c^{-1}d^{-1}a)(c^{-1}d^{-1}a)^{-1} = a^{-1}cd$ -- and then reducing the resulting word (eliminating pairs $ss^{-1}$ and $s^{-1}s$ for some generator $s$).

Because configurations in $x$ do not contain patterns in $\F_2$, the operation of reducing $w$ eliminates the same amount of up and down moves. Without loss of generality, we can replace $w$ by its reduced version. On the Cayley graph $\mC_G$, the operation of replacing a cyclical subword $u \sqsubset w$ by the inverse of its complement corresponds to decomposing the cycle $\gamma$ induced by $w$ into an elementary cycle $\gamma_0$ and the remaining cycle $\gamma'$. More precisely, if $w = w_1uw_2$, and $uv$ is an elementary cycle with $|v|<|u|$ then $w$ induces the cycle $\gamma$, $uv$ the elementary cycle $\gamma_0$ and $w_1v^{-1}w_2$ the remaining cycle $\gamma'$.

\begin{sloppypar}
We then prove the lemma by induction on the length of the chain \[{w= w_0,w_1,w_2,\dots, w_N = 1_G}.\]
In what follows, if $\zeta$ is a path in $\mC_G$ and $a_1\dots a_k$ its associated word on $S\cup S^{-1}$, we denote
\[{\upgen(\zeta):=\left|\left\{i\in\{1,\dots,k\}\mid a_i=\upgen_x\right\}\right|}\]
and
\[{\downgen(\zeta):=\left|\left\{j\in\{1,\dots,k\}\mid \exists k, a_j=\downgen_{k,x}\right\}\right|}.\]
\end{sloppypar}

If $N = 0$, then the reduced version of $w$ is the empty word. Hence $\upgen(\gamma)=\downgen(\gamma)$.

If $N\geq 1$, denote $w'=w_1v^{-1}w_2$ the word on $S\cup S^{-1}$ obtained after simplification by one cyclic permutation of $[a,b][c,d]$, $\gamma'$ the resulting cycle and $\gamma_0$ the elementary cycle corresponding to the simplification as explained above (see \cref{fig:three_cases_lemma5.4}).
Denote by $a_i$ (resp. $a_j$) the directed edge in $\gamma_0$ which is labeled by $\upgen_x$ (resp. $\downgen_{x,k}$ for some $k$) in configuration $x$. We distinguish between four cases, depending on where $a_i$ and $a_j$ are located. As no patterns from $\F_1$ appear in $x$, the elementary cycle $\gamma_0$ satisfies $\upgen(\gamma_0)=\downgen(\gamma_0)=1$, and by induction hypothesis, $\upgen(\gamma')=\downgen(\gamma')$.
Observe also that the directed edges $a_i$, $a_j$ of $\gamma_0$ are reversed if they also appear in $\gamma'$.
\begin{enumerate}
 \item If $\gamma_0\cap\gamma$ contains neither $a_i$ nor $a_j$ (see \cref{fig:lemma5.4_case1}). Then we have that
 $$\upgen(\gamma)=\upgen(\gamma')-\downgen(\gamma_0)=\downgen(\gamma')-\upgen(\gamma_0)=\downgen(\gamma).$$
 \item If $\gamma_0\cap\gamma$ contains $a_i$ and $a_j$ (see \cref{fig:lemma5.4_case2}). Then we have that
 $$\upgen(\gamma)=\upgen(\gamma')+\upgen(\gamma_0)=\downgen(\gamma')+\downgen(\gamma_0)=\downgen(\gamma).$$
 \item If $\gamma_0\cap\gamma$ contains $a_i$ but not $a_j$ (see \cref{fig:lemma5.4_case3}). Then we have that
 $$\upgen(\gamma)=\upgen(\gamma')-\downgen(\gamma_0)+\upgen(\gamma_0)=\upgen(\gamma')=\downgen(\gamma')=\downgen(\gamma).$$
 \item If $\gamma_0\cap\gamma$ contains $a_j$ but not $a_i$ (similar to case 3). Then we have that
 $$\upgen(\gamma)=\upgen(\gamma')=\downgen(\gamma')=\downgen(\gamma')+\downgen(\gamma_0)-\upgen(\gamma_0)=\downgen(\gamma).$$
\end{enumerate}
\end{proof}

Let us define $\rightgen_x^{-1}(g) = \leftgen_x(g)$, $\downgen_{1,x}^{-1}(g) = \upgen_x(g)$, and $1_G$ to be the identity of $G$.

\begin{lemma} \label{lemma:path}
  For any $g\in G$, there exists $i,j$ such that $g = \rightgen_x^j\circ\downgen_{1,x}^i(1_G)$.
\end{lemma}

\begin{proof}
  Let $x\in X$ and $g\in G$. As each symbol contains all $8$ directions in $S$, it is clear that there exist $a_1, \hdots, a_k \in \{\leftgen_x, \rightgen_x, \upgen_x, \downgen_{1,x}, \downgen_{2,x}, \hdots \}$ such that $g = a_k \circ\hdots\circ a_1 (1_G)$.

  First, we can get rid of all $\downgen$~that are not $\downgen_1$, indeed for any $l$, $\downgen_{l,x} = \rightgen_x^{6(l-1)}\circ\downgen_{1,x}$ (see \cref{fig:op-down}).
  So, by transforming all $\downgen$~ like this, we obtain $i_1, \hdots, i_l\in\Z$ such that $g = \downgen_{1,x}^{i_l}\circ\rightgen_x^{i_{l-1}} \hdots \circ \downgen_{1,x}^{i_2}\circ\rightgen_x^{i_1}(g)$.

  Let us concentrate on $\downgen_{1,x}^n \circ \rightgen_x^m(h)$ for some $m,n\in\Z$ and $h\in G$.
  Let $w$ be the word of size $m$ such that $w_i = x_{\rightgen_x^i(h)}$ for $i\in\{1\hdots m\}$.
  Then, as shown of \cref{fig:op-order}
  , $\downgen_{1,x}^n \circ \rightgen_x^m(h) = \rightgen_x^{|\texttt{s}^m(w)|} \circ \downgen_{1,x}^n(h)$.
  By doing this operation on all badly ordered operations in the sequence, and obtain $i$ and $j$ such that $g = \rightgen_x^j(\downgen_{1,x}^i(1_G))$.

  \begin{figure}[htp]
    \centering
    \scalebox{0.9}{\begin{tikzpicture}

  \foreach \i in {1, 3, 4}{
    \draw[->, >=stealth, thick] (9.7/2+.25, 3+.25) -- (0+\i*2+.25,0+.5) node [pos=0.5, above] {} ;
  }

  \foreach \i in {0}{
    \draw[->, >=stealth, ultra thick] (9.7/2+.25, 3+.25) -- (0+\i*2+.25,0+.5) node [pos=0.5, above] {$\downgen_{1,x}$} ;
  }

  \foreach \i in {2}{
    \draw[->, >=stealth, ultra thick, color=red] (9.7/2+.25, 3+.25) -- (0+\i*2+.25,0+.5) node [pos=0.5, right] {$\downgen_{l,x}$} ;
  }

  \draw[->, >=stealth, ultra thick] (0.25,-.25) -- (4.25, -.25) node [pos=0.5, below] {$\rightgen_x^{6(l-1)}$} ;


\foreach \i in {0, 1, ..., 3}{
	\draw[fill=black] (0+\i*2,0) rectangle (0.5+\i*2, 0.5) ;
  \foreach \j in {0, 1, ..., 4}{
  	\draw[fill=white] (0.3*\j+0.5+\i*2,0.1) rectangle (0.3*\j+0.3+0.5+\i*2, 0.5-0.1) ;
  }
}
\foreach \i in {4}{
	\draw[fill=black] (0+\i*2,0) rectangle (0.5+\i*2, 0.5) ;
  \foreach \j in {0, 1, ..., 3}{
  	\draw[fill=white] (0.3*\j+0.5+\i*2,0.1) rectangle (0.3*\j+0.3+0.5+\i*2, 0.5-0.1) ;
  }
}

\draw[fill=black] (9.7/2, 3) rectangle (9.7/2+.5, 3+.5) ;

\end{tikzpicture}}
    \caption{Transformation to get only $\downgen_1$ down operations.}
    \label{fig:op-down}
  \end{figure}
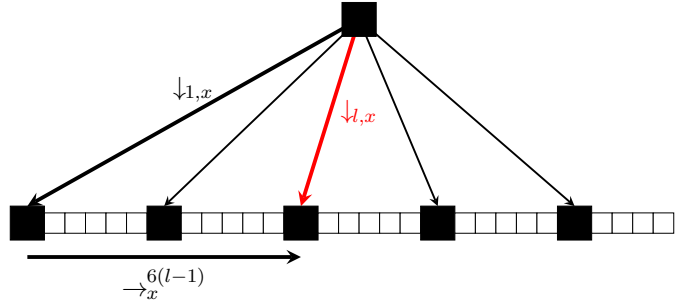

  \begin{figure}[htp]
    \centering
    \scalebox{0.8}{\begin{tikzpicture}


\node[scale=1.25] at (0,3) {\textbullet};
\node[scale=1.25] at (-0.3,3.3) {$1_G$};



\node[scale=1.25] at (4.8,0) {\textbullet};
\node[scale=1.25] at (4.8+.3,-.3) {$h$};

\draw[decorate,decoration={brace,amplitude=6pt,raise=0.1cm}] (0,3.25) -- (4,3.25);
\node[scale=1.25] at (2, 4) {$|w|$};

\draw[decorate,decoration={brace,amplitude=6pt,raise=0.1cm}] (-0.9-0.25, 0) -- (-0.25,3);
\node[scale=1.25] at (-1.3, 1.65) {$n$};

\draw[decorate,decoration={brace,amplitude=6pt,raise=0.1cm}] (4.8,-0.25) -- (-0.9,-0.25);
\node[scale=1.25] at (1.95, -0.9) {$|s(w)|$};


\foreach \i in {1, 2, ..., 4}{
  \draw[->, >=stealth, thick] (\i-1,3) -- (\i,3) node [] {} ;
}

\foreach \i in {1, 2, 3}{
  \draw[->, >=stealth, thick] (0-0.3*\i+0.3,3+1-\i) -- (0-0.3*\i,3-\i) node [] {} ;
}

\foreach \i in {1, 2, ..., 19}{
  \draw[->, >=stealth, thick] (\i*0.3-0.3-0.9,0) -- (\i*0.3-0.9,0) node [] {} ;
}

\foreach \i in {1, 2, 3}{
  \draw[->, >=stealth, thick] (4+0+0.8/3*\i-0.8/3,3+1-\i) -- (4+0.8/3*\i,3-\i) node [] {} ;
}

\end{tikzpicture}}
    \caption{Transformation to reorder the operations.}
    \label{fig:op-order}
  \end{figure}
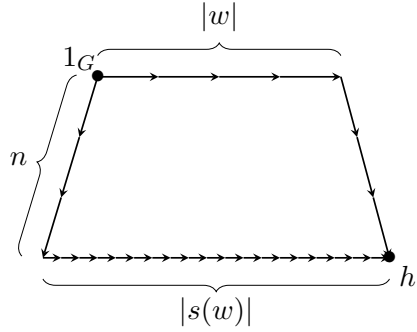

\end{proof}

\subsubsection{A bijection between \texorpdfstring{$\Z^2$}{Z2} and the surface group}
Let $x\in X$ be fixed. We define $f_x\colon \Z^2 \rightarrow G$ to be the following:
\[ f_x (i,j) = \rightgen_x^j\circ \downgen_{1,x}^i(1_G) .\]

\begin{lemma}
  \label{lemma:bijection1}
  For every $x\in X$, the function $f_x$ is a bijection.
\end{lemma}

\begin{proof}
  First, $f$ is well-defined because the operations $\rightgen_x(g)$ and $\downgen_{1,x}(g)$ are both well-defined for any $g\in G$.
  The existence of $i,j\in\Z$ such that $g = \rightgen_x^j\circ\downgen_{1,x}^i(1_G)$ is ensured by \cref{lemma:path}.

  For the uniqueness of such $i,j$, let us assume there are $i',j' \in\Z$ such that $ g = \rightgen_x^{j'}\circ\downgen_{1,x}^{i'}(1_G)$. Since $g^{-1}\cdot g = 1_G$, we get
  \[
  \downgen_{1,x}^{-i'}\circ\rightgen_x^{-j'}(g) = \downgen_{1,x}^{-i'}\circ\rightgen_x^{j-j'}\circ\downgen_{1,x}^i(1_G) = 1_G .\]
  \cref{lemma:cycles} ensures that $i=i'$.
  Then, because we only consider $\downgen_{1,x}$ operations (the first bottom neighbor and not the others), their inverses are $\upgen_x$ operations. It means that the only way of having a cycle is to have $\downgen_{1,x}^{-i'}(1_G) = \downgen_{1,x}^{-i'}\circ\rightgen_x^{j-j'}(1_G)$. Thus we have a cycle using only right operations (or only left operations), \cref{lemma:ring} ensures that $j = j'$ since the only way of having a cycle with only right (or only left) operations is to not apply any.
\end{proof}

We can moreover prove that $f_x$ also preserves locality.

\begin{lemma}
  \label{lemma:locality}
  The following equivalences are true:
  \begin{enumerate}
    \item \label{equi:next}
    $\begin{cases}
      (u,v)\in E_\Gamma\\
      L_\Gamma(u,v) = \texttt{next}
    \end{cases}
    \Leftrightarrow
    f_x(v) = \rightgen_x(f_x(u))
    $

    \item \label{equi:other}
    $\begin{cases}
      (u,v)\in E_\Gamma\\
      L_\Gamma(u,v) = k \in \{0,\hdots,M-1\}
    \end{cases}
    \Leftrightarrow
    f_x(v)=\rightgen_x^{k}\circ\downgen_{1,x}(f_x(u))
    $
    \\  where $M$ is the number of sons of $u$.
  \end{enumerate}
\end{lemma}

\begin{proof}
  \begin{enumerate}
    \item  If $L_\Gamma(u,v) = \texttt{next}$, then $(u,v) = ((i,j), (i,j+1))$, and so $f_x(v) = \rightgen_x^{j+1}\circ \downgen_{1,x}^i(1_G) = \rightgen_x(f_x(u))$.
    \\
    Conversely, assume $f_x(v) = \rightgen_x(f_x(u))$. Consider $i,j$ such that $f_x(u) = \rightgen_x^{j}\circ \downgen_{1,x}^i(1_G)$. Then $f_x(v) = \rightgen_x^{j+1}\circ \downgen_{1,x}^i(1_G)$, implying that $(u,v) = ((i,j), (i,j+1))$ by definition of $f_x$. Then, we can only have $L_\Gamma(u,v) = \texttt{next}$.

    \item \label{equi:other:proof}
    Assume that $L_\Gamma(u,v) = k$, we know that $(u,v) = ((i,j), (i+1, \Delta_{i+1}(j)+k))$ and so $f_x(v) = \rightgen_x^{\Delta_{i+1}(j)+k}\circ \downgen_{1,x}^{i+1}(1_G) = \rightgen_x^k\circ\downgen_{1,x}\circ \rightgen_x^{j}\circ \downgen_{1,x}^i(1_G) = \rightgen_x^k\circ\downgen_{1,x}(f_x(u))$ by definition of $\Delta_{i+1}(j)$.\\
    Conversely suppose $f_x(v) = \rightgen_x^k\circ\downgen_{0,x}(f_x(u))$.
    Assume also that $f_x(u) = \rightgen_x^{j}\circ \downgen_{1,x}^i(1_G)$, then $f_x(v) = \rightgen_x^k\circ\downgen_{1,x}\circ \rightgen_x^{j}\circ \downgen_{1,x}^i(1_G) = \rightgen_x^k\circ\rightgen_x^{\Delta_{i+1}(j)}\circ \downgen_{1,x}^{i+1}(1_G)=\rightgen_x^{\Delta_{i+1}(j)+k}\circ \downgen_{1,x}^{i+1}(1_G)$.
    So we get $(u,v) = ((i,j), (i+1, \Delta_{i+1}(j)+k))$ and $L_\Gamma(u,v) = k$.
  \end{enumerate}
\end{proof}

As said before, the bijection $f_x$ itself cannot be a label preserving graph isomorphism, since we lack some edges in $\mC_G$, but it nevertheless enjoys a useful property: if $\varphi$ is a label preserving graph isomorphism for $\Gamma$, then so is $f_x\circ\varphi\circ f_x^{-1}$ for $\mC_{G,x}$, and if $\varphi$ is a label preserving graph isomorphism for $\mC_{G,x}$, then so is $f_x^{-1}\circ\varphi\circ f_x$ for $\Gamma$, where $\mC_{G,x}$ is a relabeling of $\mC_g$ according to the configuration $x$. So roughly speaking, any local pattern is preserved by $f_x$ or by $f_x^{-1}$ (see~\cref{coro:bijection2} below).

\begin{coro}
  \label{coro:bijection2}
  Let $\A$ be a finite alphabet.
  For any configuration $c\in\A^{G}$, $p \sqsubset c \Rightarrow f_x^{-1}(p) \sqsubset f_x^{-1}(c)$.
  Conversely for any $d\in \A^\Gamma$, $q \sqsubset d \Rightarrow f_x(q) \sqsubset f_x(d)$.
\end{coro}

\begin{proof}
Define $\mC_{G,x}$ the oriented labeled graph obtained from $\mC_{G}$ by replacing every label in $\mC_{G}$ by the corresponding symbol in $\{\leftarrow,\rightarrow,\uparrow,\downarrow_1,\downarrow_2,\downarrow_3,\downarrow_4,\downarrow_5,\downarrow_6\}$ found in the configuration $x$: if $g_i\in S\cup S^{-1}$ labels $(g,g\cdot g_i)$ in $\mC_G$ and $(g_i,\star)\sqsubset x_g$, then $\star$ labels $(g,g\cdot g_i)$ in $\mC_{G,x}$.
Since there is no ambiguity here, let us write $f=f_x$ in this proof for lighter notations.

\medskip

We first prove that if $\varphi$ is a label preserving graph isomorphism for $\mC_{G,x}$, then so is $\psi=f^{-1}\circ\varphi\circ f$ for $\Gamma$.
Obviously, $\psi$ is a bijection as composition of bijections. Take some edge $(u,v)\in E_\Gamma$, then
  \begin{align*}
    L_\Gamma(u,v)=\texttt{next} & \Leftrightarrow f(v) = \rightgen_x(f(u)) & \text{(by \cref{lemma:locality})}  \\
     & \Leftrightarrow L_{\mC_{G,x}}\left(f(u),f(v)\right)=\rightgen_x \\
     & \Leftrightarrow L_{\mC_{G,x}}\left(\varphi\circ f(u),\varphi\circ f(v)\right)=\rightgen_x & \text{($\varphi$ is label-preserving)} \\
     & \Leftrightarrow L_{\Gamma}\left(f^{-1}\circ\varphi\circ f(u),f^{-1}\circ\varphi\circ f(v)\right)=\texttt{next}  & \text{(by \cref{lemma:locality})}\\
     & \Leftrightarrow L_{\Gamma}\left(\psi(u),\psi(v)\right)=\texttt{next}
  \end{align*}
and for $k\in\{0,\hdots,M-1\}$,
  \begin{align*}
    L_\Gamma(u,v)=k & \Leftrightarrow f(v)=\rightgen_x^{k}\circ\downgen_{1,x}(f(u)) & \text{(by \cref{lemma:locality})}  \\
     & \Leftrightarrow \varphi\circ f(v)=\rightgen_x^{k}\circ\downgen_{1,x} \left(\varphi\circ (f(u))\right) & \text{($\varphi$ is label-preserving)} \\
     & \Leftrightarrow L_\Gamma(f^{-1}\circ\varphi\circ f(u),f^{-1}\circ\varphi\circ f(v))=k & \text{(by \cref{lemma:locality})}\\
     & \Leftrightarrow L_{\Gamma}\left(\psi(u),\psi(v)\right)=k.
  \end{align*}

Assume now that $p\in A^S$, with $S$ a finite subset of $G$, is a pattern that appears in a configuration $c\in\A^{\mC_G}$. By definition, there exists $\varphi\colon S\to T$ a label preserving graph isomorphism for $\mC_G$ such that
$$\varphi(p)=c_|{T}.$$
Define $\psi:=f^{-1}\circ\varphi\circ f$. By what precedes, $\psi$ is also a label preserving graph isomorphism for $\Gamma$, and
  \begin{align*}
    \psi\left(f^{-1}(p)\right) & = f^{-1}\circ \varphi(p) \\
     & = f^{-1}(c|_{T})\\
     & =f^{-1}(c)|_{f^{-1}(T)}.
  \end{align*}
So the pattern $f^{-1}(p)$ appears in the configuration $f^{-1}(c)$.

\bigskip

Conversely, let us prove that if $\varphi$ is a label preserving graph isomorphism for $\Gamma$, then so is $\psi = f\circ\varphi\circ f^{-1}$ for $\mC_{G,x}$.
Take some edge $(u,v)\in G$, then
  \begin{align*}
    L_{\mC_{G,x}}\left(u,v\right)=\rightgen_x & \Leftrightarrow L_\Gamma(f^{-1}(u),f^{-1}(v))=\texttt{next} & \text{(by \cref{lemma:locality})}  \\
     & \Leftrightarrow L_\Gamma(\varphi\circ f^{-1}(u),\varphi\circ f^{-1}(v))=\texttt{next} & \text{($\varphi$ is label-preserving)} \\
     & \Leftrightarrow L_{\mC_{G,x}}\left(f\circ\varphi\circ f^{-1}(u),f\circ\varphi\circ f^{-1}(u)\right)=\rightgen_x  & \text{(by \cref{lemma:locality})}\\
     & \Leftrightarrow L_{\mC_{G,x}}\left(\psi(u),\psi(v)\right)=\rightgen_x.
  \end{align*}
For $l\in\{0,\hdots, 7\}$, if $v$ corresponds to the $\downgen_{l,x}$ neighbor of $u$ in $\mC_G$, it corresponds to the $6(l-1)$-th child of $f^{-1}(u)$ in $\Gamma$ and vice-versa (see \cref{fig:op-down}). Therefore,
  \begin{align*}
    L_{\mC_{G,x}}\left(u,v\right)=\downgen_{l,x} & \Leftrightarrow L_\Gamma(f^{-1}(u),f^{-1}(v))=6(l-1) & \text{(by \cref{lemma:locality})}  \\
     & \Leftrightarrow L_\Gamma(\varphi\circ f^{-1}(u),\varphi\circ f^{-1}(v))=6(l-1) & \text{($\phi$ is label-preserving)} \\
     & \Leftrightarrow L_{\mC_{G,x}}\left(f\circ\varphi\circ f^{-1}(u), f\circ\varphi\circ f^{-1}(v)\right)=\downgen_{l,x} & \text{(by \cref{lemma:locality})}\\
     & \Leftrightarrow L_{\mC_{G,x}}\left(\psi(u),\psi(v)\right)=\downgen_{l,x}.
  \end{align*}

As previously, if $q\in A^S$, with $S$ a finite subset of $\Z^2$, is a pattern that appears in a configuration $d\in\A^\Gamma$. By definition, there exists $\varphi\colon S\to T$ a label preserving graph isomorphism for $\Gamma$ such that
$$\varphi(q)=d|_{T}.$$
With $\psi:=f\circ\varphi\circ f^{-1}$, we have
  \begin{align*}
    \psi\left(f(q)\right) & = f\circ \varphi(q) \\
     & = f(d|_{T})\\
     & =f(d)|_{f(T)}.
  \end{align*}
So the pattern $f(q)$ appears in the configuration $f(d)$.
\end{proof}

\subsection{The reduction}\label{subsection.reduction_genus_2}
We now have everything in hand to prove the undecidability of the domino problem on the surface group of genus~$2$.

\begin{theorem}\label{theorem:DP_genus2}
  The domino problem is undecidable on the surface group of genus 2.
\end{theorem}

\begin{proof}
  Recall that $\Gamma$ is the orbit graph of an orbit $\Omega$ of the substitution $s$ defined on page~\pageref{definition.substitution}. Let $\A$ be a finite alphabet and $Y\subseteq\A^\Gamma$ an SFT over $\Gamma$, given by a finite set of forbidden patterns $\F_Y$. We define $Z$ the SFT over $G$ with set of forbidden patterns $\F_Z := f_x(\F_Y)$, where $f_x$ is defined in \cref{lemma:bijection1}. Clearly $\F_Z$ can be constructed effectively from $\F_Y$. We show that $Z=\emptyset$ if and only if $Y=\emptyset$, providing a reduction to $\texttt{DP}(\Gamma)$, that comes from $\texttt{s}$ which han expanding eigenvalue. By~\cref{th:DP_orbit} $\texttt{DP}(\Gamma)$ is undecidable.

  Assume $Z=\emptyset$ and consider a configuration $c\in\A^{G}$. The configuration $d:=f_x^{-1}(c)$ is thus in $\A^\Gamma$. Since $Z=\emptyset$, necessarily $c$ contains a forbidden pattern $p$ from the set $\F_Z$. Since $p\sqsubset c$, \cref{coro:bijection2} implies that $f_x^{-1}(p) \sqsubset f_x^{-1}(c) = d$.
  So a pattern $f_x^{-1}(p)$ from $\F_Y$ appears in any configuration $c\in \A^G$, i.e. the subshift $Y$ is empty.

  Conversely, if $Y=\emptyset$, take any $d\in\A^\Gamma$ and $c:=f_x(d)\in\A^G$.
  Because $Y=\emptyset$, $d$ contains a forbidden pattern $q\in\F_Y$.
  Since $q\sqsubset d$, \cref{coro:bijection2} implies that $f_x(q) \sqsubset f_x(d) = c$. Therefore, the pattern $f_x(q)\in \F_Z$ appears in any $d\in\A^\Gamma$, so $Y=\emptyset$ as well and the equivalence is proved.
\end{proof}

\begin{coro}
  The domino problem is undecidable for every surface group.
\end{coro}

\begin{proof}
 The undecidability of the domino problem is a commensurability invariant (see~\cite[Corollary 9.53]{BertheRigo2018}), and all surface groups of genus $g\geq2$ are commensurable (see~\cite[Proposition 6.7]{ClimenhagaKatok} for a recent reference). By combining these two facts with~\cref{theorem:DP_genus2}, we obtain the undecidability of domino problem for surface groups of any genus $g\geq2$. As the domino problem on $\Z^2$ --the surface group of genus $1$-- is undecidable, we obtain our result.
\end{proof}

\section{Remarks about word-hyperbolic groups}
\label{section.remarks}

The domino problem conjecture states the following.

\begin{conjecture}\label{conj_DP}
	A finitely generated group has decidable domino problem if and only if it is virtually free.
\end{conjecture}

Surface groups of genus $g \geq 2$ are special cases of a larger class of groups called \emph{word-hyperbolic} for which geodesic triangles are hyperbolic. They can be characterized as the finitely presented groups for which Dehn's algorithm solves the word problem.

An important property of the domino problem is that groups which contain subgroups with undecidable domino problem have themselves undecidable domino problem (see~\cite[Proposition 9.3.30]{BertheRigo2018}). This means that every group which contains an embedded copy of a surface group has undecidable domino problem. This is of special relevance due to the following conjecture by Gromov.

\begin{conjecture}[Gromov]
	Every one-ended word-hyperbolic group contains an embedded copy of the surface group of genus $2$.
\end{conjecture}

In particular, if Gromov's conjecture holds, every one-ended word-hyperbolic group would automatically have undecidable domino problem. A group can have either $0,1,2$ or infinitely many ends. In the case when it has $0$ ends it is finite and thus its domino problem is trivially decidable, and whenever it has $2$ ends it is virtually $\Z$ and thus it is also decidable. In the case of a finitely presented group $G$, a fundamental result by Dunwoody~\cite{Dunwoody1985} shows that if $G$ has infinitely many ends it can be expressed as the fundamental group of a finite graph of groups such that every edge is a finite group and all vertices are either finite or $1$-ended. It can also be shown (see~\cite{Karrass1973}) that $G$ is virtually free if and only if all of the vertex groups in its decomposition are finite. Therefore, if $G$ is not virtually free, it must contain a one-ended subgroup. In the case of word-hyperbolic groups, every such group in the decomposition must also be word-hyperbolic (see~\cite{
Bowditch1998}). In other words, every word-hyperbolic group which is not virtually free contains a one-ended word-hyperbolic group. This implies the following.

\begin{proposition}
	If Gromov's conjecture holds then the domino problem conjecture holds for all word-hyperbolic groups.
\end{proposition}

In fact, we could obtain the same result with a weaker version of Gromov's conjecture. We say a group $G$ acts \emph{translation-like} on a metric space $(X,d)$ if the action is free and $\sup_{x \in X}{d(x,gx)}<\infty$ for every $g \in G$. Clearly, if $H$ is a subgroup of $G$ then $H$ acts translation like on any Cayley graph of $G$ by multiplication. A theorem by Jeandel~\cite{Jeandel2015translation} shows that if a finitely presented group $H$ acts translation like on a Cayley graph of a finitely generated group $G$, then the domino problem of $H$ is many-one reducible to the domino problem on $G$, in particular, we obtain that any group on which the surface group of genus $2$ acts translation-like has undecidable domino problem.

\begin{proposition}
	If every $1$-ended word-hyperbolic group admits a translation-like action of the surface group of genus $2$, then the domino problem conjecture holds for all word-hyperbolic groups.
\end{proposition}

\section*{Acknowledgments}
The authors would like to thank Chaim Goodman-Strauss for fruitful discussions which pointed us in the right direction in~\cref{section.undecidability_orbit_graphs}. We would also grateful to Yann Ollivier for his pictures of the Cayley graph of the surface group~\cite{Ollivier}, the one in the front page being one of them. They were the starting point of many ideas that led to this paper. This project was partially supported by the ANR project CoCoGro (ANR-16-CE40-0005).

\bibliographystyle{abbrv}
\bibliography{biblio}

\end{document}